\documentclass[sn-mathphys,Numbered]{sn-jnl}

\usepackage{graphicx}%
\usepackage{amsmath,amssymb,amsfonts,bm}%
\usepackage{epstopdf}
\usepackage{multirow}
\usepackage{tikz}
\usepackage{amsthm}%
\usepackage{float}
\usepackage{tabularx}
\usepackage{hyperref}
\usepackage{subcaption}
\usepackage[title]{appendix}%
\usepackage{textcomp}%
\usepackage{manyfoot}%
\usepackage{algpseudocode}%
\usepackage{listings}%
\usepackage{bigints}
\usepackage{outlines}
\usepackage{geometry}
\usepackage{siunitx}
\usepackage{booktabs}
\usepackage{setspace}
\usepackage{lineno}

\usepackage[utf8]{inputenc}
\usepackage{graphicx}
\usepackage[ruled,vlined]{algorithm2e}
\usepackage{amsmath}
\usepackage[section]{placeins}
\usepackage{enumitem}
\usepackage{color,soul}
\usepackage{xfrac}

\begin{document}

\title[Article Title]{A finite volume Simo-Reissner beam method for moored floating body dynamics}

\author[1]{\fnm{Amirhossein} \sur{Taran}}
\author[1]{\fnm{Seevani} \sur{Bali}}
\author[2]{\fnm{\v{Z}eljko} \sur{Tukovi\'{c}}}
\author[1]{\fnm{Vikram} \sur{{Pakrashi}}}
\author*[1]{\fnm{Philip} \sur{Cardiff}}\email{philip.cardiff@ucd.ie}

\affil[1]{\orgdiv{School of Mechanical and Materials Engineering}, \orgname{University College Dublin}, \orgaddress{\country{Ireland}}}
\affil[2]{\orgdiv{Faculty of Mechanical Engineering and Naval Architecture, University of Zagreb}, \orgaddress{\country{Croatia}}}

\abstract
{
This paper presents a novel finite volume mooring line model based on the geometrically exact Simo-Reissner beam model for analysing the interaction between a floating rigid body and its mooring lines. The coupled numerical model is implemented entirely within a finite volume-based discretisation framework using a popular computational fluid dynamics C++ toolbox, OpenFOAM. Unlike existing methods for modelling mooring lines, which rely on lumped mass models or finite element-based approaches, this work simulates the mooring cables using non-linear beam models implemented in a finite volume framework to account for bending, tensile, and torsional loading. This advancement makes the current work particularly valuable for simulating extreme sea conditions. The coupled model developed in this study has been validated and verified using experimental and numerical data for a floating box moored with four catenary mooring lines under regular wave conditions featuring different wave heights and periods. The results demonstrate strong agreement with both experimental and numerical data, highlighting the model’s accuracy in capturing mooring dynamics and floating body motion.
}


\maketitle

\section{Introduction}\label{sec:intro}
The importance of mooring lines extends beyond mere attachment and station-keeping; they fundamentally contribute to the safety and operational viability of floating offshore structures, especially in environments where uncontrolled movement could lead to catastrophic consequences, such as complete structural failure \citep{yang_investigation_2021}.
An accurate model of mooring lines is not only essential for understanding the hydrodynamic response amplitude of floating offshore structures but is also critical for assessing long-term fatigue and ensuring the longevity of both the mooring system and the floating structure it supports.

\noindent Mooring line analysis models are generally categorised into static, quasi-static, and dynamic types. Static models \citep{oppenheim_static_1982,smith_statics_2001} assume that the system is in static equilibrium, accounting only for constant forces such as steady wind or current loads while neglecting any time-dependent or dynamic effects. This simplification allows for basic assessments of mooring line tension and position \citep{karnoski_validation_1988} under calm or steady-state conditions but may overlook critical forces in more variable environments. Quasi-static models, on the other hand, enhance the realism of static analyses by incorporating slowly varying environmental forces and accounting for the mooring line's stiffness response to changing loads \citep{dunbar_development_2015,cheng_numerical_2019}. However, they still neglect high-frequency dynamics associated with wave action. The quasi-static approach offers a balance between simplicity and accuracy \citep{liu_establishing_2017}. It is well-suited for conditions dominated by low-frequency forces. In contrast, full dynamic modelling is more appropriate for scenarios where wave-induced motions play a significant role. When the floating structures experience large amplitude motions, the dynamic behaviours of the
mooring system become non-negligible \citep{masciola_implementation_2013}. Thus, quasi-static models will underestimate the tension on the lines \citep{chen_cfd_2022,hall_validation_2015}.
Incorporating dynamic effects into the equations of motion for mooring lines 
allows dynamic analysis models to predict mooring loads and the resulting responses of 
floating structures with greater accuracy compared to static and quasi-static models.
Recent dynamic mooring models include Mooring Design and Dynamics \citep{dewey_mooring_1999}, which addresses dynamic equilibrium equations under environmental loads; SEAWAY \citep{jmj_journee_ljm_adegeest_theoretical_2003}, a frequency-domain model; and OrcaFlex \citep{randolph_non-linear_2010}, a non-linear framework for time-domain analysis.
Ferri and Palm \citep{ferri_implementation_2015} developed a finite element model (named Moody) with high-order polynomial basis functions and spatial discretisation using the local discontinuous Galerkin method. Following this higher-order approach, Moody was developed with the aim of resolving snap loads, and it incorporates a range of hydrodynamic forces within the model \citep{davidson_mathematical_2017}. The bending stiffness has later been added to the model to investigate the bending stiffness on snap-loads \citep{palm_influence_2020}. 
Hall et al. \citep{hall_moordyn_2015} developed an open-source dynamic model called MoorDyn based on a lumped mass approach \citep{kreuzer_mooring_2002}, which discretises the mooring line into a finite number of mass nodes connected with mass-less springs. 
While the initial version of MoorDyn neglects bending stiffness, this feature was later added in the model \citep{hall_implementation_2021}.
In all the aforementioned references, the consistent conclusion is that incorporating bending stiffness into mooring cable analysis diversifies the usability of the numerical simulation to low-tension and slack conditions, avoiding any ill-conditioning of the system characteristics.

\noindent Recently, some researchers have focused on coupling dynamic mooring analysis models with computational fluid dynamics (CFD) simulations to enhance the accuracy of mooring system behaviour predictions under complex fluid-structure interactions. This approach allows for the detailed simulation of forces from waves and currents on mooring lines, improving the understanding of both static and dynamic responses in challenging offshore environments \citep{darling_role_2024}. CFD models can utilise either mesh-based or mesh-free approaches. Mesh-free approaches are primarily suitable for free-surface wave structure interaction problems due to their speed-up when using graphical computing units and their ability to capture interfaces without any additional processing \citep{chen_cfd_2022,lind_review_2020}. Dominguez et al. \citep{dominguez_sph_2019} was the first to couple DualSPHysics \citep{crespo_dualsphysics_2015}, which is an open-source parallel mesh-free CFD solver based on smoothed particle hydrodynamics (SPH) with MoorDyn. Later, Liu and Wang \citep{liu_numerical_2020} extended this coupling to study different submerged box-type breakwaters. Although these mesh-free methods are popular among researchers, they can become highly complex regarding particle refinement and turbulence modelling \citep{macia_truncated_2021,vacondio_grand_2021}. 
Islam et al. \citep{islam_openfoam_2019} studied the wave radiation due to heave, surge, and pitch by a box-type floating structure using OpenFOAM \citep{weller_tensorial_1998}, a widely-used CFD toolbox. 
Haider et al. \citep{haider_comprehensive_2024} carried out a fully coupled fluid-structure interaction analysis for a floating offshore wind turbine, employing MoorDyn to model the mooring system and OpenFOAM as the CFD platform.
Wei et al. \citep{wei_coupled_2024} performed a coupled simulation of a very large floating structure (VLFS) moored with catenary mooring lines using OpenFOAM and MBDyn \citep{masarati_efficient_2014}. 
Palm et al. \citep{palm_coupled_2016} conducted a coupled CFD-mooring analysis of a wave energy converter (WEC), coupling OpenFOAM with Moody to simulate the behaviour of the mooring lines. Similarly, Chen and Hall \citep{chen_cfd_2022} integrated MoorDyn with OpenFOAM to simulate a moored floating box. More recently, they extended their work to model multiple floating objects \citep{CHEN2024117697}. Beyond the studies mentioned, additional research has explored coupling mooring line models with various CFD software \citep{qiao_dynamic_2024,huang_new_2022,martin_numerical_2021,chen_coupled_2018}. However, OpenFOAM's customizability makes it a preferred choice among researchers as a CFD solver. Furthermore, the reviewed studies highlight OpenFOAM's capability to accurately simulate the hydrodynamic response of floating rigid bodies.

\noindent The literature on moorings highlights the importance of incorporating bending stiffness into numerical simulations for accurate analysis. For instance, commercial software such as OrcaFlex \citep{orcina_orcaflexdocumentation_2020}, which employs a lumped mass approach, explicitly calculates bending moments based on the bending angle at the intersections of the lumped nodes in the mooring line model. Rather than using lumped mass or finite element-based approaches, one can consider the mooring line as a beam or rod model undergoing deformations due to its slender shape and flexibility \citep{cottanceau_finite_2018}.
In the existing literature, beam models range from linear formulations like the Euler-Bernoulli model, which neglects shear deformation, to the geometrically exact Simo-Reissner beam model. The latter are non-linear and account for large deformations and rotations, offering a more comprehensive representation of beam behaviour \citep{reissner_one-dimensional_1972,reissner_finite_1981,simo_finite_1985}. While numerous finite element formulations \citep{simo_three-dimensional_1986,zupan_quaternion-based_2009,jelenic_geometrically_1999,ibrahimbegovic_computational_1995} have been introduced in the literature using geometrically exact Simo-Reissner beam theory, Tukovi\'{c} et al. \citep{tukovic_finite_2019} were the first to propose a finite volume formulation of geometrically exact Simo-Reissner beams. More recently, Bali et al. \citep{bali_cell-centered_2022} expanded the finite volume formulation for shear-deformable geometrically exact beams and conducted a comparative analysis with finite element-based classical benchmark cases using OpenFOAM. Additionally, they extended this formulation to handle beam-to-beam contact interactions \citep{bali_finite_2024}. The primary motivation for adapting the Simo–Reissner formulation to finite volume discretisation within OpenFOAM is to unify numerical methods for solids and fluids, allowing complex multi-physics problems, like fluid-solid interactions, to be addressed within a unified framework. Furthermore, \citet{bali_cell-centered_2022} noted that, unlike the finite element approach, the finite volume formulation for beams is not subject to the shear-locking phenomenon.

\noindent In the current work, the finite volume Simo–Reissner beam formulation introduced by \citet{bali_cell-centered_2022} is extended by fully integrating it with a multiphase CFD solver and a six-degree-of-freedom motion solver. This results in a unified, fully coupled fluid–structure interaction framework within the finite volume method. The proposed framework is capable of capturing large nonlinear deformations in mooring lines—including bending, shear, and torsion—while simultaneously resolving the free surface dynamics of the fluid. While many existing approaches simplify mooring dynamics using lumped-mass or quasi-static models, the present work introduces a unified, comprehensive simulation environment for moored floating structures.
A key novelty of this study lies in the application of a geometrically exact dynamic beam theory in a finite volume context, which allows for consistent coupling with the CFD solver and effective simulation of nonlinear mooring behaviour. In addition, the model incorporates seabed contact and friction effects to capture realistic boundary interactions. The coupled framework is validated and verified against both experimental measurements \citep{wu_experimental_2019} and numerical results from MoorDyn \citep{chen_cfd_2022} for a floating box moored with catenary lines under regular wave conditions, as part of an introductory test campaign for wave energy converters. The rest of the paper is organised as follows. Section \ref{sec:math_model} provides the details of the mathematical model of the moored floating body problem and the coupling strategy. Section \ref{sec:benchmarking} elaborates on the benchmarking of the proposed algorithm. Finally, some conclusions are drawn in Section \ref{sec:conclusion}.

\section{Mathematical Model and Numerical Methods}\label{sec:math_model}
The proposed approach utilises OpenFOAM \citep{weller_tensorial_1998} as the main framework for CFD analysis, while the beam solver, also implemented within the OpenFOAM framework, is used for mooring line dynamics modelling. Consequently, the equation of motion for the floating rigid body and its interaction with mooring lines must be coupled with the CFD model. This section elaborates on the mathematical models that have been used for each component of a moored floating body problem.
\subsection{Free surface flows}
The volume of fluid (VOF) formulation based on the surface capturing method with Reynolds Averaging Navier-Stokes (RANS) equations is employed here 
for solving a two-phase air-water system. The continuity and the conservation of momentum equations are given by: 
\begin{equation}
	\nabla \cdot\mathbf{U}=0
	\label{eqn:1}
	\end{equation}
	\begin{equation}
	\frac{\partial \rho \mathbf{U}}{\partial t}+ \nabla \cdot(\rho \mathbf{UU}) - \nabla \cdot(\mu_{eff}\nabla \mathbf{U})= -\nabla p^{*} - \mathbf{g} \cdot\mathbf{X}\nabla \rho + \nabla \mathbf{U}\cdot \nabla \mu_{eff}
	\label{eqn:2}
\end{equation}
where $\mathbf{U}$ is the fluid velocity vector in Cartesian coordinates, $\rho$ is the density of the mixture, $p^{*}=p-\rho g \cdot \mathbf{X}$ is the dynamic pressure while $p$ is the total pressure, $\mathbf{g}$ is the gravitational acceleration, $\mathbf{X}$ represents the position vector of the points on the computational grid and $\mu_{eff}=\mu + \mu_{t}$ is the effective dynamic viscosity which is the sum of molecular dynamic viscosity ($\mu$) and turbulent dynamic viscosity ($\mu_{t}$). The core concept of the VOF method is to introduce the volumetric fraction (of the primary phase),
$\alpha$, to model an immiscible two-phase system. Then, the two-phase system can be treated as a single-phase system by introducing the mixture density and mixture viscosity:
\begin{equation}
\rho=\alpha \rho_{1}+(1-\alpha)\rho_{2}
\label{eqn:3}
\end{equation}
\begin{equation}
\mu=\alpha \mu_{1}+(1-\alpha)\mu_{2}
\label{eqn:4}
\end{equation}
\noindent
where subscripts 1 and 2 denote the density/viscosity of primary and secondary phases, respectively.
To calculate the volume fraction, a transport equation is solved as follows:
\begin{equation}
\frac{\partial \alpha}{\partial t} + \nabla .(\mathbf{U} \alpha) =0
\label{eqn:5}
\end{equation}

\noindent Physically, Eqn. \ref{eqn:5} represents the conservation of mass for one of the phases in a multiphase flow, ensuring that the volume fraction evolves correctly over time.
The Multidimensional Universal Limiter with Explicit Solution (MULES) algorithm \citep{damian_extended_2014} is applied here to reduce interface diffusion. This method introduces an artificial compression (or anti-diffusion) term $\nabla \cdot [\mathbf{U_{r}}\alpha(1-\alpha)]$ on the left side of Eqn. \ref{eqn:5}. In this term, $\mathbf{U_{r}}$ represents the compression velocity. The final expression for the volume fraction advection equation then becomes:
\begin{equation}
\frac{\partial \alpha}{\partial t} + \nabla .(\mathbf{U} \alpha) + \nabla \cdot [\mathbf{U_{r}}\alpha(1-\alpha)]=0
\label{eqn:6}
\end{equation}
Further details on the MULES formulation can be found in \citet{damian_extended_2014}.

\subsection{Mooring model}
\label{subsec:beamSection}
\begin{figure}[!h]
\includegraphics[width=0.9\textwidth]{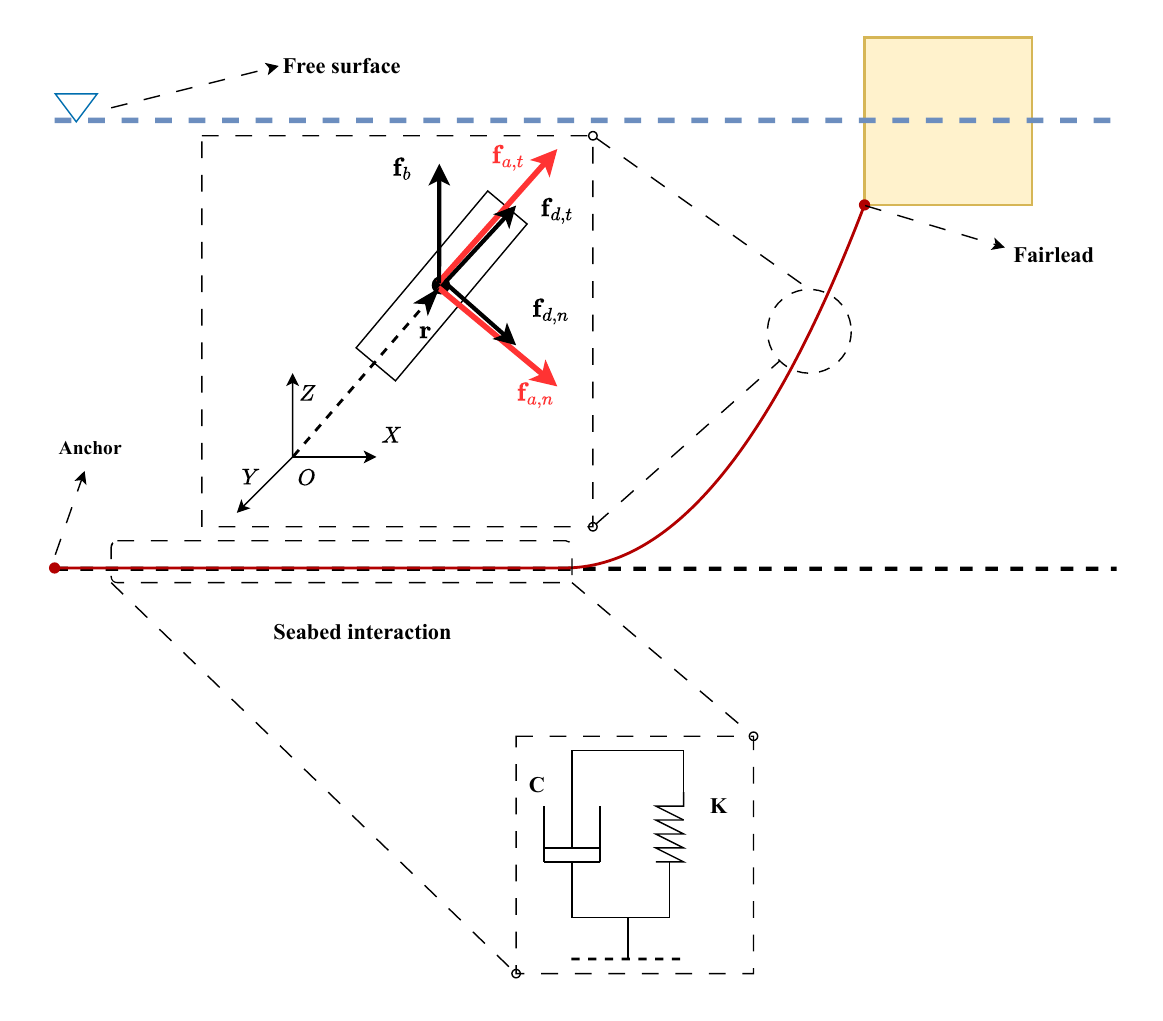}
\caption{Forces acting on a beam segment for a typical mooring line setup}
\label{fig:forcesonbeam}
\centering
\end{figure}

A beam of length $L$, which represents the mooring line, is described by a position vector $\mathbf{r} = \{ r_{x}(s), r_{y}(s), r_{z}(s) \} $ (see Figure \ref{fig:forcesonbeam}) in a cartesian coordinate system. Here, $0 \leq s \leq L$ represents the arc length along the undeformed configuration of the beam. In differential form, the conservation of linear and angular momentum can be written as
\begin{equation}
    \frac{\partial\mathbf{n}}{\partial s} + \mathbf{f}_{ext} = \rho_{b} A \frac{\partial^2 \mathbf{r}}{\partial t^2}
    \label{eqn:forcebalanceEqn}
\end{equation}
\begin{equation}
    \frac{\partial \mathbf{m}}{\partial s}
    \;+\; \mathbf{r}' \times \mathbf{n}
    \;+\; \mathbf{t}
    \;=\; \mathbf{I}_\rho \,\frac{\partial \boldsymbol{\omega}}{\partial t}
    \;+\; \boldsymbol{\omega} \times \!\left( \mathbf{I}_\rho \,\boldsymbol{\omega} \right)    \label{eqn:momentBalanceEqn}
\end{equation}
where in Eqn. \ref{eqn:forcebalanceEqn}, $\mathbf{n}(s,t)$ represents the internal force resultant vector at position $s$ and time $t$, $\mathbf{f}_{ext}(s,t)$ is the external force per unit length acting on the beam at position $s$ and time $t$, $\rho_{b}$ is the density of the beam, and $A$ is the cross-sectional area of the beam. Consequently, $\rho_{b} A$ is the mass per unit length of the beam, $\frac{\partial^2 \mathbf{r}}{\partial t^2}$ is the linear acceleration of the point on the centre line at position $s$.
In the moment equation (Eqn. \ref{eqn:momentBalanceEqn}), $\mathbf{m}$ is the internal moment vector, $\mathbf{r}'$ is the tangent vector along the beam's centre line, $\boldsymbol{\tau}_{\text{ext}}(s,t)$ represents the external torques, $\mathbf{I}_\rho$ is the spatial time-dependent inertia tensor of the beam's cross-section and $\boldsymbol{\omega}$ is the angular velocity vector. 
Equation \ref{eqn:forcebalanceEqn} represents the conservation of linear momentum, where the net internal and external forces acting on a beam segment result in its linear acceleration. Similarly, Eqn. \ref{eqn:momentBalanceEqn} represents the conservation of angular momentum, accounting for the net internal moments and external torques acting on the segment.

\noindent The external forces, represented by the term $\mathbf{f}_{ext}(s,t)$ in Eqn. \ref{eqn:forcebalanceEqn}, can be expressed as:
\begin{equation}
    \mathbf{f}_{ext} = \mathbf{f}_{d} + \mathbf{f}_{a} + \mathbf{f}_{b} + \mathbf{f}_{gc}
\end{equation}
in which, $\mathbf{f}_{d}$ represents the drag force acting on the beam segment, $\mathbf{f}_{a}$ is the added mass force, which is the inertial force from the surrounding fluid accelerating with a moving body, $\mathbf{f}_{b}$ is the buoyancy force, $\mathbf{f}_{gc}$ accounts for ground contact force. These forces for a mooring line with a circular cross-section are also illustrated in Figure \ref{fig:forcesonbeam}.

\noindent The drag force $\mathbf{f}_{d}$ based on Morrison's equation \citep{morison_force_1950}, which describes the total hydrodynamic forces acting on a structure submerged in or interacting with fluid flow, comprises normal and tangent components expressed as \cite{palm_coupled_2016}
\begin{equation}
    \mathbf{f}_{d} = \frac{1}{2}\rho d\sqrt{1+\epsilon}\left( \vphantom{\frac{1}{2}} 
    C_{D,t} \left\|  \mathbf{V}_{r,t} \right\| \mathbf{V}_{r,t} + 
    C_{D,n} \left\|  \mathbf{V}_{r,n} \right\| \mathbf{V}_{r,n} 
    \right)
    \label{eqn:drag}
\end{equation}
where $d$ is the beam's diameter, $C_{D,t}$ and $C_{D,n}$ are the drag coefficients in tangential and normal directions, respectively,  $\boldsymbol{\epsilon}$ is the unit tangential vector and can be expressed based on the position vector as
\begin{equation}
    \epsilon =  \left|\left| \frac{\partial \mathbf{r}}{\partial s} \right|\right| - 1 
     \label{eqn:epsilonEquation}
\end{equation}
$\mathbf{V}_{r,t}$ and $\mathbf{V}_{r,n}$ are the tangential and normal components of the relative velocity between the beam and the surrounding fluid, i.e. $\mathbf{V}_{r} = \mathbf{V} - \mathbf{V}_{fluid}$. In this study, a quiescent fluid is assumed around the mooring line, i.e. $\mathbf{V}_{fluid} = 0$.
Similar to the drag force, the added mass force is one of the components in the Morison equation \citep{morison_force_1950}. The added mass force arises because when a body moves through a fluid, it must push aside or accelerate the surrounding fluid. This extra inertia, which the body feels due to the fluid's resistance, is called added mass. This term increases the effective mass of the beam because it includes both the actual mass of the beam and the mass of the fluid that the beam displaces during motion. Following the similar notation in Eqn. \ref{eqn:drag}, the added mass force can be expressed as
\begin{equation}
    \mathbf{f}_{a} = \rho_{f} A (C_{M,t} \mathbf{a}_{r,t} + C_{M,n}\mathbf{a}_{r,n} + \mathbf{a}_{f})
\end{equation}
Where $C_{M,t}$ and $C_{M,n}$ are the tangential and normal added mass coefficient, respectively, and $\mathbf{a}_r$ is the relative acceleration between the beam and fluid, i.e. $\mathbf{a}_{r} = \mathbf{a} - \mathbf{a}_{f}$.

\noindent The total buoyancy force acting on the beam, proportional to the displaced fluid, can be expressed as
\citep{palm_coupled_2016}
\begin{equation}
\mathbf{f}_{b} = \rho_{b} A\left(\frac{\rho_{b} - \rho_{f}}{\rho_{b}}\right) \mathbf{g}
\label{eqn:buoyancy}
\end{equation}
Since the beam is more dense than the fluid, its weight exceeds the buoyant force, resulting in a net downward force.
As illustrated in Figure \ref{fig:forcesonbeam}, the linear spring-damper model has been used for the interaction between the seabed and the beams. For any segment of the beam in contact with the seabed, a vertical reaction force is applied as follows
\begin{equation}
    \mathbf{f}_{gc}^{n} =
  \begin{cases}
    \left[ K_{n}d \left(r_{z} - z_{g}\right) - C d \max \left(V_{z}, 0 \right)\right] \mathbf{n}_{\text{g}} & \text{if } z_{b} \leq z_{g} \\
    \mathbf{0} & \text{otherwise}
  \end{cases}
  \label{eqn:groundContact}
\end{equation}
Where $K_{n}$ and $C$ are the seabed stiffness in the normal direction and damping coefficients, respectively, and $z_{g}$ represents the seabed elevation. The unit vector $\mathbf{n}_{\text{g}} = \left(0,0,1 \right)$ denotes the seabed normal direction and ensures that the resulting contact force acts upward from the seabed. In addition to the seabed contact normal force, a tangential friction force is applied to any segment of the beam in contact with the seabed. This force follows the Coulomb friction model, opposing the relative tangential velocity between the beam and the seabed. The magnitude and direction of the friction force depend on the tangential velocity $\mathbf{V}_{t}$ of the beam segment, the normal contact force $\mathbf{f}_{gc}^{n}$, and the friction coefficient $\mu$, as defined by the following piecewise expression:

\begin{equation}
\mathbf{f}_{gc}^{t} =
\begin{cases}
- \mu \left| \mathbf{f}_{gc}^{n} \right| \dfrac{\mathbf{V}_{t}}{\left| \mathbf{V}_{t} \right|}, & \text{if } K_{t}d \left| \mathbf{V}_{t} \right| \geq \mu \left| \mathbf{f}_{gc}^{n} \right| \\
- K_{t}d \mathbf{V}_{t}, & \text{otherwise}
\end{cases}
\label{eqn:tangentialFriction}
\end{equation}

\noindent Here, $K_{t}$ is a tangential stiffness parameter that enables a transition between the sticking and slipping regimes. When the relative tangential velocity is small, the friction force behaves in a sticking manner. Once the velocity exceeds a threshold defined by the Coulomb limit, the model transitions to a constant-magnitude slipping friction force oriented opposite to the motion.

\noindent
The total seabed contact force is then the sum of the normal and tangential components:

\begin{equation}
\mathbf{f}_{gc} = \mathbf{f}_{gc}^{n} + \mathbf{f}_{gc}^{t}
\end{equation}

\begin{figure}[!h]
\includegraphics[width=0.8\textwidth]{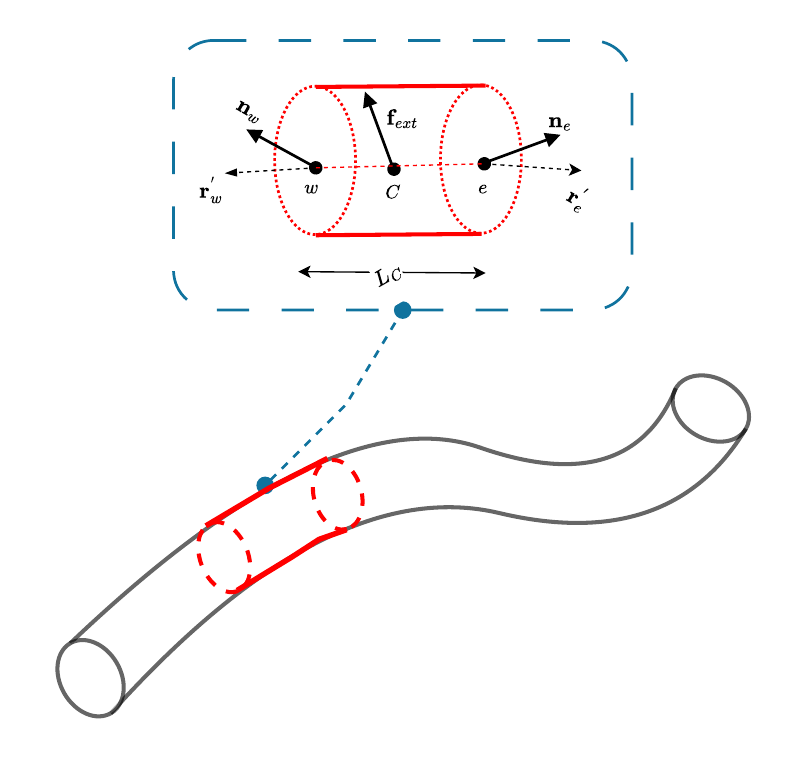}
\caption{Balance of internal and external forces on a control volume}
\label{fig:beamelement}
\centering
\end{figure}

\noindent For an isolated control volume shown in the Figure 
\ref{fig:beamelement}
The integral form of the conservation of linear momentum (Eqn. \ref{eqn:forcebalanceEqn}) can be discretised over a control volume as follows
\begin{equation}
    \begin{aligned}
    \int_{L} \mathbf{n}dL = n_{e} - n_{w} \quad ; \quad & \int_{L} \mathbf{f}_{ext}dL \approx  \mathbf{f}_{ext,c}L_{c} \\
    & \Rightarrow \mathbf{n}_{e} - \mathbf{n}_{w} + \mathbf{f}_{ext,c} L_{c} = \rho A \ddot{r_{c}}
    \end{aligned}
\label{eqn:beamdisc}
\end{equation}
Where $\mathbf{n}_{e}$ and $\mathbf{n}_{w}$ are the force values evaluated at cell face $e$ and $w$, respectively, and the subscript $c$ represents the values at the cell-centre. The term $\mathbf{f}_{ext}$ is assumed to have a linear variation across the control volume and can therefore be approximated using the midpoint rule by evaluating the function at the midpoint of each control volume.

\noindent The discretised equations for each control volume result in a system of the form  
\begin{equation}
\mathbf{A}_w \mathbf{X}_w + \mathbf{A}_c \mathbf{X}_c + \mathbf{A}_e \mathbf{X}_e = \mathbf{B}_c,
\label{eqn:finalForm}
\end{equation}
where $\mathbf{A}_{w}$, $\mathbf{A}_{c}$, and $\mathbf{A}_{e}$ are $6 \times 6$ coefficient matrices; 
$\mathbf{X}_{w}$, $\mathbf{X}_{c}$, and $\mathbf{X}_{e}$ are $6 \times 1$ vectors of unknown displacements and rotations, and $\mathbf{B}_c$ is the corresponding $6 \times 1$ source vector. This form is assembled for each beam cell, leading to a global system of $6 \times N$ coupled equations. The resulting global matrix is sparse and block tri-diagonal, and the non-linear system is solved using the Newton–Raphson method. The detailed formulation of the Simo-Reissner beam theory, alongside its discretisation, can be found in Bali et al. \citep{bali_cell-centered_2022}.

\subsection{Rigid body motion dynamics}
The motion of the rigid floating body is determined using Newton’s second law of motion, which governs the translational and rotational dynamics of the body. The external forces and moments considered include hydrodynamic forces, hydrostatic forces, gravitational forces, mooring line forces, and any additional damping or restoring forces that may be present. All these forces except gravity and mooring forces are integrated over the body’s surface to account for the influence of the surrounding fluid, while gravity is treated as a body force and integrated over the volume of the body. Mooring loads are applied as concentrated external forces and moments at the attachment points.
To ensure stability and accuracy, the equations of motion are solved in a coupled manner, fully accounting for the interaction between the rigid body and the surrounding fluid. This approach captures the mutual influence of fluid forces on the body’s motion and the body’s motion on the fluid flow.
In OpenFOAM, the \texttt{sixDoFRigidBodyMotion} library is used to solve the motion of a body in six degrees of freedom (DoF). The conservation of linear and angular momentum, with rotation occurring around the centre of mass (COM), leads to the following equations for the floating body,
\begin{equation}
\frac{\partial \mathbf{V}_{rb}}{\partial t}=\frac{\mathbf{F}_{rb}}{m_{rb}}
\label{eqn:6}
\end{equation}
\begin{equation}
\frac{\partial \boldsymbol{\omega}_{rb}}{\partial t}= \mathbf{I}^{-1}_{rb}\cdot[\mathbf{M}_{rb}-\boldsymbol{\omega}_{rb} \times (\mathbf{I}_{rb}\boldsymbol{\omega}_{rb})]
\label{eqn:7asli}
\end{equation}
where $\mathbf{V}_{rb}$ and $\boldsymbol{\omega}_{rb}$ denotes the linear and angular velocity of the floating rigid body, $m_{rb}$ is the mass, $\mathbf{I}_{rb}$ is the moment of inertia of the body, $\mathbf{F}_{rb}$ and $\mathbf{M}_{rb}$ are the external forces and moments which act on the floating rigid body:
\begin{equation}
\mathbf{F}_{rb}=\oint_{S} (p\mathbf{\mathbb{I}} + \boldsymbol{\tau}) \cdot d\mathbf{S} \;+\; \mathbf{F}_{moor} \;+\; m_{f}\mathbf{g}
\label{eqn:7}
\end{equation}
\begin{equation}
\mathbf{M}_{rb}=\oint_{S} \mathbf{r}_{1}\times(p\mathbf{\mathbb{I}}+\boldsymbol{\tau})\cdot d\mathbf{S} \;+\; \mathbf{r}_{2}\times \mathbf{F}_{moor} \;+\; \mathbf{r}_{3}\times m_{f}\mathbf{g}
\label{eqn:8}
\end{equation}
In which $\mathbf{\mathbb{I}}$ denotes the $3\times 3$ identity matrix, $p$ is the normal pressure acting perpendicular to the surface of the floating body, and $\boldsymbol{\tau}$ represents the viscous stress term, which is the tangential force due to the fluid's viscosity, opposing the motion of the body, $\mathbf{S}$ is rigid body surface while 
$\mathbf{F}_{moor}$ is the resultant restraining force transmitted by the mooring system to the floating body, obtained from the beam solver described in Section {\ref{subsec:beamSection}}. Note that this is distinct from $\mathbf{f}_{ext}$ in the beam equations, which denotes the distributed external load acting along the mooring line, $\mathbf{r}_{1}$, $\mathbf{r}_{2}$, $\mathbf{r}_{3}$ represent the moment arms of the hydrodynamic force, mooring restraining force, and gravitational force acting on the floating object, respectively.
\subsection{Morphing mesh method}
\begin{figure}[!h]
\includegraphics[width=0.9\textwidth]{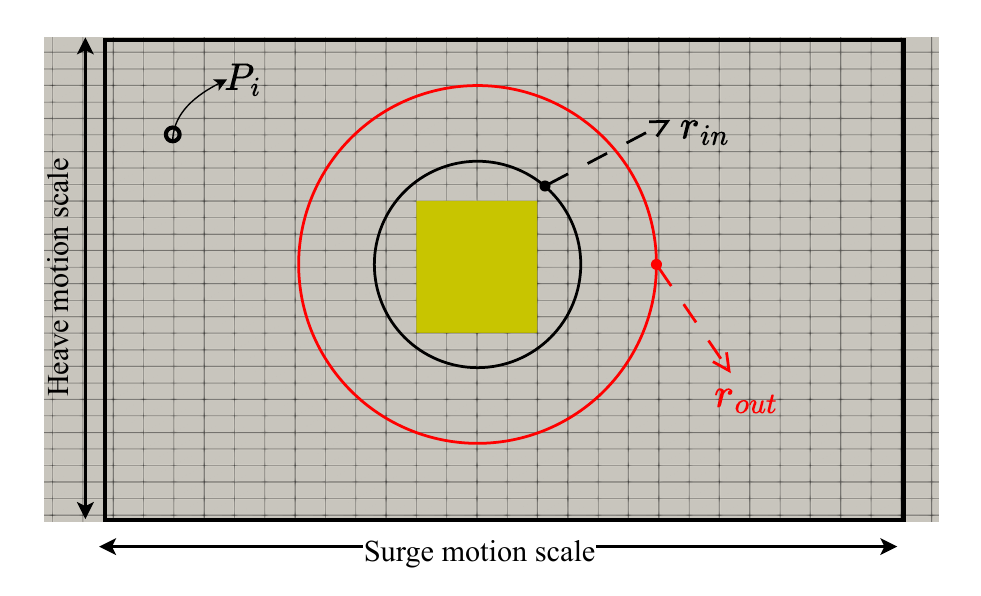}
\caption{Parameters of modified mesh morphing method }
\label{fig:mmm-param}
\centering
\end{figure}

In the present study, the deforming mesh and dynamic motion of the floating rigid body is controlled using a modified version of the morphing mesh suggested by \citet{palm_facilitating_2022}, in which the state of the floating rigid body is defined by its position, \( \mathbf{B}(t) \), and orientation, \( \mathbf{q}(t) \). The orientation is expressed in quaternion notation, assuming that the body has the initial orientation of $\mathbf{q}_{0}$. For an arbitrary point $\mathbf{P}_{i}$ in the CFD mesh (see Figure \ref{fig:mmm-param}), treated as a pure quaternion with zero scalar part for the purpose of rotation, the point displacement $\boldsymbol{\delta p}_i$ can be computed as described in \citep{palm_facilitating_2022}.
\begin{equation}
\mathbf{\boldsymbol{\delta} p}_{i} = \boldsymbol{\alpha}_{i} \odot \mathbf{b} + \mathbf{q}_{\beta_{i}}\mathbf{P}_{i}\mathbf{q}_{\beta_{i}}^{*}
\label{eqn:mmm1}
\end{equation}
\begin{equation}
    \beta_{i} =
  \begin{cases}
    1 & \text{if $W(\mathbf{P}_{i}) < r_{in}$ } \\
    1-\frac{W({\mathbf{P}_{i}})-r_{in}}{r_{out}-r_{in}} & \text{otherwise} \\
    0 & \text{if $W(\mathbf{P}_{i}) > r_{out}$ }
  \end{cases}
  \label{eqn:mmm2}
\end{equation}
\begin{equation}
\mathbf{q}_{\beta_{i}} = (\mathbf{q}\mathbf{q}_{0}^{-1})^{\beta_{i}}\mathbf{q_{0}} 
\label{eqn:rotationNative}
\end{equation}
where $\beta_i$ is the interpolation weight of point $i$, $\boldsymbol{\alpha}_{i} = [\alpha_{i}^{(x)}, \alpha_{i}^{(y)}, \alpha_{i}^{(z)}]^\top$ is the vector of scaling factors for translational motions in surge ($x$), sway ($y$), and heave ($z$), respectively, and $\mathbf{b} = [b_{x}, b_{y}, b_{z}]^\top$ represents the body displacement vector. The operator $\odot$ denotes element-wise (Hadamard) multiplication, the superscript asterisk (*) indicates the quaternion conjugate, and $W(\mathbf{P}{i})$ is the minimum distance from $\mathbf{P}{i}$ to the surface of the body. Two parameters $r_{in}$ and $r_{out}$ have been employed to define the inner and outer regions, respectively. Within the $r_{in}$ region, $\mathbf{P}_{i}$ rigidly follows the body, while for the points outside the outer region ($W(\mathbf{P}_{i}) > r_{out}$), the mesh remains static.

\begin{figure}[!h]
    \centering
    \subfloat[Modified mesh morphing\label{fig:mmm}]{\includegraphics[width=0.8\textwidth]{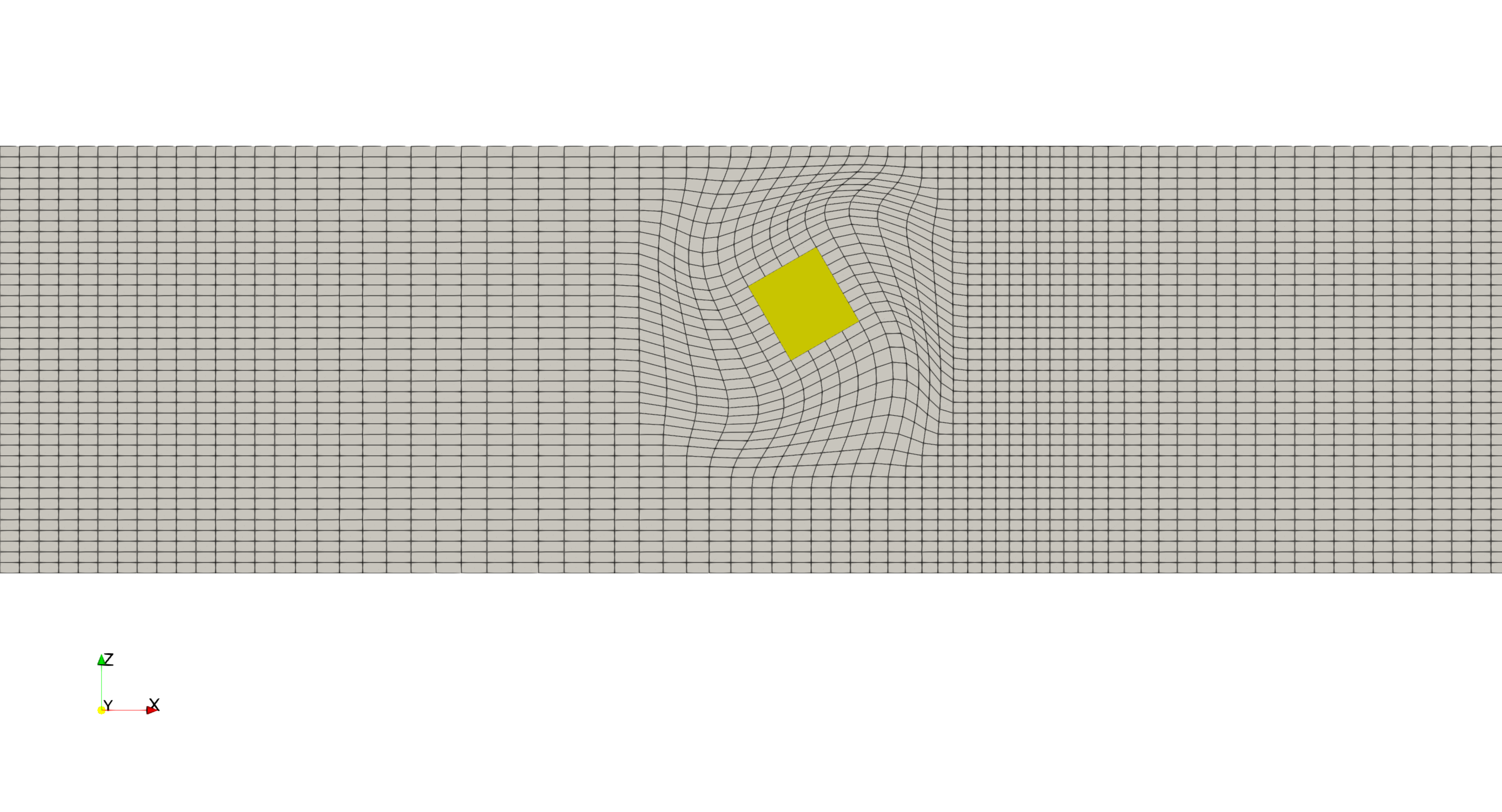}} \hspace{1em}
    \subfloat[Native mesh morphing\label{fig:omm}]{\includegraphics[width=0.8\textwidth]{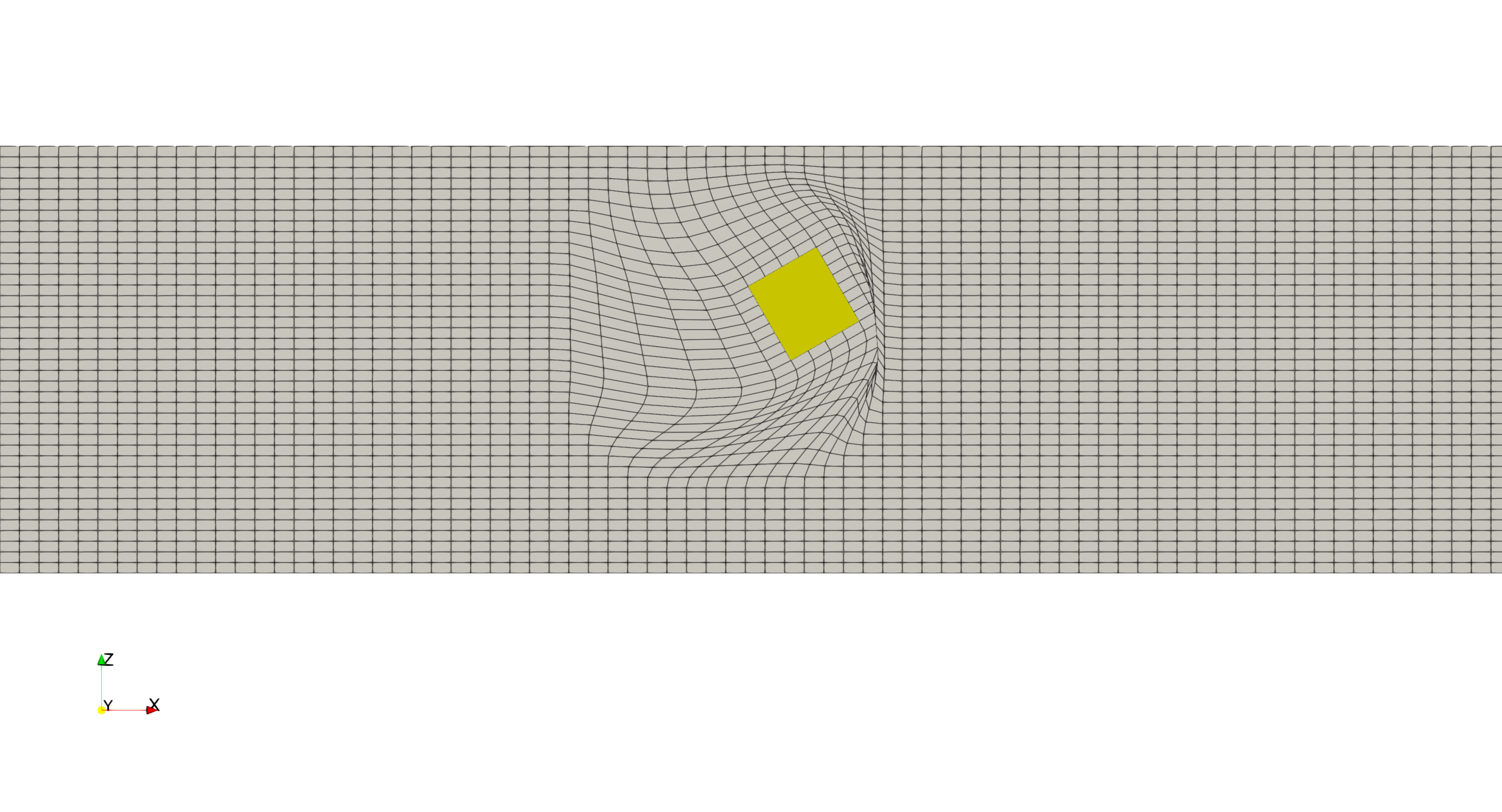}} \\[1em]
    \caption{Comparison between modified mesh morphing and native mesh morphing techniques}
    \label{fig:parallelStudy}
\end{figure}

\noindent The native mesh morphing class in OpenFOAM results in poor mesh quality and grid distortion due to high-amplitude rotations (see Figure \ref{fig:omm}), while the modified mesh morphing (Eqns. \ref{eqn:mmm1}, \ref{eqn:mmm2} and \ref{eqn:rotationNative}) improves upon the original by adding an additional mesh morphing region to accommodate one or two translational degrees of freedom, specifically surge and/or heave. As a result, the morphing region defined by \(r_{in}\) and \(r_{out}\) translates the rigid body horizontally. By excluding translational motion from the inner deformation zone, the method generates deformation patterns that accurately reflect the combined effects of heave and pitch, without introducing spurious distortions (see Figure \ref{fig:mmm}). This improvement enhances the stability of numerical simulations, allowing for the successful simulation of the box's response to the incoming wave. For a detailed discussion, see \citet{palm_facilitating_2022}.
 
\subsection{Coupling Strategy}
The coupling between the moored floating body and the multiphase free-surface flow is achieved through the PIMPLE algorithm. PIMPLE is a hybrid pressure–velocity coupling algorithm that combines features of the Pressure Implicit with Splitting of Operators (PISO) algorithm and the Semi-Implicit Method for Pressure-Linked Equations (SIMPLE) algorithm. It is specifically designed to handle transient simulations with larger time steps than those used in PISO, while still maintaining numerical stability and accuracy.

\begin{figure}[!h]
\includegraphics[width=\textwidth]{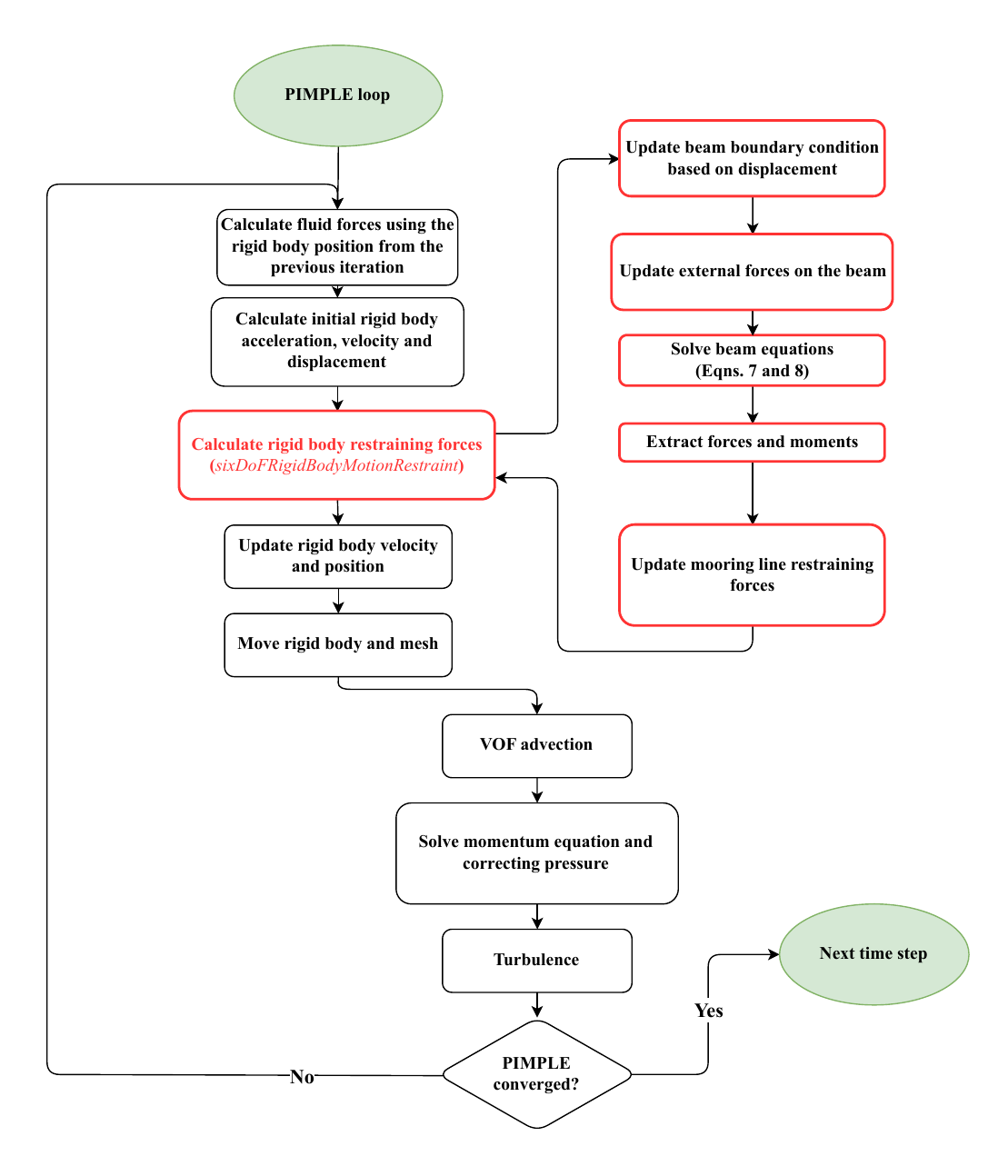}
\caption{Procedure for coupling the beam solver and fluid solver}
\centering
\label{fig:algo}
\end{figure}

\noindent Within each PIMPLE iteration, the beam equations are solved first using the position of the floating body from the previous iteration (or from the previous time step, in the first PIMPLE iteration). The beam solver calculates the mooring forces and updates the position and velocity of the mooring lines. These restraining forces and moments are applied at the fairlead location and passed to the six-degree-of-freedom (6DOF) solver to update the acceleration of the floating body.
The motion of the floating body is then solved by integrating the equations of motion (Eqns. {\ref{eqn:6}}, {\ref{eqn:7asli}}). Based on the updated position, the computational mesh is deformed using the morphing strategy described in Eqn. \ref{eqn:mmm1}.
The position of the water-air interface is then updated by solving the advection equation for the volume fraction (Eqn. \ref{eqn:5}). Finally, the velocity-pressure coupling is resolved within the PISO sub-loop, including updates to turbulence quantities if applicable.
Each time step consists of 3 PIMPLE iterations to update the fluid solver, the floating body motion, and the mooring system state. A summary of this coupling strategy is shown in Figure \ref{fig:algo}.
\section{Benchmark Case}\label{sec:benchmarking}
To validate the proposed approach, a numerical model is created based on
the experimental 
work of \citet{wu_experimental_2019}. 
This experimental setup
was part of the European MaRINET2-EsflOWC project \citep{kisacik_efficiency_2020}, in which responses 
and motions of mooring lines attached to a solid floating box are investigated. Details of the numerical model will be provided later in this section, once the experimental setup is described.

\subsection{Experimental Setup}
\begin{figure}[!h]
    \centering
    \begin{subfigure}[b]{0.55\textwidth}
        \centering
        \includegraphics[width=\linewidth]{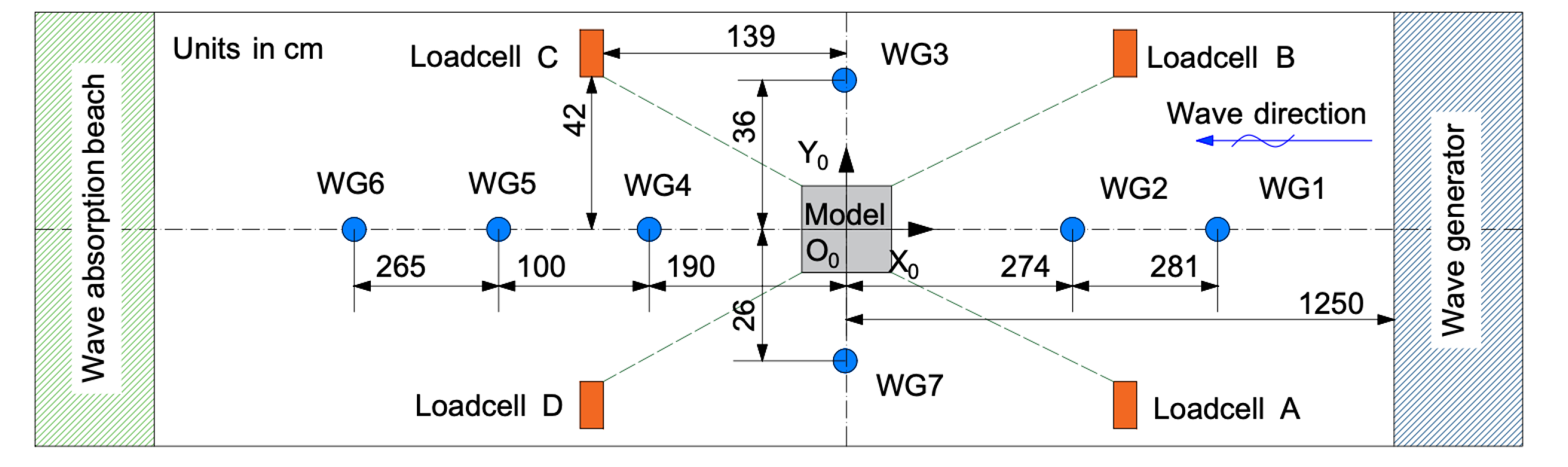}
        \caption{}
    \end{subfigure}%
    \begin{subfigure}[b]{0.55\textwidth}
        \centering
        \includegraphics[width=0.5\linewidth]{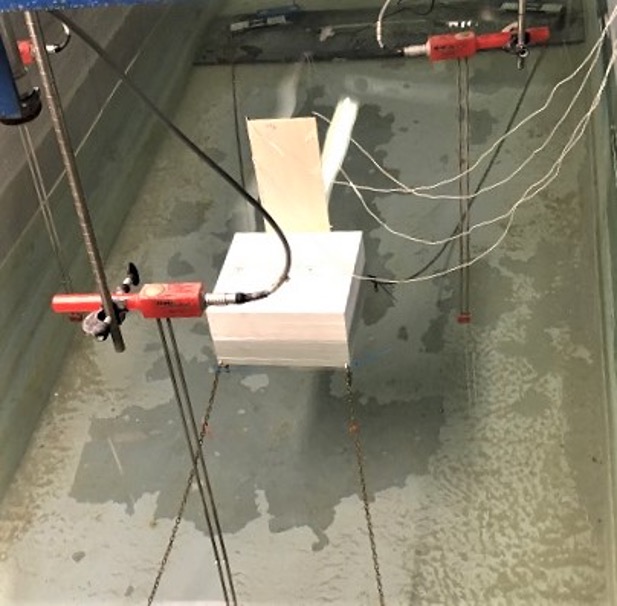}
        \caption{}
    \end{subfigure}%
    \caption{Overview of the experimental domain, with both images taken from the experimental work by \citet{wu_experimental_2019}: (a) top view of the domain, including the locations of wave gauges, and (b) floating box with mooring lines and attached plate.}
    \label{fig:experimentaldomain}
\end{figure}

The experimental model consists of a box constructed from lightweight PVC, moored symmetrically to the bottom of a $30.0$ m long, $1.0$ m wide, and 
$1.2$ m high wave flume belonging to the Coastal Engineering Research Group at the Department of Civil Engineering 
of Ghent University. During this project, several tests with different wave conditions were carried out. 
The rigid box made of light PVC material is $0.2$ m long, $0.2$ m wide, and $0.132$ m high and has a $3.148$ kg mass. To track the box's motion, a $0.324$ m plate was attached to the front face of the box, facing the incoming waves, on which reflective markers were installed.
The mass of the plate and markers is included in the reported properties of the box (Table {\ref{table:properties}}), and their influence is further discussed in Section {\ref{subsec:boxmotions}}. These reflective markers were used by a 6-DOF optical motion tracking system. Initially, the box had a draft length of $0.0786$ m, which meant 60\% of its height was submerged.
Figure \ref{fig:experimentaldomain} depicts a schematic view of the experimental domain.

\noindent The mooring system employed in this experiment includes four symmetrical iron chains initially slack and postured in a catenary layout. Four load cells have been mounted
on the anchor end of the mooring lines to record the tensions. On the fairlead side, each line is connected at a point 0.5 cm above the box’s corner. \citep{dominguez_sph_2019}. 
A detailed list of box and mooring line properties is provided in Table \ref{table:properties}. 
The reported centre of gravity values in Table 1 already account for the attached wooden plate ($0.324$ m, front face) and reflective markers, ensuring that the listed values correspond to the full experimental model.
\begin{table}[!h]
    \caption{Properties experimentally tested box model with mooring lines \citep{chen_cfd_2022} }
    \begin{tabular}{l l l}
        \textbf{Properties} & \textbf{Value} &\textbf{Units} \\
        \midrule
        Box dimension ($L \times W \times H$) & $0.2 \times 0.2 \times 0.132$ & m\\
        Box mass (+ Connections) & $3.16$ & kg \\
        Centre of gravity ($x$,$y$,$z$) & $(0,0,-0.126)$ & m \\
        Box initial draft & $0.0786$ & m \\
        Box moment of inertia ($\mathbf{I}_{xx}, \mathbf{I}_{yy}, \mathbf{I}_{zz}$) & $(0.015, 0.015, 0.021)$ & kg m$^{2}$ \\
        Box material & PVC lightweight & - \\
        Attached wooden plate & 0.324 & m \\
        Mooring line diameter & $0.003656$ & m \\
        Mooring line length & $1.455$ & m \\
        Mooring line density & $5782.10$ & kg m$^{-3}$ \\
		Mooring line Young's modulus & $1.82$ & MPa \\
		Mooring line Poisson's ratio & $0.5$ & - \\
    \end{tabular}
    \label{table:properties}
\end{table}
The coordinates for anchors and fairleads are presented in Table \ref{table:fairleadAnchorCoords}. The experimental data include the motion of the moored body, the tension in the mooring lines, and the surface elevation at several locations in the wave tank. These data are used to validate the numerical model presented in Section \ref{sec:math_model}.
The surface elevation of the water is recorded using resistive wave gauges (WGs). Table \ref{table:wavegaugescoords} lists the WG positions within the wave tank. WG1 is located between the upstream and the box,
WG2 and WG3 are located between the box and the side walls of the wave tank, and WG4, WG5 and WG6 are located downstream and before the wave absorption area.

\begin{table}[!ht]
    \caption{Fairleads and anchors coordinates as reported in \citet{wu_experimental_2019} }
    \begin{tabular}{l l}
        \textbf{Fairlead/Anchor} & \textbf{Global Coordinates (m)}  \\
        \midrule
        Fairlead 1 & $(-0.1,0.1,0.0736)$ \\
        Fairlead 2 & $(-0.1,-0.1,-0.0736)$ \\ 
        Fairlead 3 & $(0.1,0.1,-0.0736)$ \\
        Fairlead 4 & $(0.1,-0.1,-0.0736)$ \\
        Anchor 1 & $(-1.385,0.423,-0.5)$ \\
        Anchor 2 & $(-1.385,-0.423,-0.5)$ \\
        Anchor 3 & $(1.385,0.423,-0.5)$ \\
        Anchor 4 & $(1.385,-0.423,-0.5)$ 
    \end{tabular}
    \label{table:fairleadAnchorCoords}
\end{table}

\begin{table}[!ht]
    \caption{Wave gauges coordinates within the wave tank as reported in \citet{wu_experimental_2019}} 
    \begin{tabular}{l  l}
        \textbf{WG number} & \textbf{Global Coordinates (m)}  \\
        \midrule
        Wave gauge 1 & $(-2.74 , 0)$ \\
        Wave gauge 2 & $(-0.05 , 0.26)$ \\
        Wave gauge 3 & $(0.07 , -0.36)$ \\
        Wave gauge 4 & $(0.55 , 0)$ \\
        Wave gauge 5 & $(1.90 , 0)$ \\
        Wave gauge 6 & $(2.90 , 0)$
    \end{tabular}
    \label{table:wavegaugescoords}
\end{table}

\subsection{Numerical model}
\begin{figure}[!ht]
     \centering
     \begin{subfigure}[b]{0.8\textwidth}
         \centering
         \includegraphics[width=\textwidth]{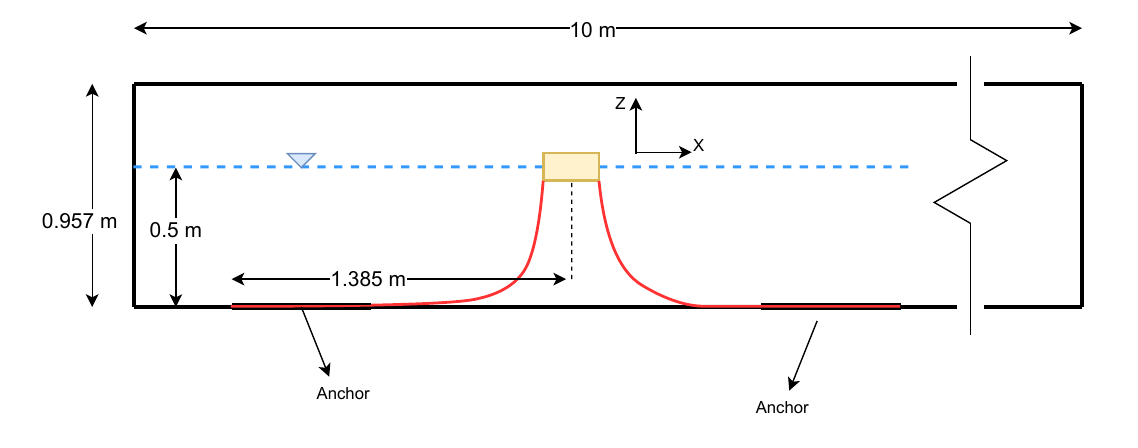}
         \caption{Side view}
         \label{fig:sideView}
     \end{subfigure}
     \hfill
     \begin{subfigure}[b]{0.8\textwidth}
         \centering
         \includegraphics[width=\textwidth]{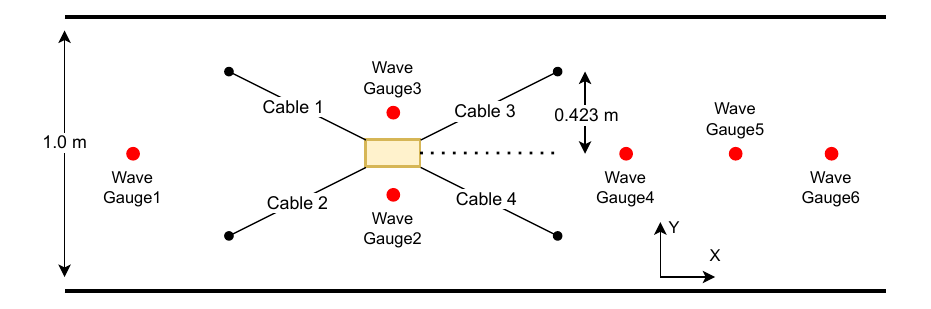}
         \caption{Plan view}
         \label{fig:planView}
     \end{subfigure}
        \caption{Numerical model of the European MaRINET2 EsflOWC \citep{kisacik_efficiency_2020} project benchmark case}
        \label{fig:NomDom}
\end{figure}

\begin{figure}[!ht]
\includegraphics[width=0.8\textwidth]{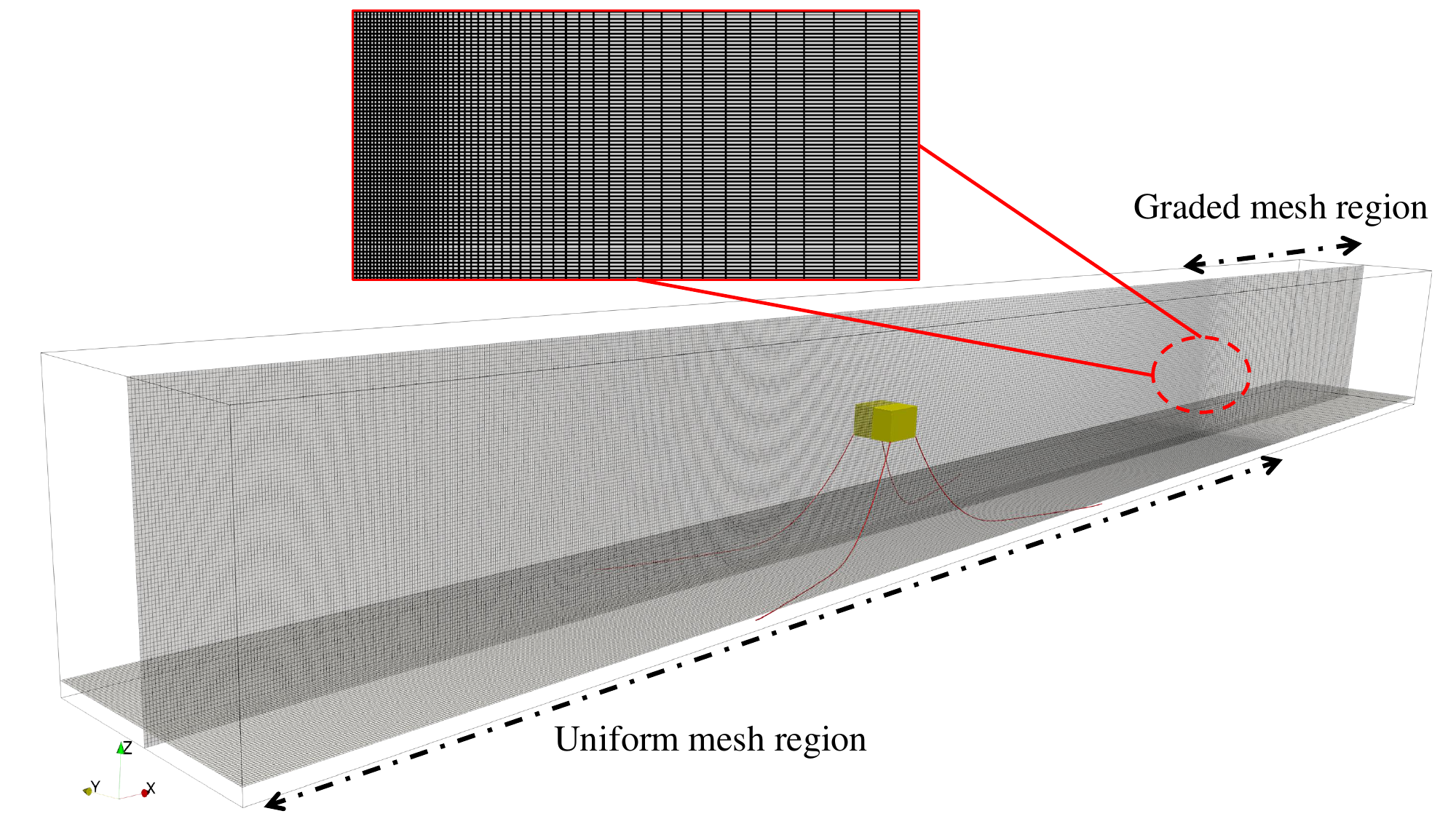}
\caption{Visualisation of the mesh structure, showing both uniform and graded regions}
\centering
\label{fig:meshOverview}
\end{figure}

Figure \ref{fig:NomDom} gives an overview of the numerical setup and the boundary conditions used in the numerical model. 
Following the approach employed in \citet{chen_cfd_2022}, 
rather than simulating the wave absorption beach in the physical wave flume, a cell stretching zone is designated in 
the numerical flume's last 3 meters, which decreases the total grid count and aids in wave damping.
The numerical grid is represented by a $10$ m $\times$ $1$ m $\times$ $0.957$ m (L $\times$ W $\times$ H) domain. Following the approach outlined in \citet{dominguez_sph_2019}, the numerical domain used in this study is shorter than the experimental domain to enhance computational efficiency. A grid convergence study was performed to identify an appropriate mesh resolution that balances accuracy with computational efficiency (see Appendix {\ref{sec:gridConvergence}} for details). Based on this study, a uniform grid (expansion ratio = $1$) with the size of $0.01$, $0.01$ and $0.011$ is applied in $x$,$y$ and the $z$ direction for the 70\% of the entrance length $x \in [-3.5,3.5]$ using multi-grading functionality of OpenFOAM that can divide a block in a given direction and apply different grading within each division. A graded mesh with an expansion ratio of $20$ was applied to the last 30\% of the domain. The resulting fluid mesh consists of approximately 5.8 million hexahedral cells. The structure of the numerical mesh is visualised in Figure \ref{fig:meshOverview}. 

\noindent Following previous numerical studies \citep{dominguez_sph_2019}, a Young’s modulus of $E = 2.78$ MPa was used for the mooring lines, resulting in an axial stiffness of $EA = 19$ N. Each mooring line was discretised into $60$ cells. The perpendicular and tangential added mass coefficients have been set to $1.6$ and $0$, respectively \citep{dominguez_sph_2019}. The drag coefficients in perpendicular and tangential directions have been set to $1.6$ and $0.5$, respectively \citep{dominguez_sph_2019}. For the normal component of ground contact with the seabed, A penalty stiffness of $1000\,\mathrm{N/m^2}$ and a damping coefficient of $1\,\mathrm{N\cdot s/m^2}$ are applied in the model, with a tangential friction coefficient of $0.01$; these values were chosen to significantly reduce penetration without compromising convergence.

\begin{figure}[!ht]
\includegraphics[width=0.8\textwidth]{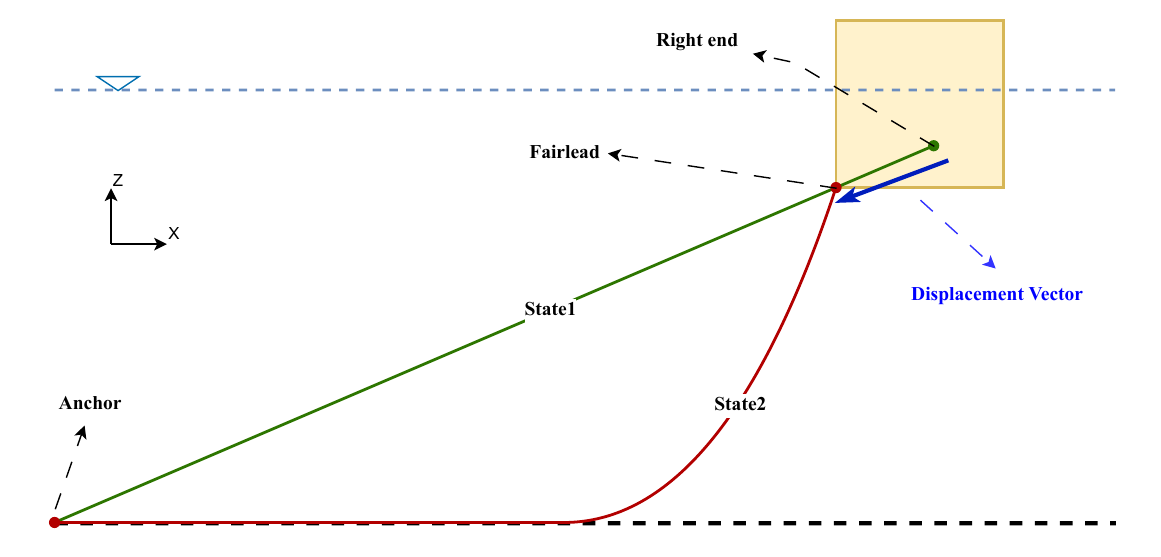}
\caption{Initilising the mooring lines, from State1 to State2}
\centering
\label{fig:beamInitialization}
\end{figure}

\noindent In order to initialise the mooring lines to have the correct initial posture and also match the pretension observed in the experimental work \citep{wu_experimental_2019}, a straight and inclined beam starting from the relevant anchor point and pointing toward the fairlead is considered. The beam equation is solved in a separate beam simulation using the beam solver, where the anchor end of the beam is constrained with fixed displacement but free to rotate, and a displacement boundary condition is applied to the opposite end. This displacement condition is defined by a vector extending from the free end of the beam to the fairlead location. Once the beam adopts its equilibrium posture, the resulting configuration is mapped to the coupled CFD framework and used as the initial condition for the beams. The initialisation process is illustrated in Figure \ref{fig:beamInitialization}, and results in an initial pretension of $0.38$ N in each mooring line.

\noindent A second-order Stokes wave generation and absorption are applied to the inlet and outlet boundaries, respectively. The parameters of the regular wave generation
are listed in Table \ref{tab:wave}, following the naming convention proposed by Chen and Hall \citep{chen_cfd_2022}. For both cases, the water depth is $0.5$ m. 
\begin{table}[!h]
    \caption{Wave parameters of the validation test cases }
    \label{tab:wave}
    \setlength{\tabcolsep}{4pt} 
    \renewcommand{\arraystretch}{1.1} 
    \begin{tabularx}{\textwidth}{l *{4}{>{\centering\arraybackslash}X}}
    \toprule
    \textbf{Case name} & \textbf{Nominal wave height (m)} & \textbf{FFT‐based wave height (m)} & \textbf{Wave period (s)} & \textbf{Wavelength (m)}\\
    \midrule
    H12T20 & 0.12 & 0.1088 & 2.0 & 4.06\\
    H15T18 & 0.15 & 0.1208 & 1.8 & 3.57\\
    \bottomrule
    \end{tabularx}
\end{table}

\noindent Rather than using the nominal heights 0.12 m (H12T20) and 0.15 m (H15T18) straight from the test report, the experimental wave-gauge signals were first analysed with a fast Fourier transform. The magnitude of the dominant spectral peak was taken as the true regular-wave height. This gives $0.1088$ m for H12T20 and $0.1208$ m for H15T18. These FFT-based heights were then prescribed at the inlet so that the simulated waves match the measured free-surface records.

\noindent The convection terms in this study were discretised using a second-order central scheme with Van Leer interpolation \citep{van_leer_towards_1974}, which has been used for the momentum and the volume fraction. This combination helps reduce numerical diffusion while limiting spurious oscillations near steep gradients or interfaces.
Gradient and diffusion operators were discretised using a linear second-order central difference scheme. The velocity gradient is computed using the cell-limited scheme, which constrains the gradient to ensure that when cell values are extrapolated to the faces, the resulting face values remain within the range of the surrounding cell values. 
The simulations use a surface-compression coefficient of $c_{\alpha} = 1$.
The time-derivative term is discretised using the Crank-Nicolson scheme with an off-centring coefficient of $0.9$. This coefficient, ranging between $0$ and $1$, controls the balance between implicit and explicit weighting. A coefficient of $1$ results in a fully centred, second-order accurate scheme, while a coefficient of $0$ reduces the scheme to the Euler-implicit method.  For the beam solver, a first-order Euler time discretisation scheme is employed.
The overall simulation domain has been decomposed into $120$ sub-domains for parallel processing, where each sub-domain is assigned to a different CPU core. The adjustable time step based on a maximum momentum Courant number of $0.5$ was used for time integration, as detailed in Appendix {\ref{sec:timeStepStudy}}. The interface Courant number has also been set to $0.5$. The momentum Courant number determines the main time step size, while the interface Courant number governs the sub-time-step used by the MULES algorithm to solve the volume fraction equation, ensuring accurate tracking of the phase interface without excessive movement per step. For the mesh morphing around the floating box, an outer distance of $0.45$ m and an inner distance of $0.05$ m have been considered. In addition, the modified morphing method utilises a $0.9$ m region in the x-direction and a $0.1$ m region in the z-direction.
The acceleration of the floating box was under-relaxed by a factor of $0.7$ to mitigate numerical instabilities in the rigid body motion procedure, providing a balance between solution accuracy and stability.

\noindent To test the parallel efficiency of the proposed coupled algorithm, a scaling study is conducted with different numbers of CPU cores for the mesh with $5.8$ million cells, using the Meluxina supercomputer equipped with 2x AMD EPYC Rome 7H12 64c {@} 2.6GHz CPUs. Details of the scalability study are presented in Table \ref{tab:parallelTab}. 
\begin{table}[!h]
    \caption{Details of Parallel strong scaling test cases}
    \begin{tabular}{l l l l}
        \textbf{Number of CPU cores} & \textbf{Execution Time (s)} &\textbf{Speedup} &\textbf{Cells per CPU core} \\
        \midrule
        $1$ & $2004.9$ & $-$ & $5\,800\,000$\\
        $2$ & $985.21$ & $2$ & $2\,900\,000$\\
        $4$ & $495.78$ & $4$ & $1\,450\,000$\\
        $8$ & $244.15$ & $8.2$ & $725\,000$\\
        $16$ & $139.85$ & $14.3$ & $362\,500 $\\
        $32$ & $93.33$ & $21.4$ & $181\,250$\\
        $64$ & $73.39$ & $27.3$ & $90\,625$\\
        $128$ & $53.92$ & $37.2$ & $45\,312$\\
        $256$ & $26.07$ & $76.9$ & $22\,656$\\
        $512$ & $18.77$ & $106.7$ & $11\,328$\\
        $1024$ & $16.35$ & $122.5$ & $5\,664$\\
    \end{tabular}
    \label{tab:parallelTab}
\end{table}

\noindent Parallel speed-up $\frac{t_{1}}{t_{N}}$ refers to the improvement in execution time achieved by using multiple CPU cores ($t_{N}$) compared to a single core ($t_{1}$). It quantifies how effectively a parallelised task reduces computation time as more cores are added. Strong scaling, on the other hand, measures how the execution time of a fixed-size problem changes as the number of CPU cores increases. It helps evaluate the efficiency of a parallel algorithm when the problem size remains constant, highlighting how well the workload is distributed across increasing resources. Ideally, as the number of CPU cores increases, the execution time should scale linearly; i.e., for the same problem, using 100 CPU cores should result in an execution time that is 1\% (or $0.01$) of the time required when using a single CPU core. 

\begin{figure}[!h]
    \centering
    \subfloat[Speed-up\label{fig:speedup}]{\includegraphics[width=0.45\textwidth]{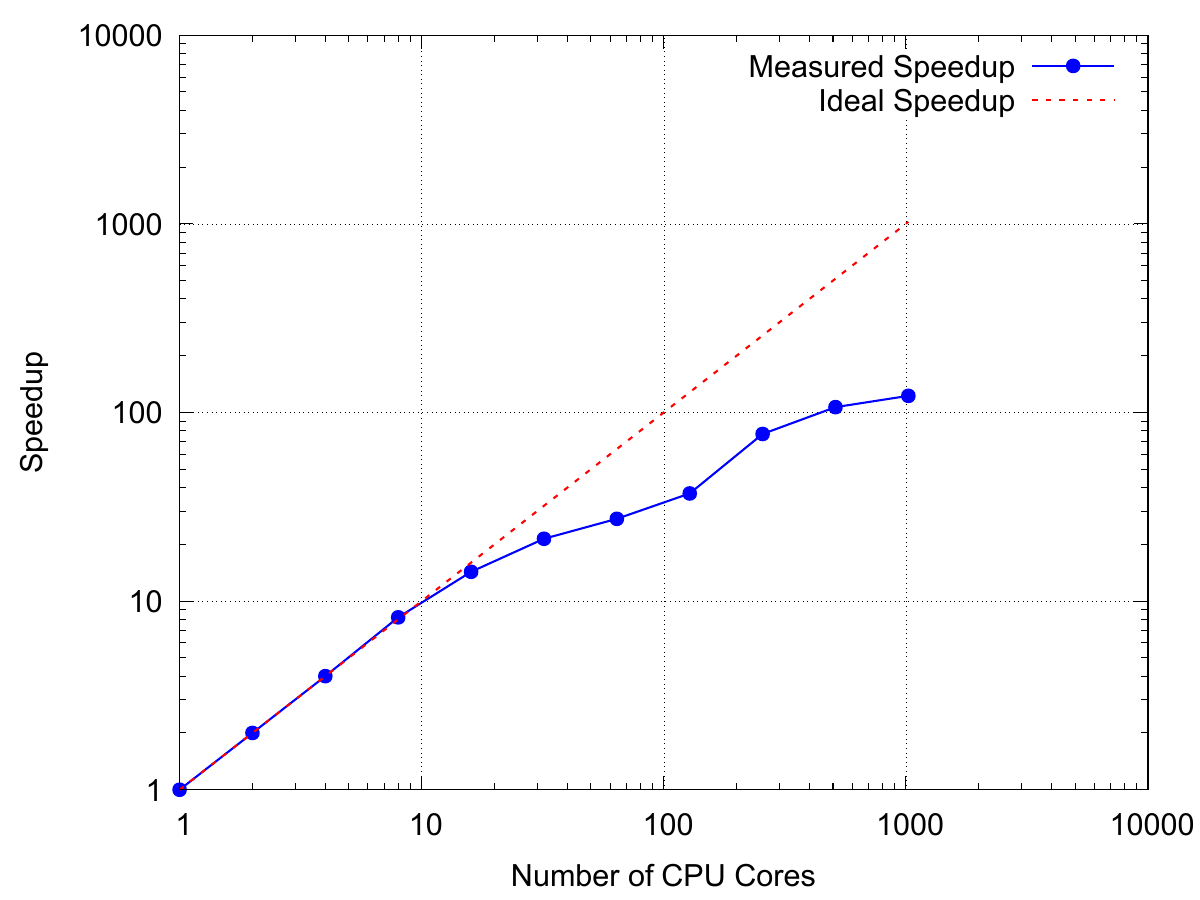}} \hspace{1em}
    \subfloat[Strong Scaling\label{fig:strongScaling}]{\includegraphics[width=0.45\textwidth]{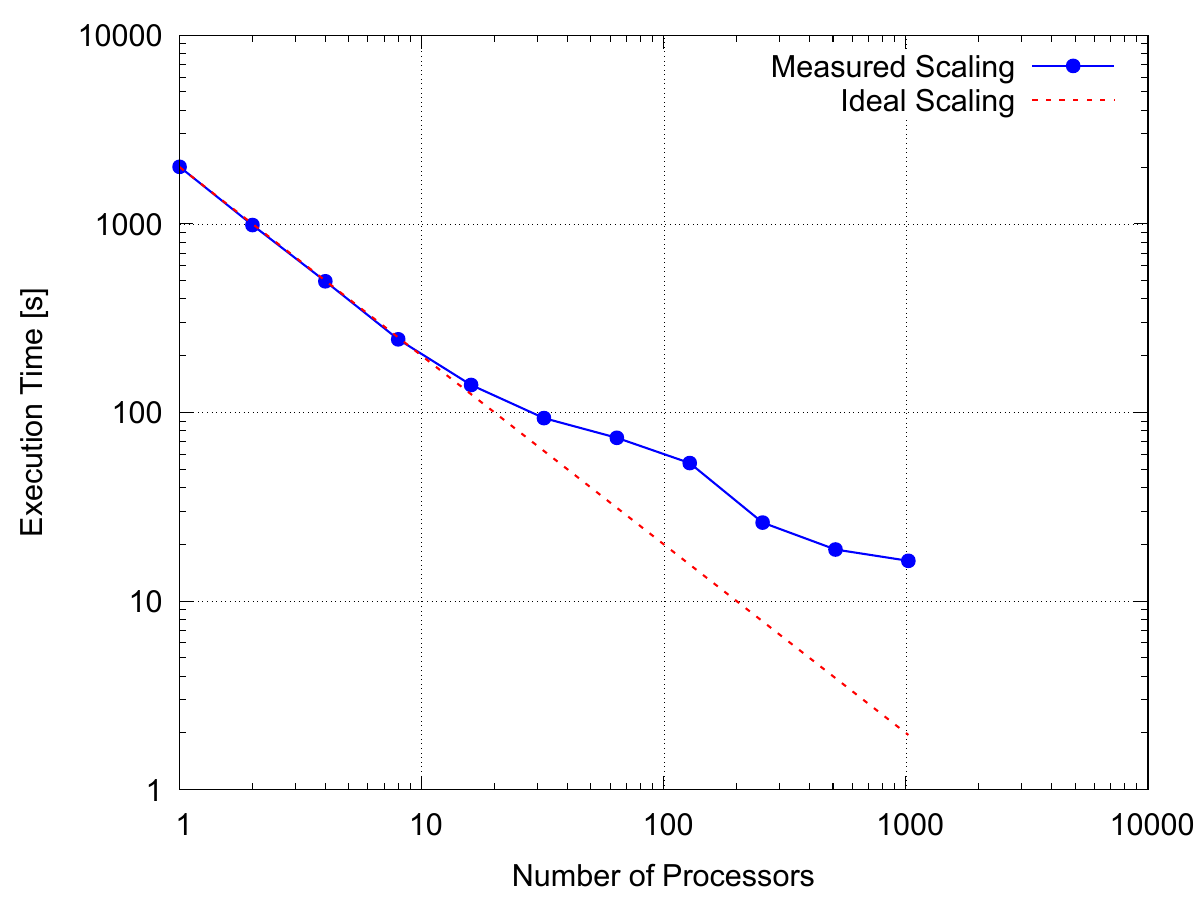}} \\[1em]
    \caption{Parallel performance scalability}
    \label{fig:parallelStudy}
\end{figure}

\noindent The results of the parallel performance study are illustrated in Figure \ref{fig:parallelStudy}. Based on Figure \ref{fig:speedup},  initially, as the number of cores increases, the speed-up grows almost linearly, with notable improvements in execution time. However, the speed-up starts to diminish as more cores are added, particularly after $128$ cores, where the rate of performance improvement slows down. By $512$ cores, the speed-up is about $106.72$, and with $1024$ cores, it reaches $122.51$, showing diminishing returns. This suggests that the method’s parallel efficiency decreases as the number of cores increases, likely due to communication overhead and decreasing workload per core. As the number of cores increases, the execution time decreases, but the reduction slows down significantly after $128$ cores. 

\noindent Figure \ref{fig:strongScaling} shows that initially, the algorithm scales efficiently, with notable reductions in execution time. However, after $256$ cores, the performance gains become less significant, indicating diminishing returns. By the time $1024$ cores are used, the execution time is only slightly improved compared to $512$ cores, highlighting the impact of communication overhead and load imbalance at higher core counts.
\clearpage

\begin{figure}[!h]
     \centering
     \begin{subfigure}[b]{\textwidth}
         \centering
         \includegraphics[width=\textwidth]{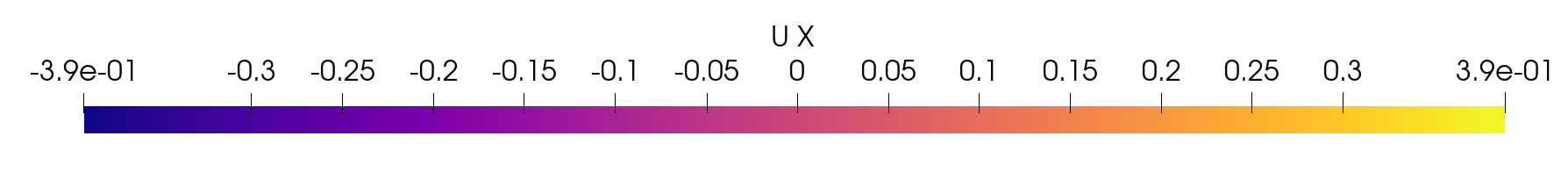}
     \end{subfigure}
     \begin{subfigure}[b]{0.45\textwidth}
         \centering
         \includegraphics[width=\textwidth]{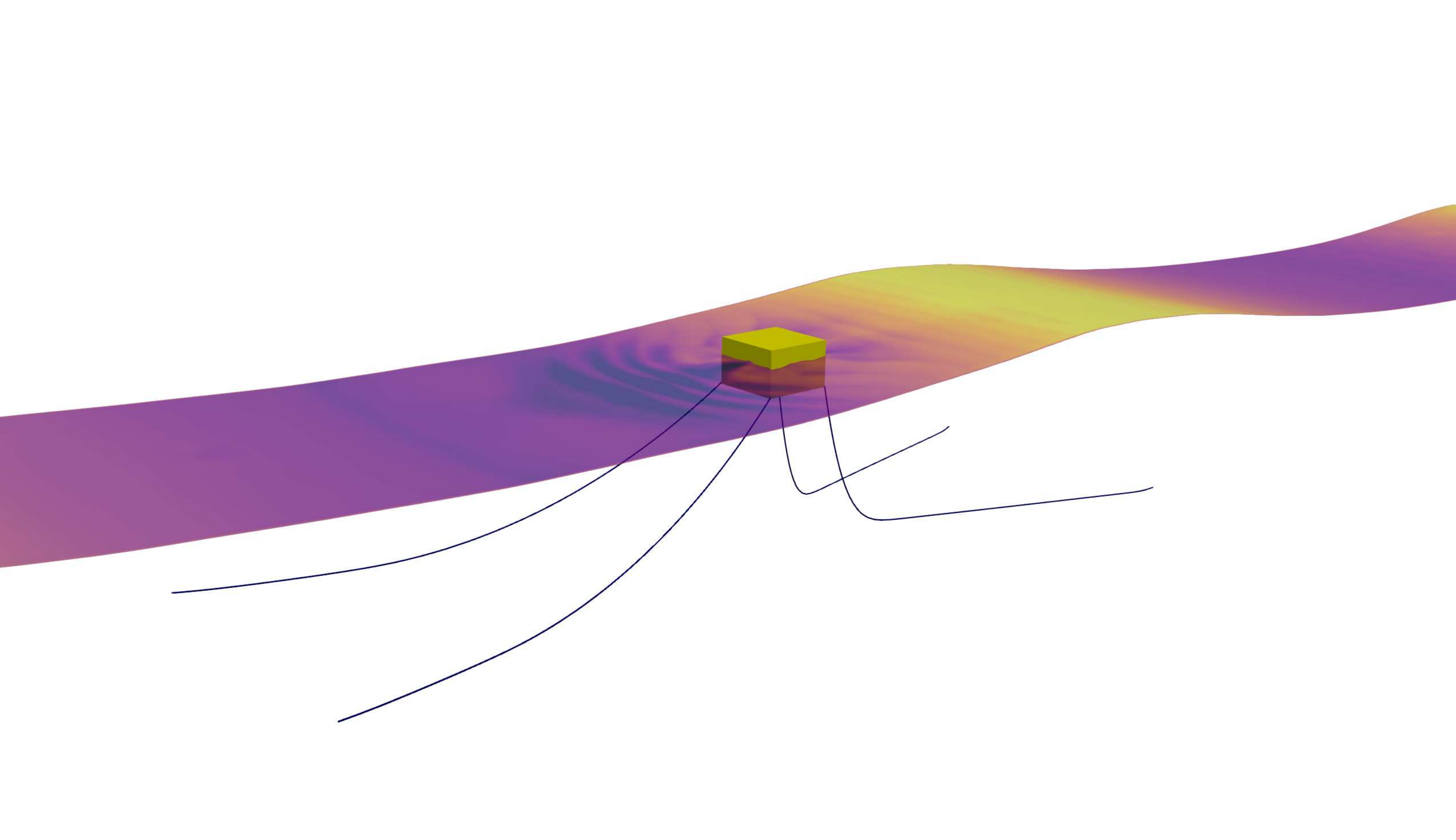}
         \caption{T = 10 s}
     \end{subfigure}
     \hfill
     \begin{subfigure}[b]{0.45\textwidth}
         \centering
         \includegraphics[width=\textwidth]{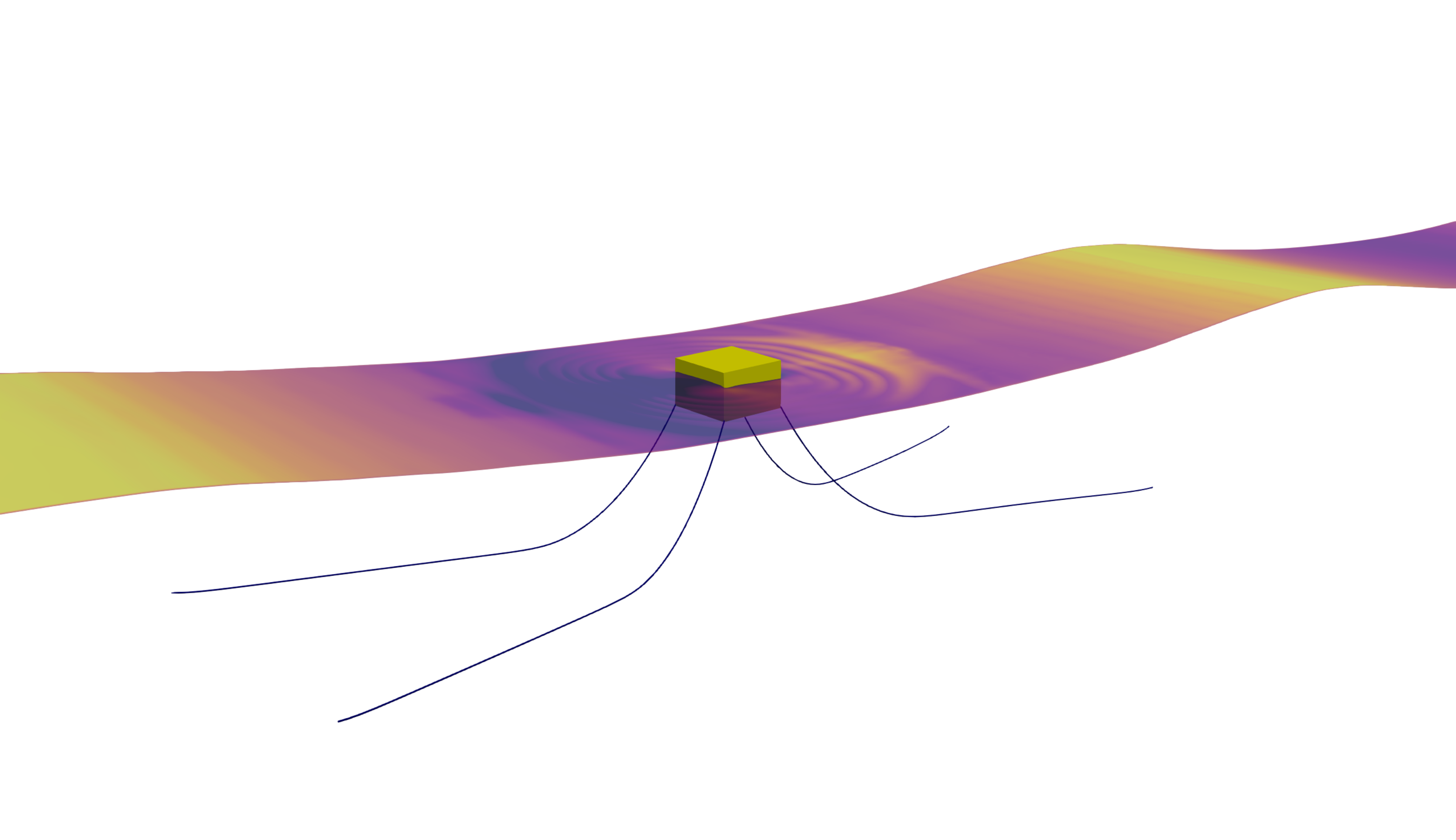}
         \caption{T = 10.5 s}
     \end{subfigure}
     \hfill
     \begin{subfigure}[b]{0.45\textwidth}
         \centering
         \includegraphics[width=\textwidth]{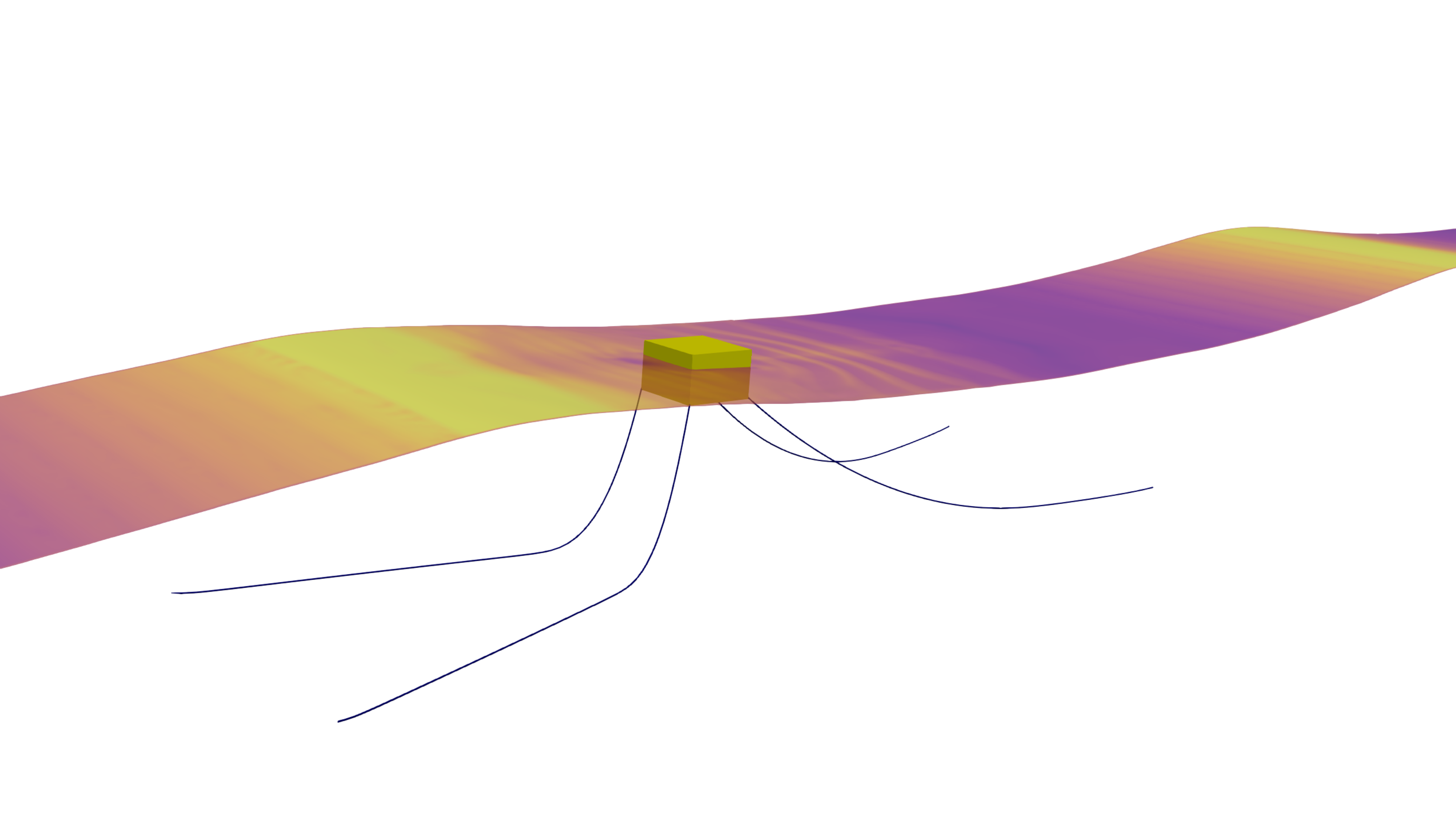}
         \caption{T = 11 s}
     \end{subfigure}
     \hfill
     \begin{subfigure}[b]{0.45\textwidth}
         \centering
         \includegraphics[width=\textwidth]{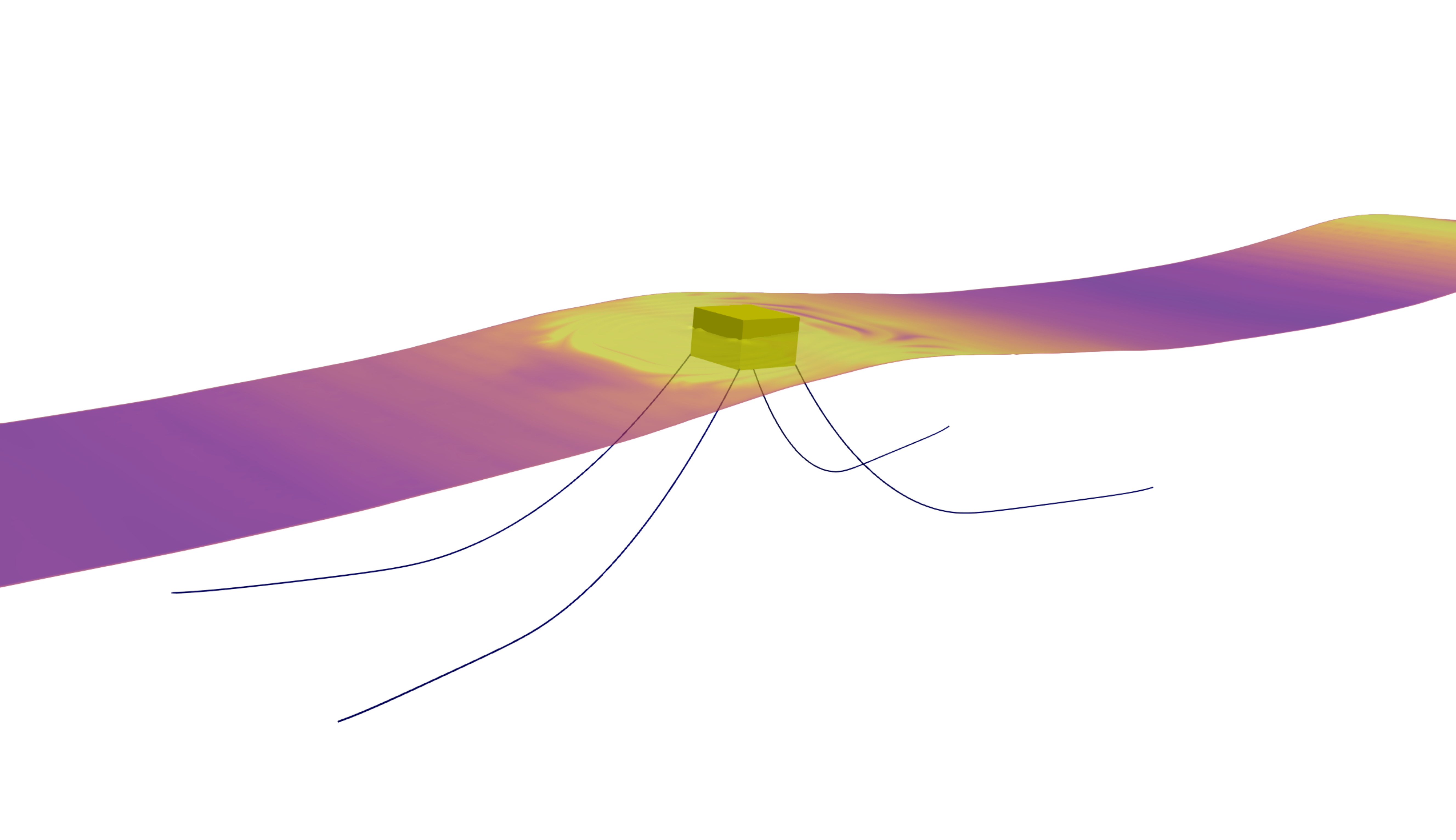}
         \caption{T = 11.5 s}
     \end{subfigure}
     \hfill
     \begin{subfigure}[b]{0.45\textwidth}
         \centering
         \includegraphics[width=\textwidth]{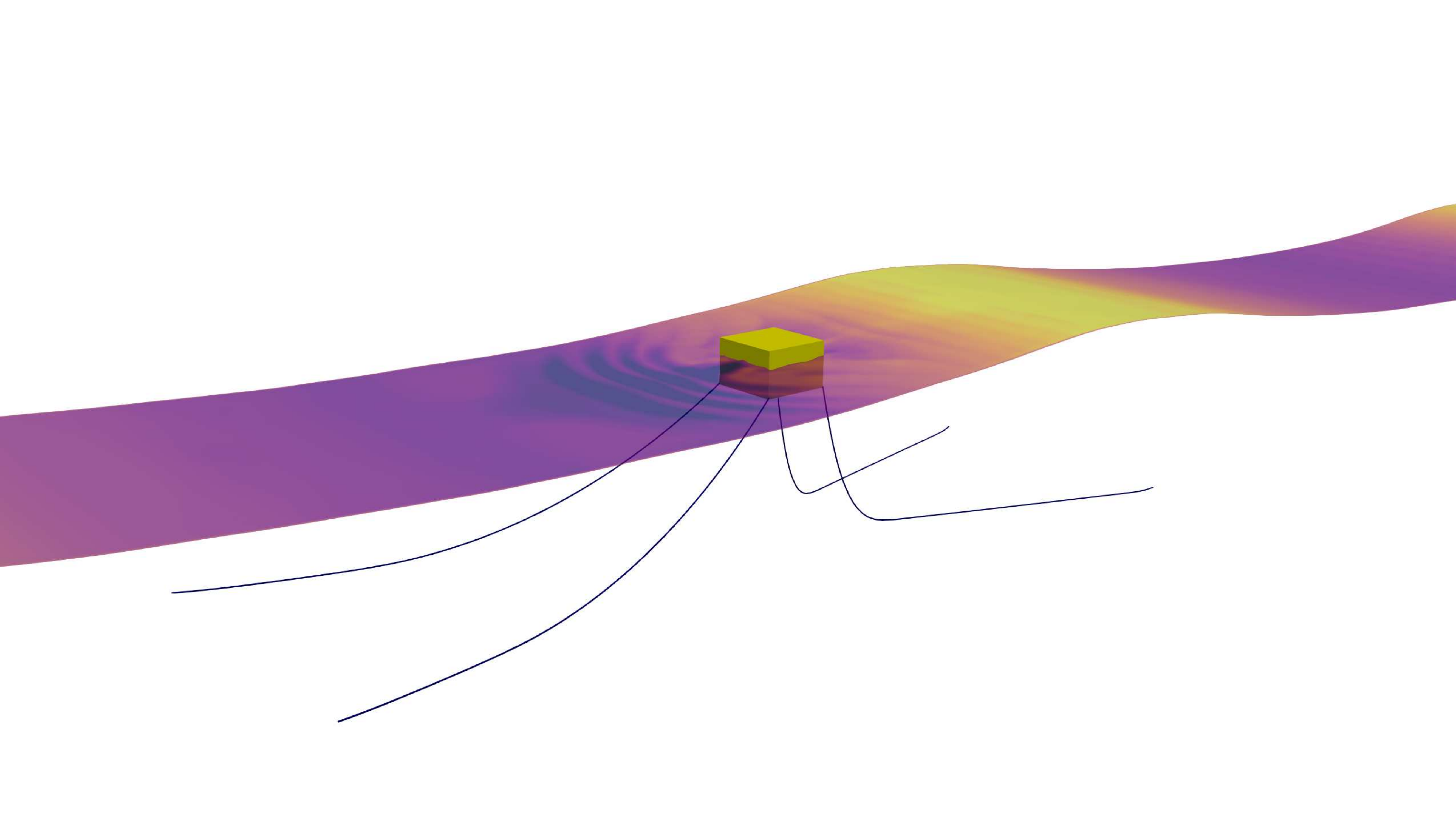}
         \caption{T = 12 s}
     \end{subfigure}
     \hfill
     \begin{subfigure}[b]{0.45\textwidth}
         \centering
         \includegraphics[width=\textwidth]{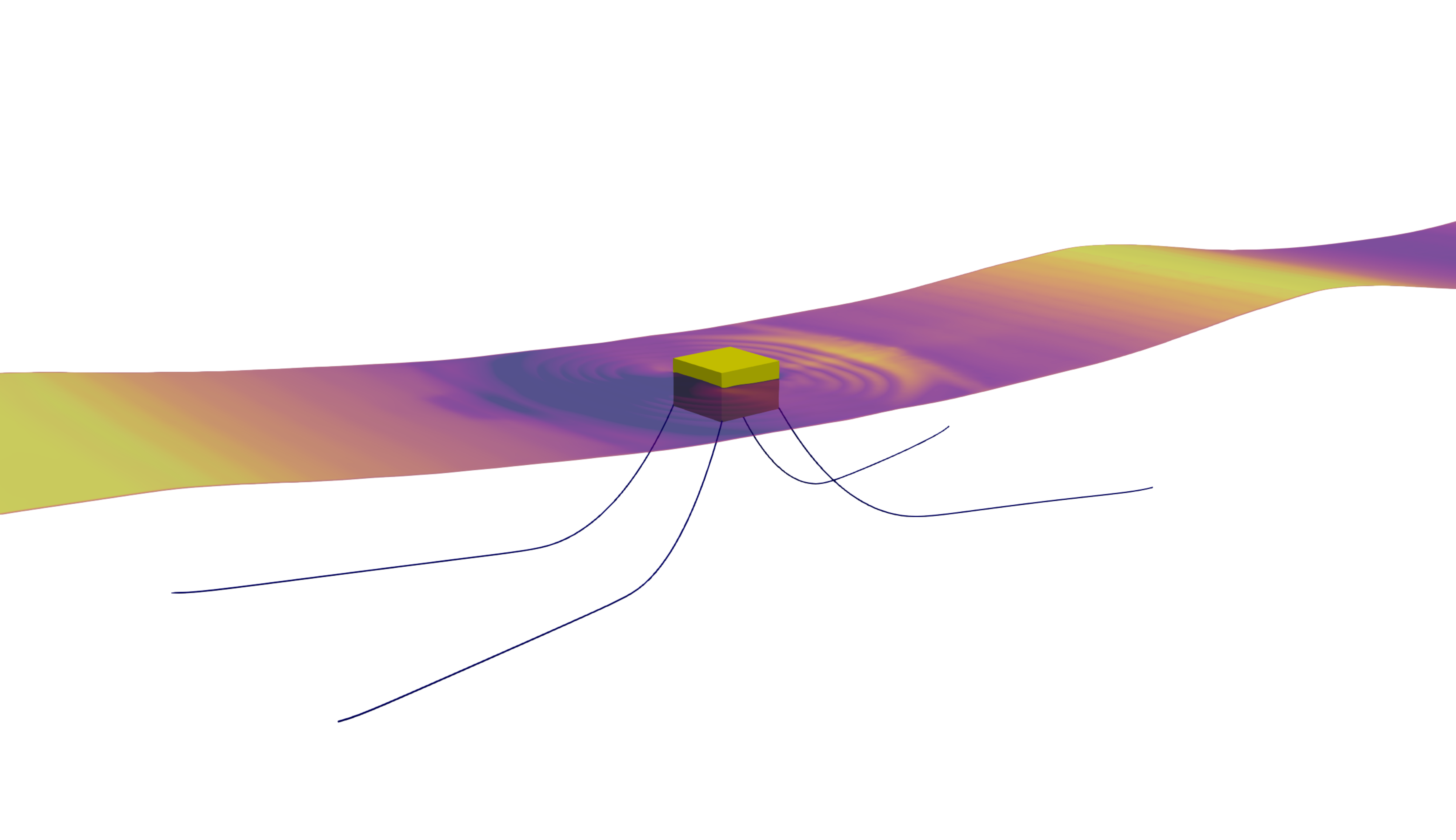}
         \caption{T = 12.5 s}
     \end{subfigure}
     \hfill
        \caption{Temporal snapshots capturing wave-box-mooring interactions during multiple wave cycle for H15T18, wave elevation iso-surface ($\alpha = 0.5$) coloured by flow horizontal velocity($U_{x}$)}
        \label{fig:SnapShots}
\end{figure}

\noindent Figure \ref{fig:SnapShots} shows the free surface interaction with incoming waves and mooring lines for the case H12T20. The iso-surface for $\alpha = 0.5$ represents the free surface and is coloured by the horizontal flow velocity. Each snapshot represents a successive instance of incoming waves so that an entire cycle of surge and heave motions can be observed.

\noindent The numerical and experimental wave profiles at the locations listed in Table \ref{table:wavegaugescoords} for both the H12T20 and H15T18 cases are presented in Figure \ref{fig:waveProfiles}. While the amplitude agreement between the numerical and experimental data is generally good, an evident phase lag was initially observed. This lag likely arises from differences in how waves are generated in the two setups. In the experiments, the wave maker may take some time to reach the target amplitude and frequency, and initial disturbances can affect the wave shape. In contrast, the numerical simulations assume an ideal and instantaneous start to wave generation, which leads to a timing mismatch. To account for this and allow a clearer comparison between the wave profiles, a constant time shift was applied to align the numerical and experimental data. The shift was determined by matching the first wave peak and was applied uniformly to all wave gauges within each test case (i.e., a single shift for H12T20 and another for H15T18). This adjustment makes it easier to identify differences in wave shape and amplitude.
\begin{figure}[!h]
    \centering
    \subfloat[H12T20 - wave gauge 2]{\includegraphics[width=0.45\textwidth]{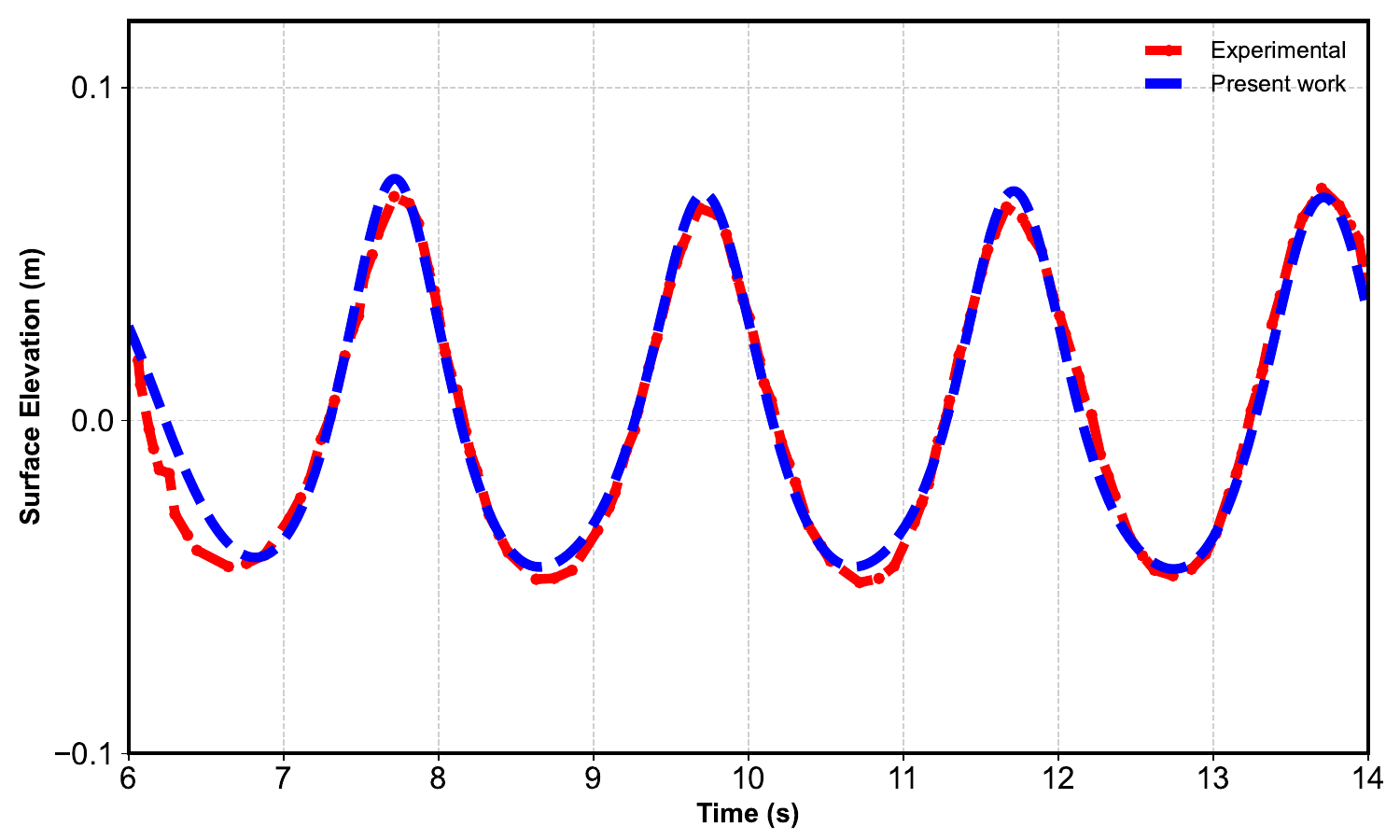}} \hspace{1em}
    \subfloat[H15T18 - wave gauge 2]{\includegraphics[width=0.45\textwidth]{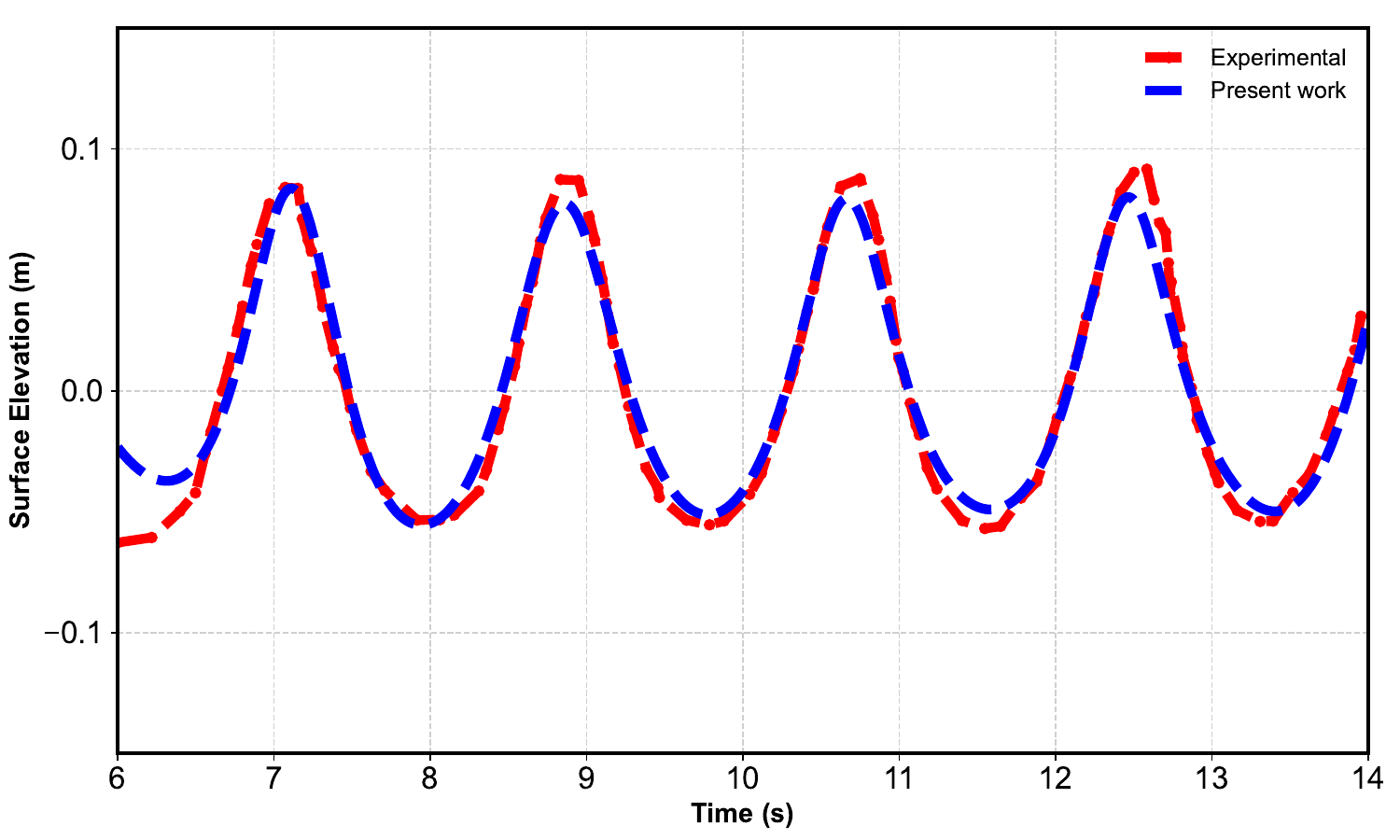}} \\[1em]
    
    \subfloat[H12T20 - wave gauge 3]{\includegraphics[width=0.45\textwidth]{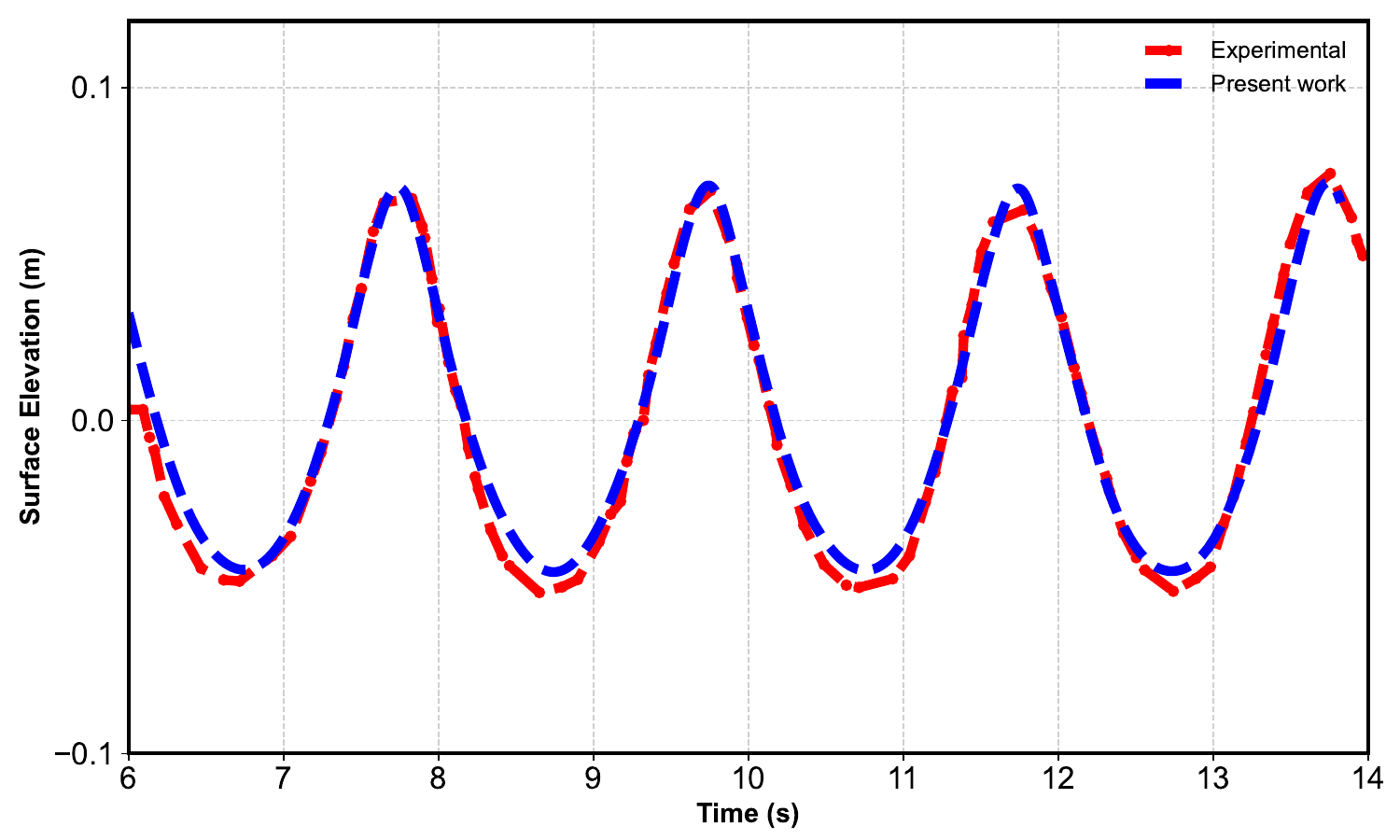}} \hspace{1em}
    \subfloat[H15T18 - wave gauge 3]{\includegraphics[width=0.45\textwidth]{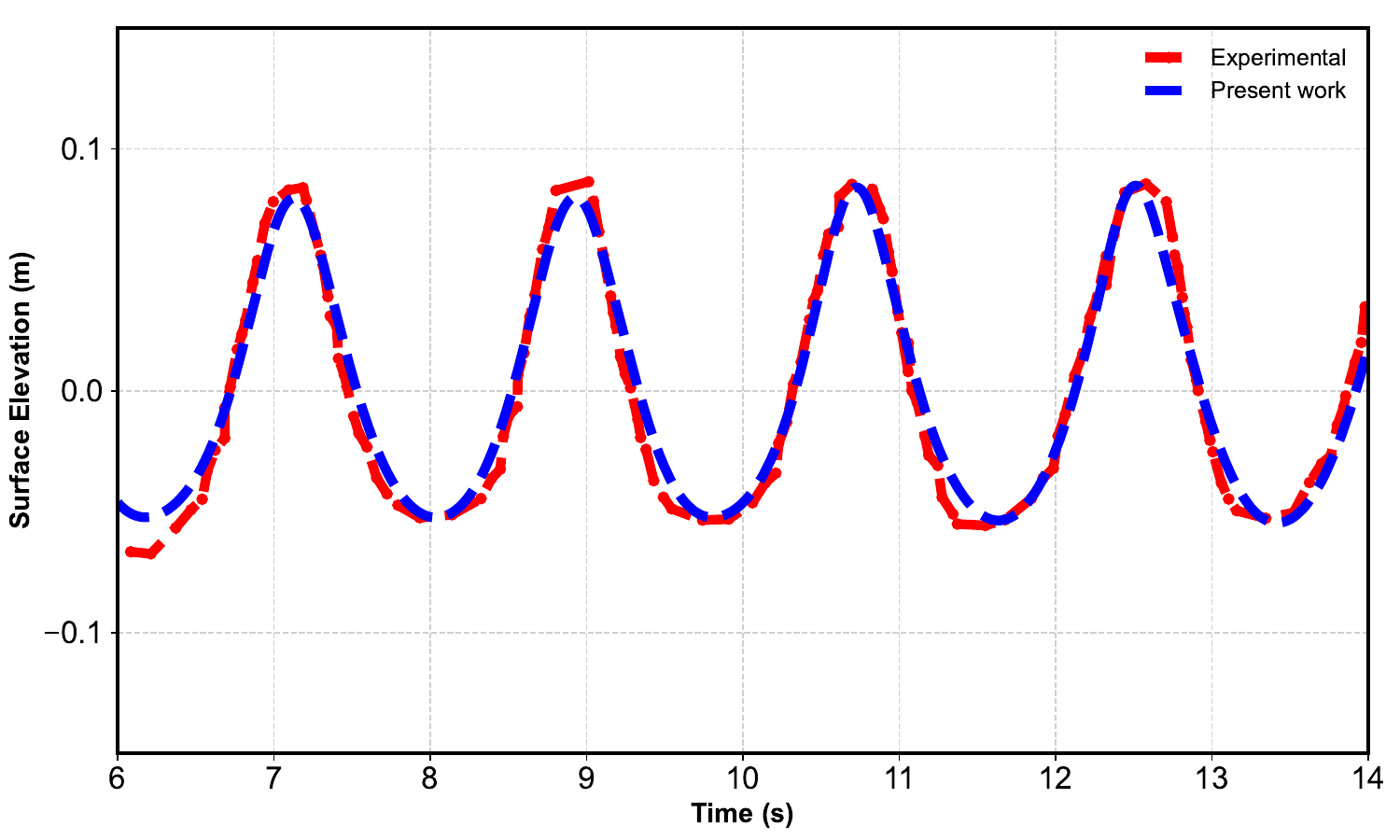}} \\[1em]
\end{figure}
\begin{figure}[!h]
    \centering
    \ContinuedFloat
    \subfloat[H12T20 - wave gauge 4]{\includegraphics[width=0.45\textwidth]{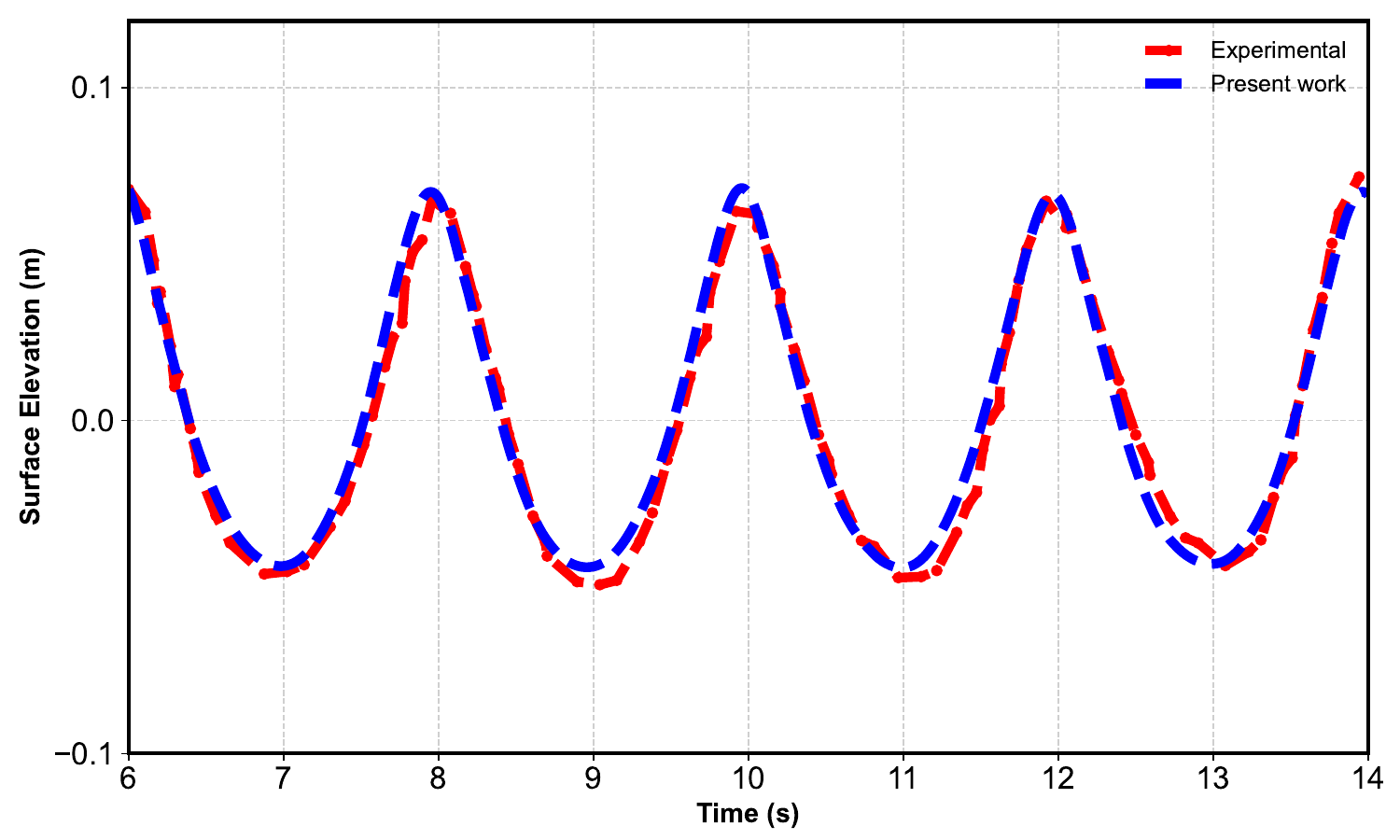} \label{fig:wg412}} \hspace{1em}
    \subfloat[H15T18 - wave gauge 4]{\includegraphics[width=0.45\textwidth]{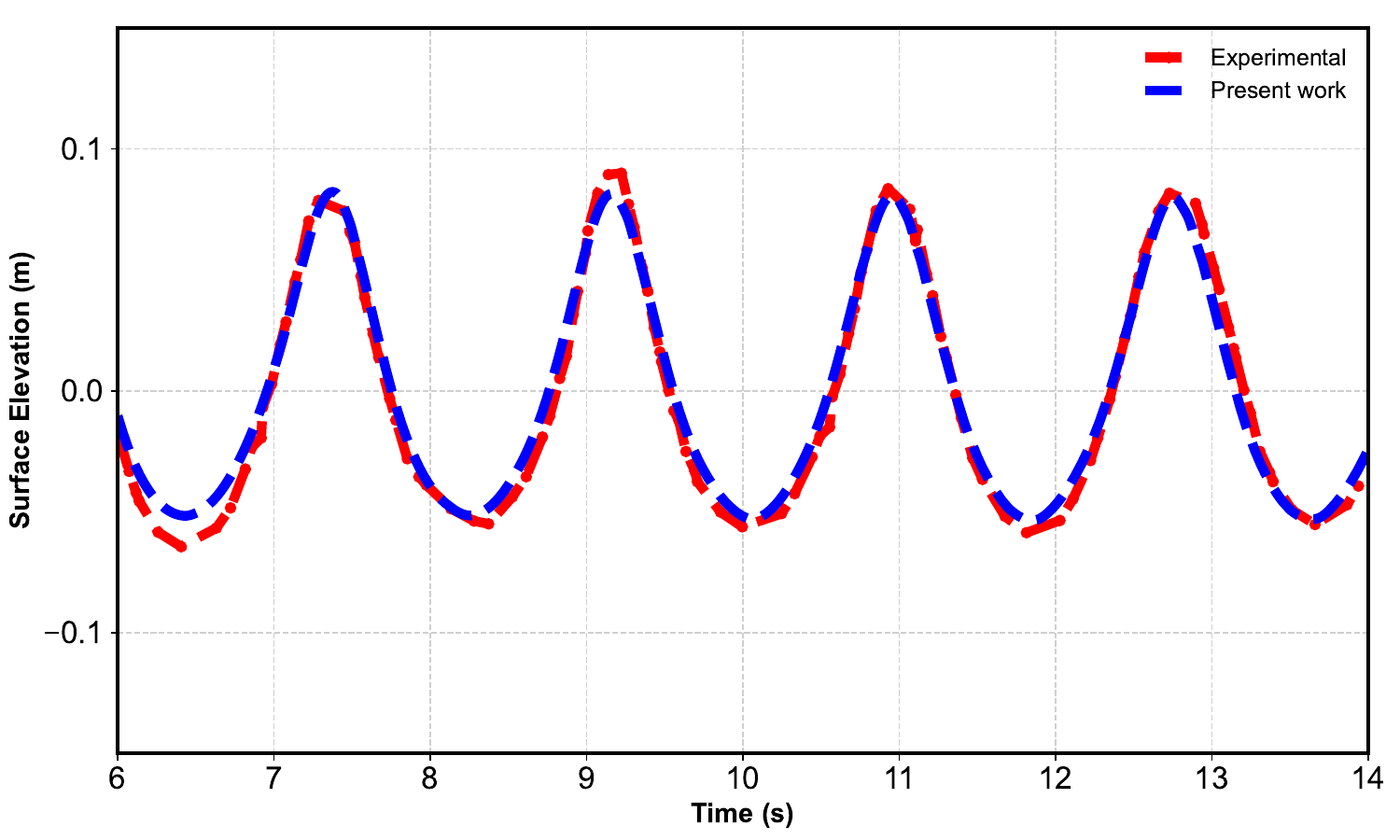} \label{fig:wg415}}\\[1em]
    
    \subfloat[H12T20 - wave gauge 5]{\includegraphics[width=0.45\textwidth]{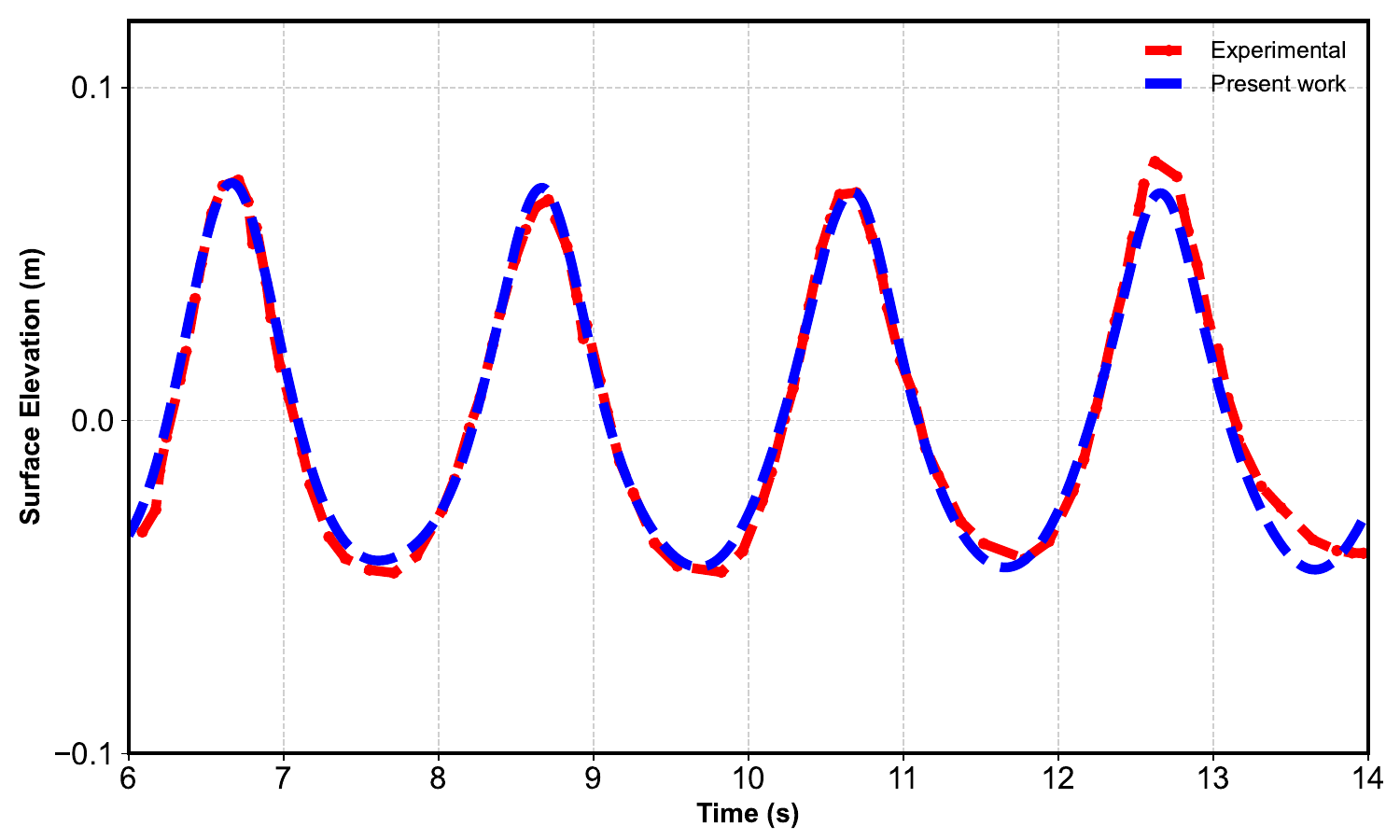} \label{fig:wg512}} \hspace{1em}
    \subfloat[H15T18 - wave gauge 5]{\includegraphics[width=0.45\textwidth]{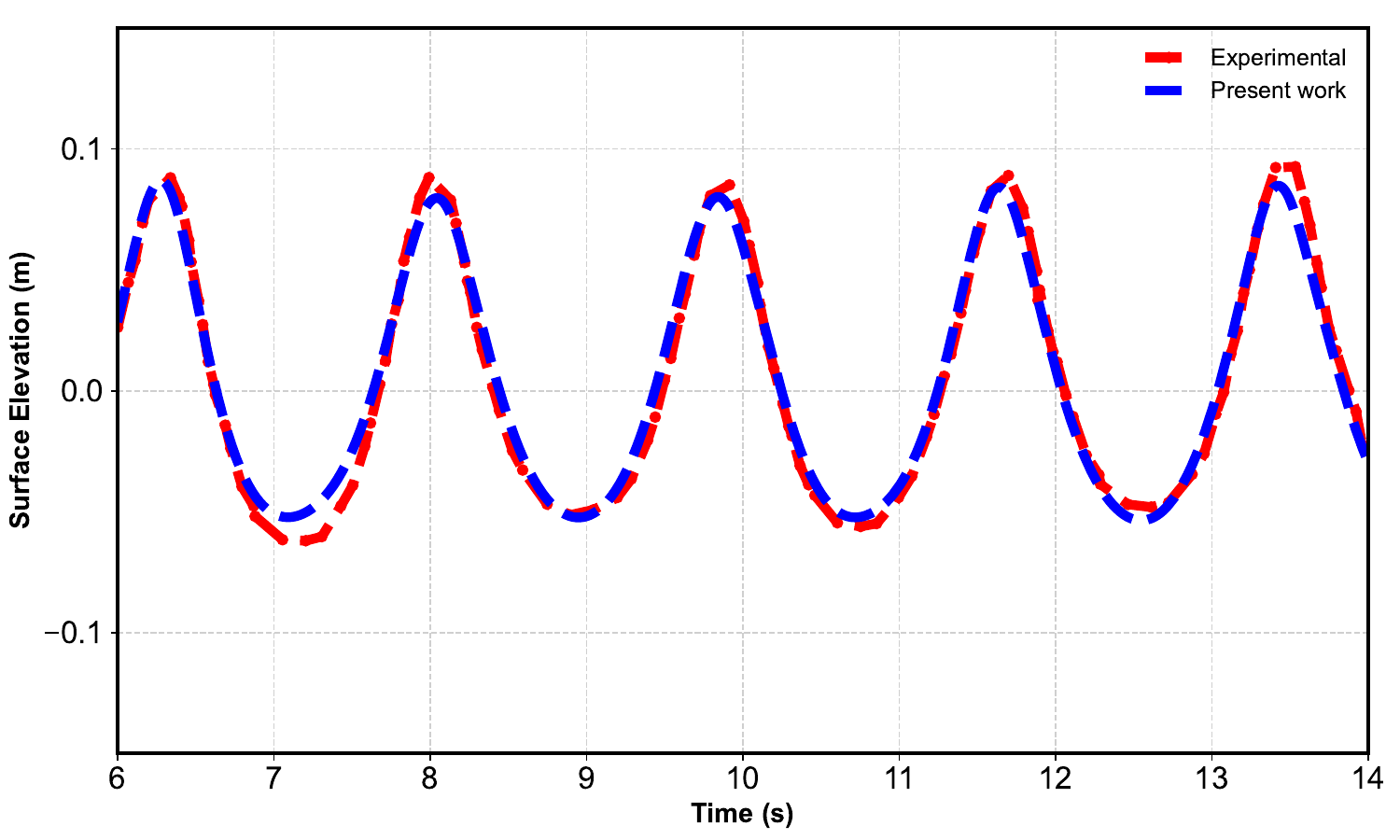} \label{fig:wg515}} \\[1em]
    
    \subfloat[H12T20 - wave gauge 6]{\includegraphics[width=0.45\textwidth]{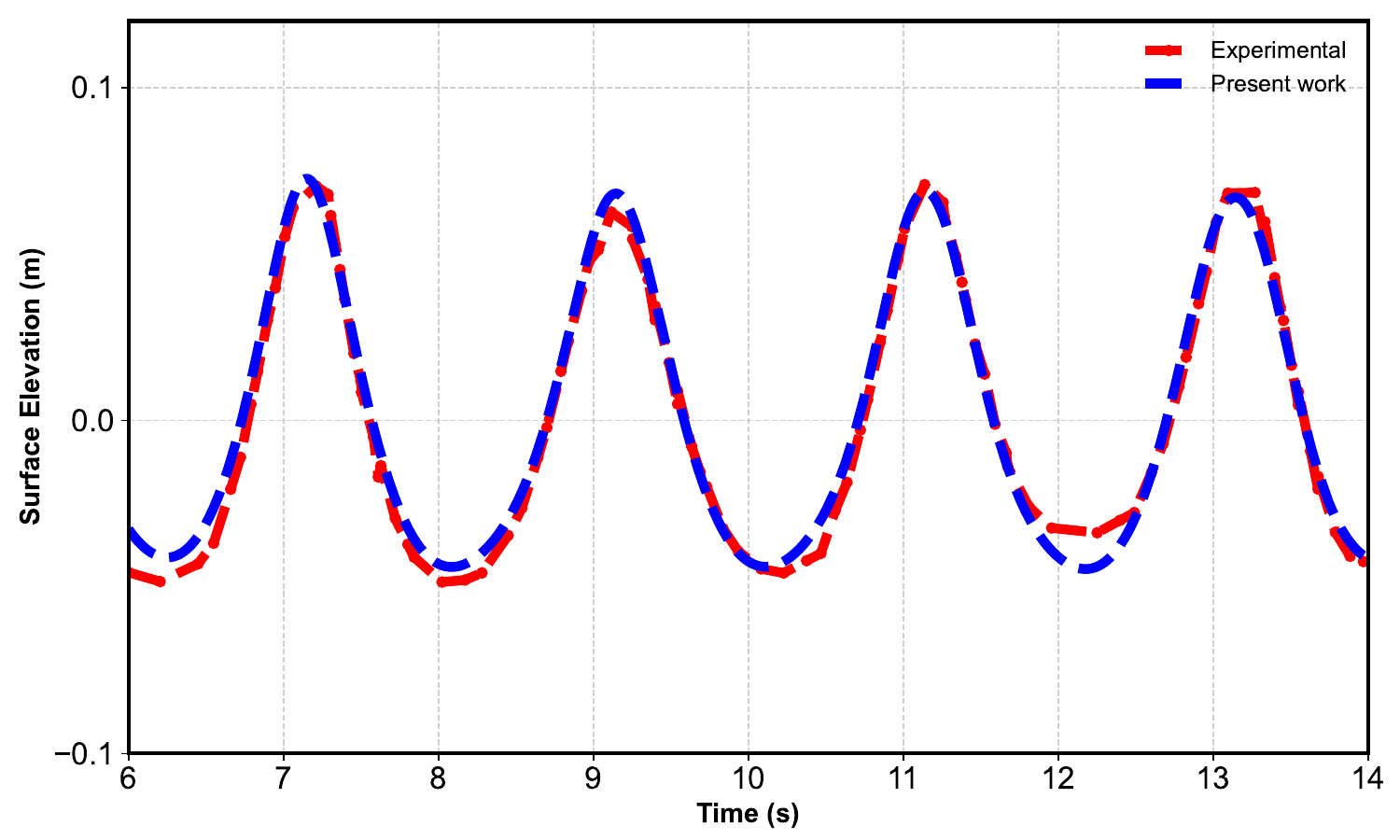} \label{fig:wg612}} \hspace{1em}
    \subfloat[H15T18 - wave gauge 6]{\includegraphics[width=0.45\textwidth]{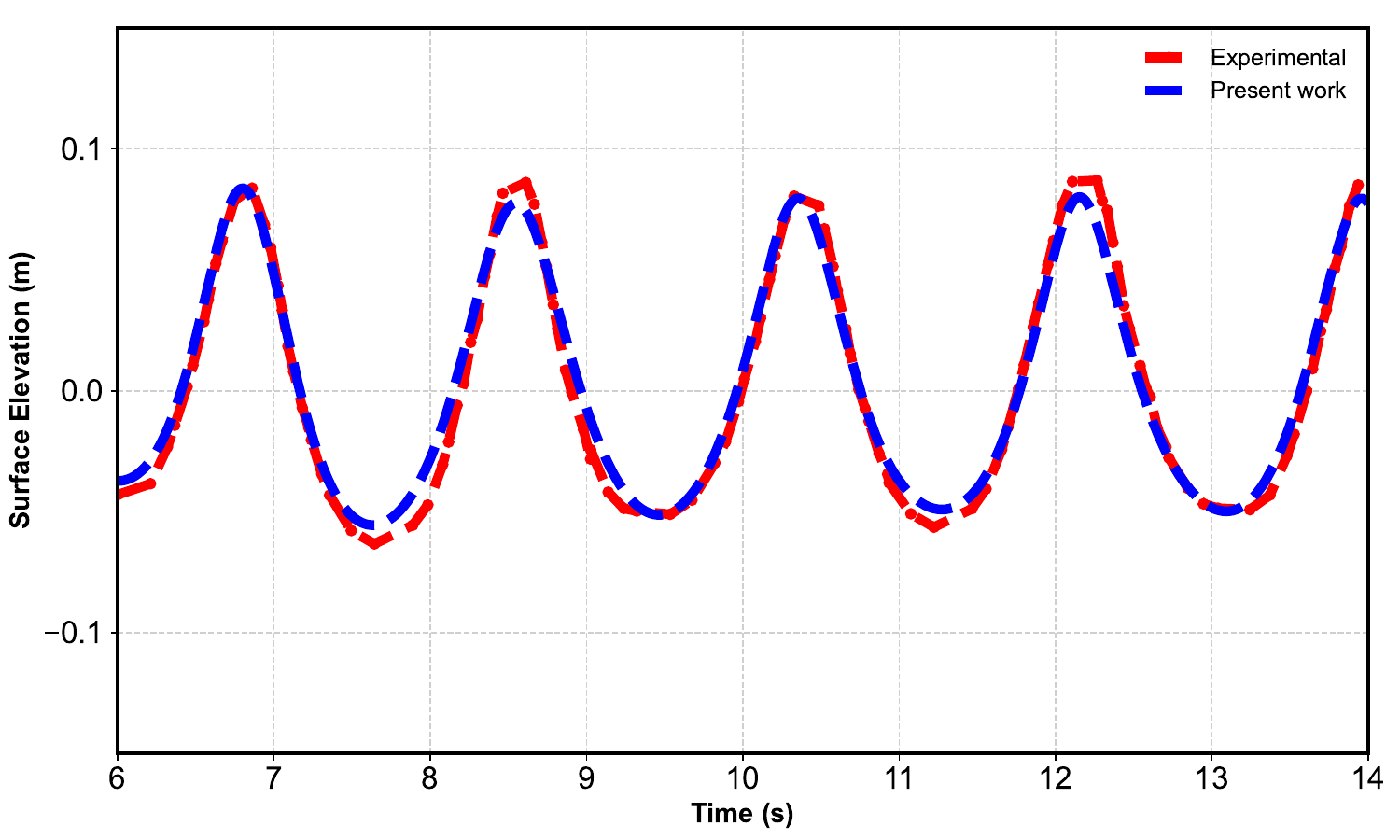} \label{fig:wg615}}
    
    \caption{Wave height profiles on wave gauge locations}
    \label{fig:waveProfiles}
\end{figure}
\clearpage

\noindent Table \ref{tab:wave_amplitudes} provides a quantitative comparison of wave height amplitudes between experimental measurements and the present numerical model. These values are extracted using FFT analysis over the interval [6\,\text{s}, 14\,\text{s}], where the amplitude is calculated by identifying the dominant frequency in the spectrum and taking the corresponding magnitude, which is corrected for windowing effects and the length of the signal. The simulation shows strong agreement with the experimental amplitudes, particularly near the floating body, where the mesh resolution is finer. At wave gauges 2 and 3, the differences between the present model and the experimental results are minor. Further downstream, however, a gradual divergence appears at gauges 4–6 (Figures {\ref{fig:wg412}}, {\ref{fig:wg415}}, {\ref{fig:wg512}}, {\ref{fig:wg515}}, {\ref{fig:wg612}}, and {\ref{fig:wg615}}). This trend reflects the combined influence of weak residual currents produced by partial reflections at the glass walls, higher-order harmonics introduced by the hinged-flap wavemaker that are absent from the Stokes-II field imposed in the numerical relaxation zone, and mesh coarsening in the absorption region, which increases numerical dispersion and damps the higher-frequency components of the wave train.
Among these possible sources, the dominant factor is the gradual divergence at the downstream gauges, which is mainly caused by residual reflections in the experimental wave basin and numerical dispersion introduced by mesh coarsening in the absorption zone. These effects outweigh the smaller contributions from wavemaker startup transients and higher-order harmonics.
Overall, the average amplitude difference compared to the experimental values is 0.86\% for H12T20, with the largest deviation of 2.27\% observed at wave gauge 6. For H15T18, the average difference increases slightly to 1.30\%, with a maximum deviation of 4.0\% also occurring at wave gauge 6. Similar discrepancies have been reported in previous studies using this benchmark. For example, \citet{chen_cfd_2022} and \citet{dominguez_sph_2019} observed similar amplitude deviations in downstream gauges. This highlights common challenges in maintaining wave accuracy over distance rather than limitations specific to the present model.

\begin{figure}[!h]
    \centering
    \subfloat[H12T20 - wave gauge 2]{\includegraphics[width=0.45\textwidth]{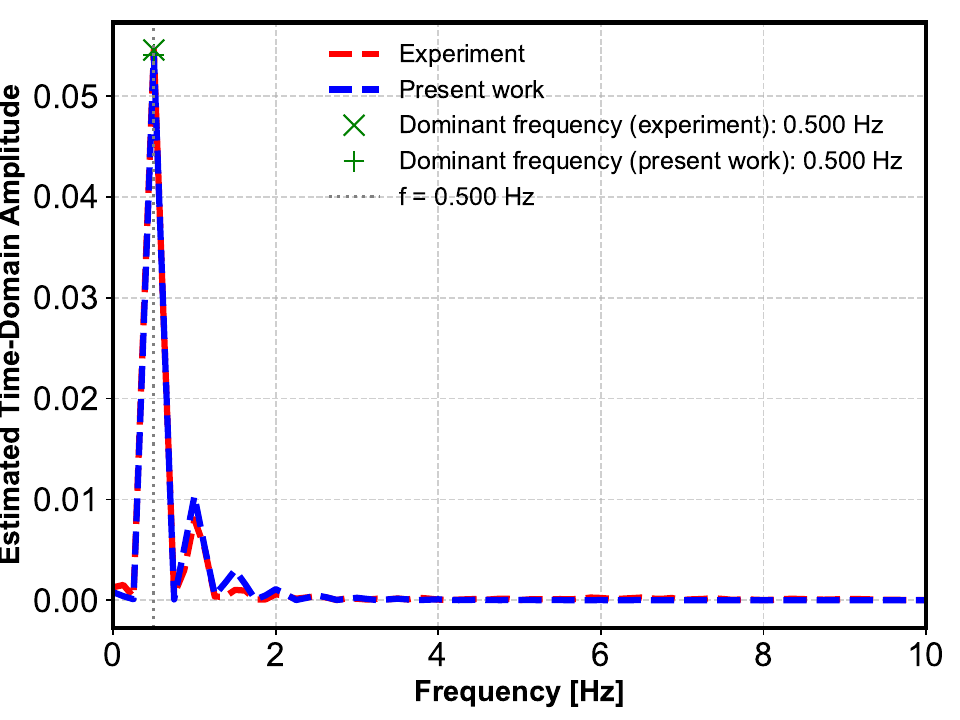}} \hspace{1em}
    \subfloat[H15T18 - wave gauge 2]{\includegraphics[width=0.45\textwidth]{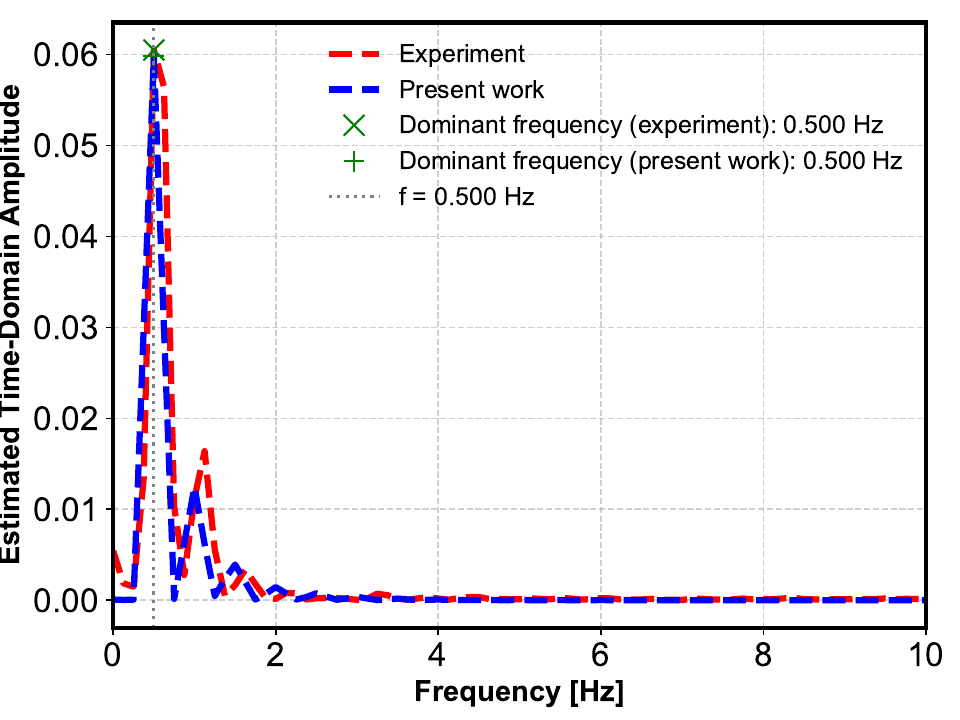}} \\[1em]
    
    \subfloat[H12T20 - wave gauge 6]{\includegraphics[width=0.45\textwidth]{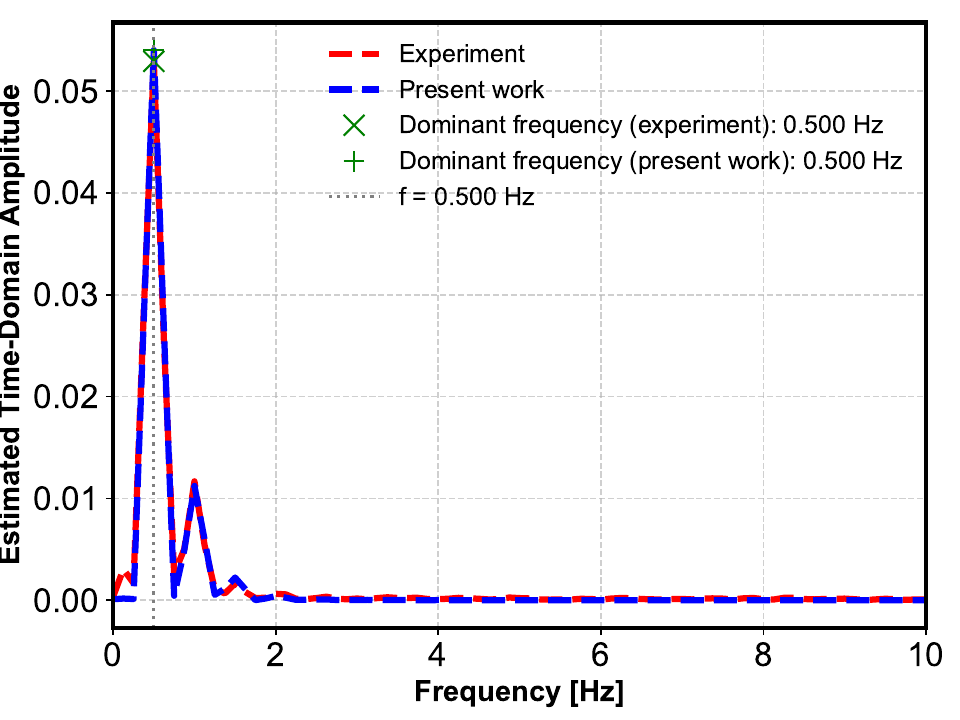}} \hspace{1em}
    \subfloat[H15T18 - wave gauge 6]{\includegraphics[width=0.45\textwidth]{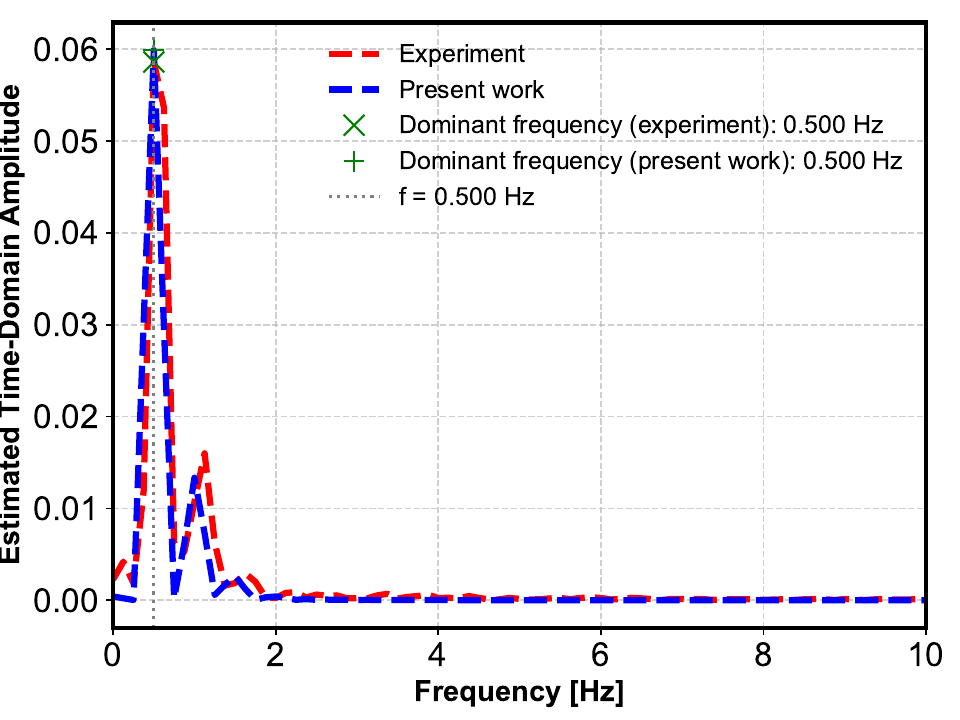}} \\[1em]
    \caption{FFT spectrum of surface elevation for H12T20 and H15T18 cases at wave gauge 2 and 6}
    \label{fig:fftofWaves}
\end{figure}

\begin{table}[!h]
    \centering
    \small
    \caption{Wave height amplitude comparison at wave gauges for H12T20 and H15T18}
    \begin{tabularx}{\linewidth}{l*{4}{>{\centering\arraybackslash}X}}
        \toprule
        \multirow{2}{*}{\textbf{Wave Gauge}} & \multicolumn{2}{c}{\textbf{H12T20}} & \multicolumn{2}{c}{\textbf{H15T18}} \\
        \cmidrule{2-5}
        & Experimental (m) & Present model (m) & Experimental (m) & Present model (m) \\
        \midrule
        Wave gauge 2 & 0.0546 & 0.0541 & 0.0604 & 0.0601 \\
        Wave gauge 3 & 0.0539 & 0.0535 & 0.0605 & 0.0599 \\
        Wave gauge 4 & 0.0534 & 0.0533 & 0.0591 & 0.0588 \\
        Wave gauge 5 & 0.0533 & 0.0534 & 0.0593 & 0.0590 \\
        Wave gauge 6 & 0.0528 & 0.0540 & 0.0576 & 0.0599 \\
        \midrule
        \midrule
        \textbf{Average diff. (\%)} & \multicolumn{1}{c}{--} & 0.86 & \multicolumn{1}{c}{--} & 1.30 \\
        \textbf{Max diff. (\%)} & \multicolumn{1}{c}{--} & 2.27 & \multicolumn{1}{c}{--} & 4.00 \\
        \bottomrule
    \end{tabularx}
    \label{tab:wave_amplitudes}
\end{table}
\clearpage
\noindent The FFT spectrum is presented in Figure \ref{fig:fftofWaves} for both H12T20 and H15T18 cases. These plots show that the dominant wave frequency observed in the experimental data is also predicted by the numerical model, indicating that the primary wave characteristics are well reproduced. Among all the wave gauge locations, only wave gauges 2 and 6 are shown here for brevity. These two locations were selected to highlight the performance of the model in regions with different mesh resolutions and wave conditions.

\clearpage
\subsection{Floating Box Motions}\label{subsec:boxmotions}
\begin{figure}[!ht]
     \centering
     \begin{subfigure}[b]{0.45\textwidth}
         \centering
         \includegraphics[width=\textwidth]{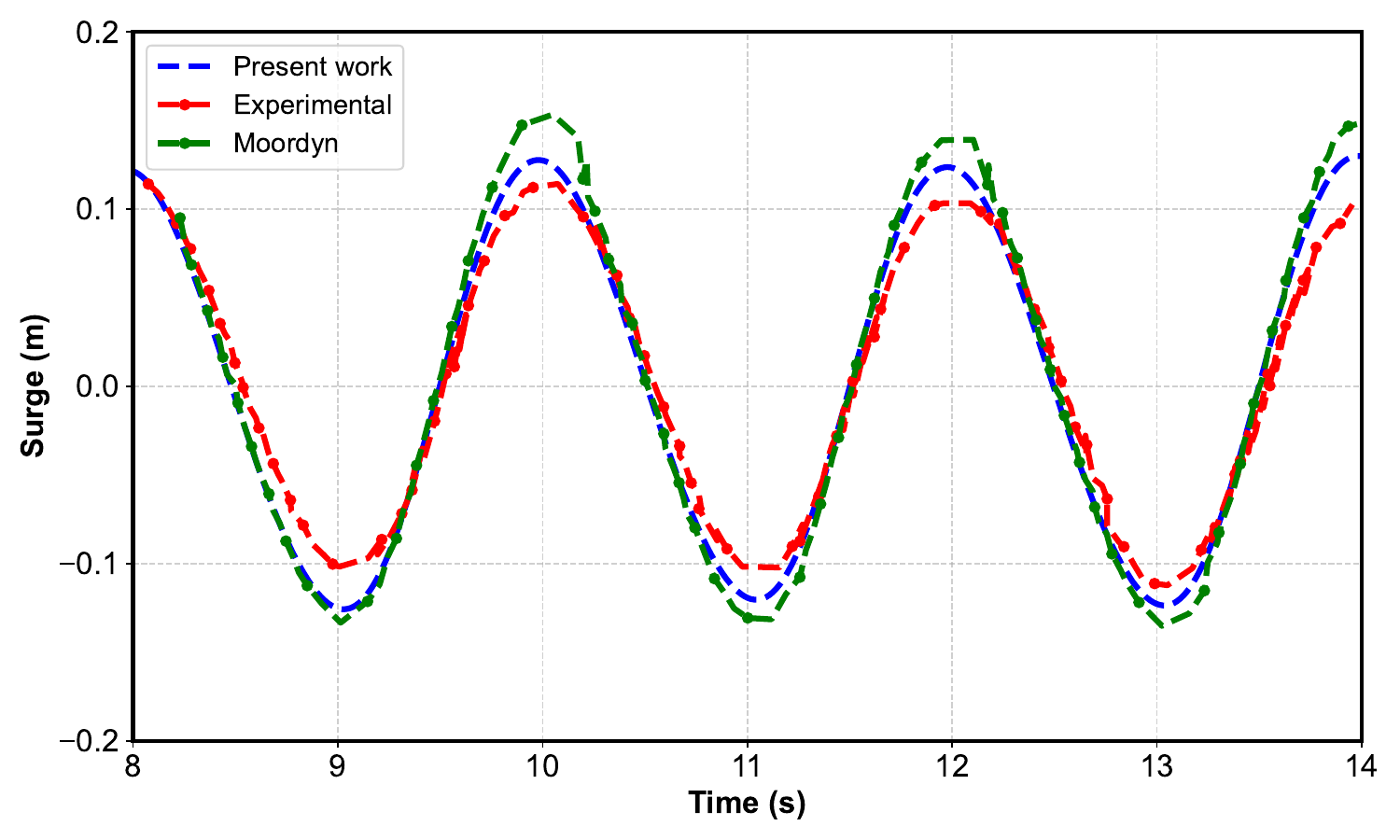}
         \caption{H12T20-Surge}
         \label{fig:H12T20-Surge}
     \end{subfigure}
     \hfill
     \begin{subfigure}[b]{0.45\textwidth}
         \centering
         \includegraphics[width=\textwidth]{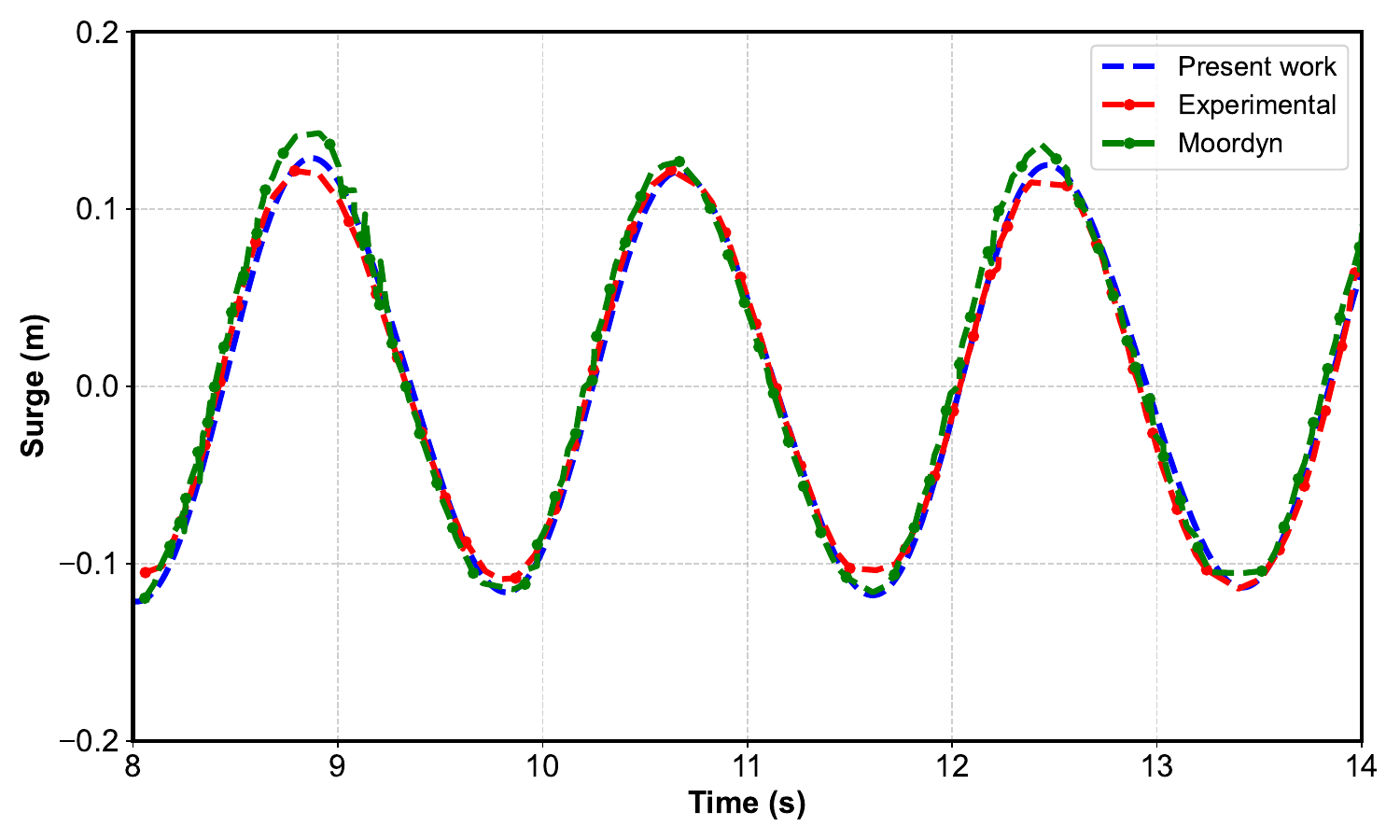}
         \caption{H15T18-Surge}
         \label{fig:H15T18-Surge}
     \end{subfigure}
     \hfill
     \begin{subfigure}[b]{0.45\textwidth}
         \centering
         \includegraphics[width=\textwidth]{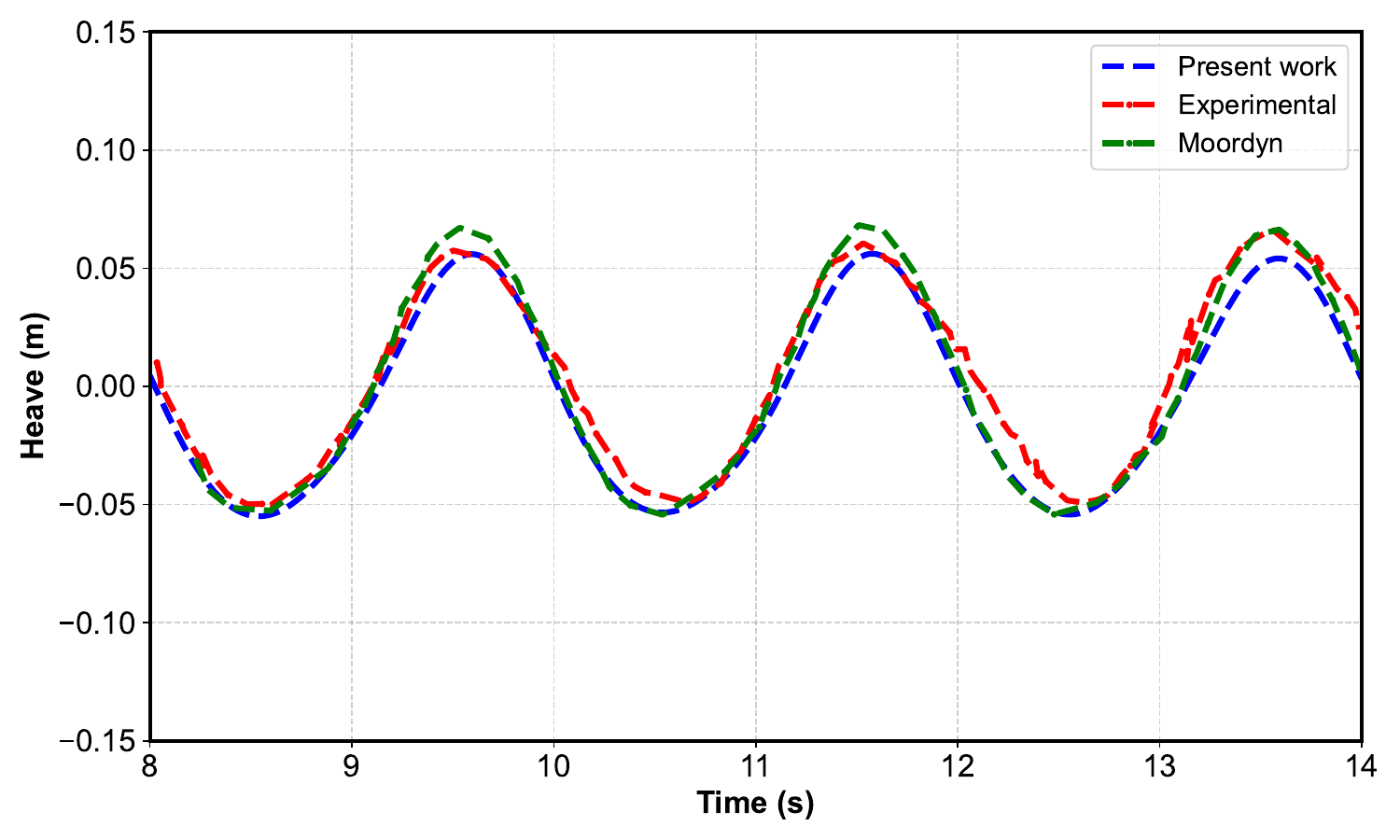}
         \caption{H12T20-Heave}
         \label{fig:H12T20-Heave}
     \end{subfigure}
     \hfill
     \begin{subfigure}[b]{0.45\textwidth}
         \centering
         \includegraphics[width=\textwidth]{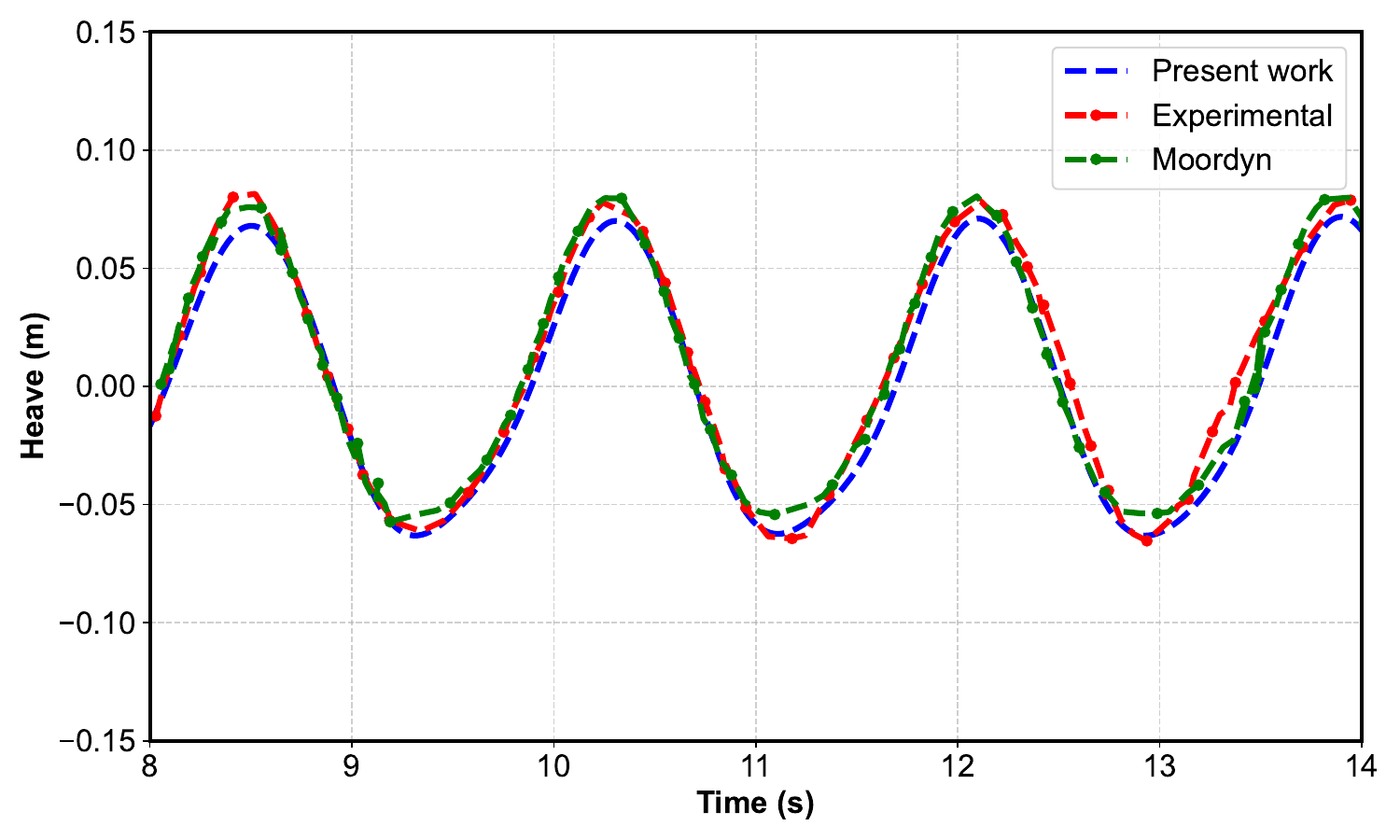}
         \caption{H15T18-Heave}
         \label{fig:H15T18-Heave}
     \end{subfigure}
     \hfill
     \begin{subfigure}[b]{0.45\textwidth}
         \centering
         \includegraphics[width=\textwidth]{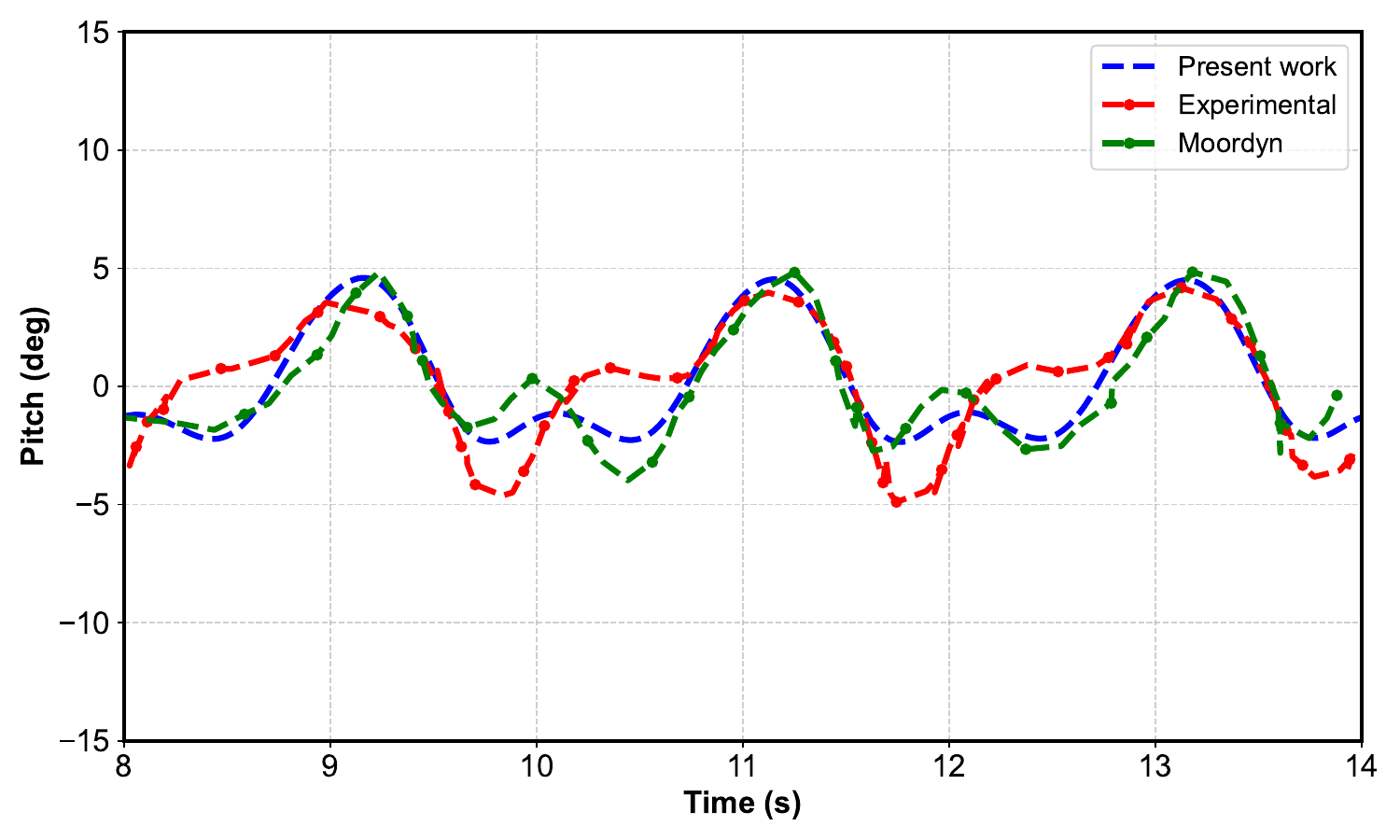}
         \caption{H12T20-Pitch}
         \label{fig:H12T20-Pitch}
     \end{subfigure}
     \hfill
     \begin{subfigure}[b]{0.45\textwidth}
         \centering
         \includegraphics[width=\textwidth]{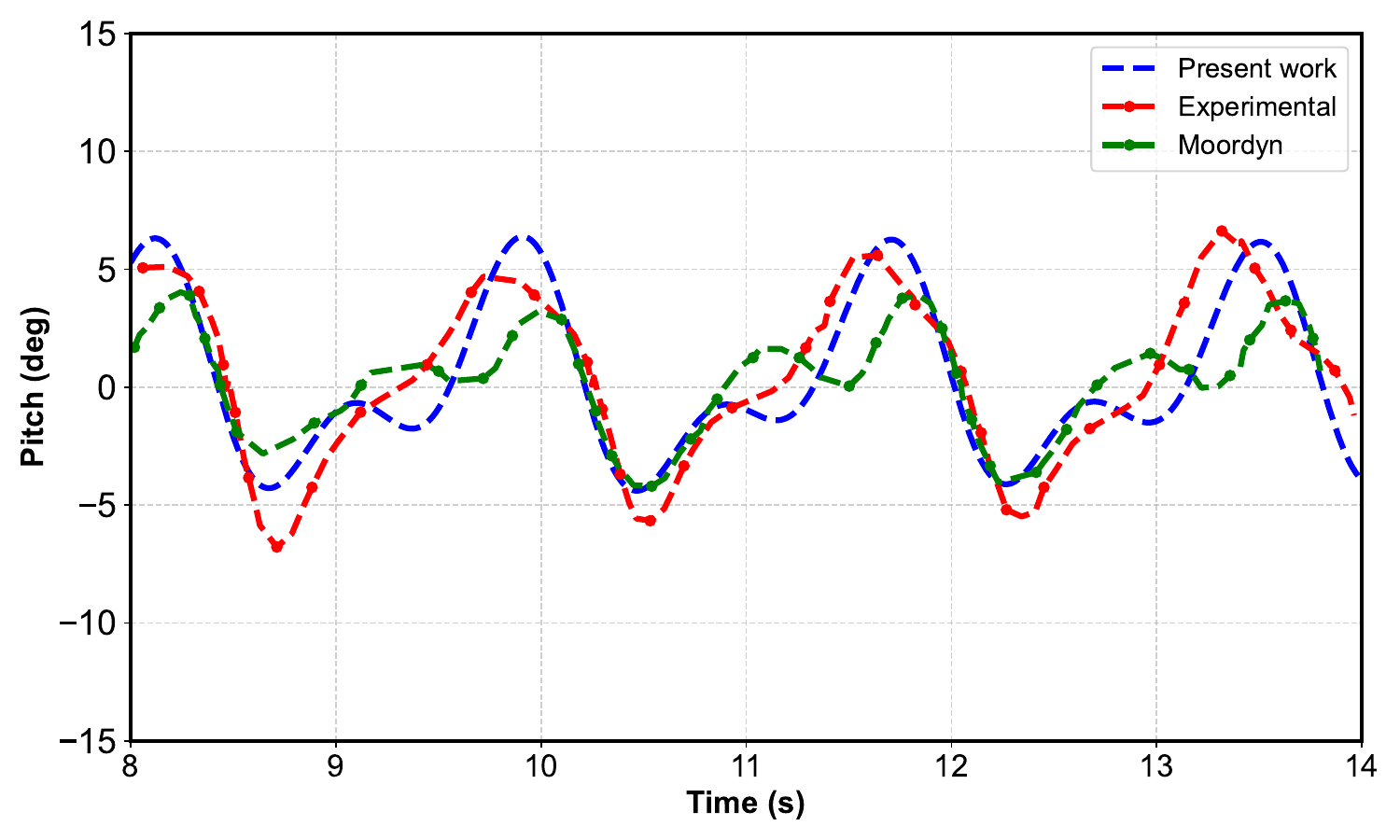}
         \caption{H15T18-Pitch}
         \label{fig:H15T18-Pitch}
     \end{subfigure}
     \hfill
        \caption{Box Motion}
        \label{fig:boxmotions}
\end{figure}

Figure \ref{fig:boxmotions} compares the surge, heave, and pitch responses of the floating box for cases H12T20 and H15T18. It compares results from the present model against MoorDyn \citep{chen_cfd_2022}, and experimental data. The numerical time series were shifted in time so that their first peak coincides with the experimental record, thereby removing the initial phase lag and making the comparison clearer. A broad consistency exists between the numerical and experimental results for heave and surge motion. Although surge motion is overpredicted in both wave cases, the steeper wave condition (H15T18) shows a smaller phase shift relative to the experiments for both numerical models. The interpretation also applies to heave motions (Figures \ref{fig:H12T20-Heave} and \ref{fig:H15T18-Heave}). 

\noindent Amplitude is computed as the average of multiple peak-to-trough values in the time interval between 8 and 16 seconds, corresponding to approximately 4.4 wave periods for case H15T18 and four wave periods for case H12T20. This provides a consistent estimate of the characteristic oscillation magnitude. The percentage error is reported relative to the experimental amplitude. For H12T20, both the present model and MoorDyn overpredict the surge amplitude (Figure \ref{fig:H12T20-Surge}), with the present model deviating by 9.09\% and MoorDyn by 27.40\%. In heave, the present model matches the experimental data more closely than MoorDyn (Figure \ref{fig:H12T20-Heave}), with a deviation of 3\%, compared with 9.08\% for MoorDyn. For pitch, both numerical models slightly underestimate the amplitude; the present model is closer to the experimental value (Figure \ref{fig:H12T20-Pitch}), with a 2.55\% error compared to {3.93\%} from the MoorDyn. In the steeper wave case (H15T18), the present model performs better in all surge, heave and pitch predictions. For surge, the present model differs from the experiments by 8.6\%, while MoorDyn shows a larger 13.04\% error.
 In heave, the present model has a deviation of 6.49\%, compared to 7.14\% for MoorDyn (Figure \ref{fig:H15T18-Heave}). For pitch, the present model shows 4.32\% error, while MoorDyn underpredicts the amplitude with a larger deviation of 38.60\% (Figure \ref{fig:H15T18-Pitch}). These findings indicate that the present model captures pitch responses more consistently under steeper wave conditions (H15T18) than MoorDyn. Stronger nonlinear coupling and closer alignment between the wave period and the system’s natural pitch frequency likely contribute to the lower pitch prediction error in H15T18 compared to H12T20. Overall, the present model delivers smaller amplitude errors than MoorDyn across both wave cases, demonstrating improved predictive capability. Nevertheless, further investigation would help confirm and refine these observations.


\noindent Discrepancies between the numerical and experimental motion responses, particularly in pitch and surge, can be attributed to several factors. 
One likely contributor is the wooden plate with reflective markers mounted on the front face of the floating box in the experimental setup {\citep{wu_experimental_2019}}, which was not included in the numerical model. This addition introduces asymmetry that could influence the added mass and hydrodynamic loading, especially in pitch. Moreover, splashing observed on this plate during the experiments {\citep{dominguez_sph_2019, wu_experimental_2019}} may further affect the pitch response and tracking accuracy. In the present work, the plate has been explicitly accounted for in the reported mass properties (centre of gravity and moments of inertia) to represent its static contribution, but no additional empirical damping was introduced, since such effects cannot be quantified without dedicated decay tests. This choice preserves the physical transparency of the numerical model, avoiding case specific calibration.
The numerical model also assumes ideal second-order Stokes waves, while in the physical setup, the wave maker introduces startup transients and small irregularities, potentially contributing to phase shifts or amplitude deviations. Additional simplifications, such as the assumption of quiescent flow around the mooring lines and approximated hydrodynamic coefficients, may also influence the simulated responses. Similar discrepancies have been reported in previous numerical studies using this benchmark case \citep{chen_cfd_2022, dominguez_sph_2019}.

\begin{table}[h]
    \centering
    \caption{Comparison of box motion amplitudes and percentage differences from experimental values}
    \begin{tabular}{lcccccc}
        \hline
        & \multicolumn{2}{c}{Experimental} & \multicolumn{2}{c}{MoorDyn \citep{chen_cfd_2022}} & \multicolumn{2}{c}{Present Model} \\
        & Amplitude & Diff. & Amplitude & Diff. (\%) & Amplitude & Diff. (\%) \\
        \hline
        \textbf{H12T20} & & & & & & \\
        Surge (m) & 0.22 & - & 0.28 & 27.40 & 0.24 & 9.09 \\
        Heave (m) & 0.11 & - & 0.12 & 9.08 & 0.1133 & 3  \\
        Pitch (°) & 9.0  & - & 8.71   & 3.93 & 8.77  & 2.55 \\
        \hline
        \textbf{H15T18} & & & & & & \\
        Surge (m) & 0.23 & - & 0.26 & 13.04 & 0.25 & 8.6 \\
        Heave (m) & 0.14 & - & 0.13 &  7.14 & 0.131 & 6.49 \\
        Pitch (°) & 13.40  & - & 8.23   & 38.60 & 13.98  & 4.32 \\
        \hline
    \end{tabular}
    \label{tab:box_motion_amplitudes}
\end{table}
\clearpage
\subsection{Mooring Line tensions}
\begin{figure}[!ht]
     \centering
     \begin{subfigure}[b]{0.45\textwidth}
         \centering
         \includegraphics[width=\textwidth]{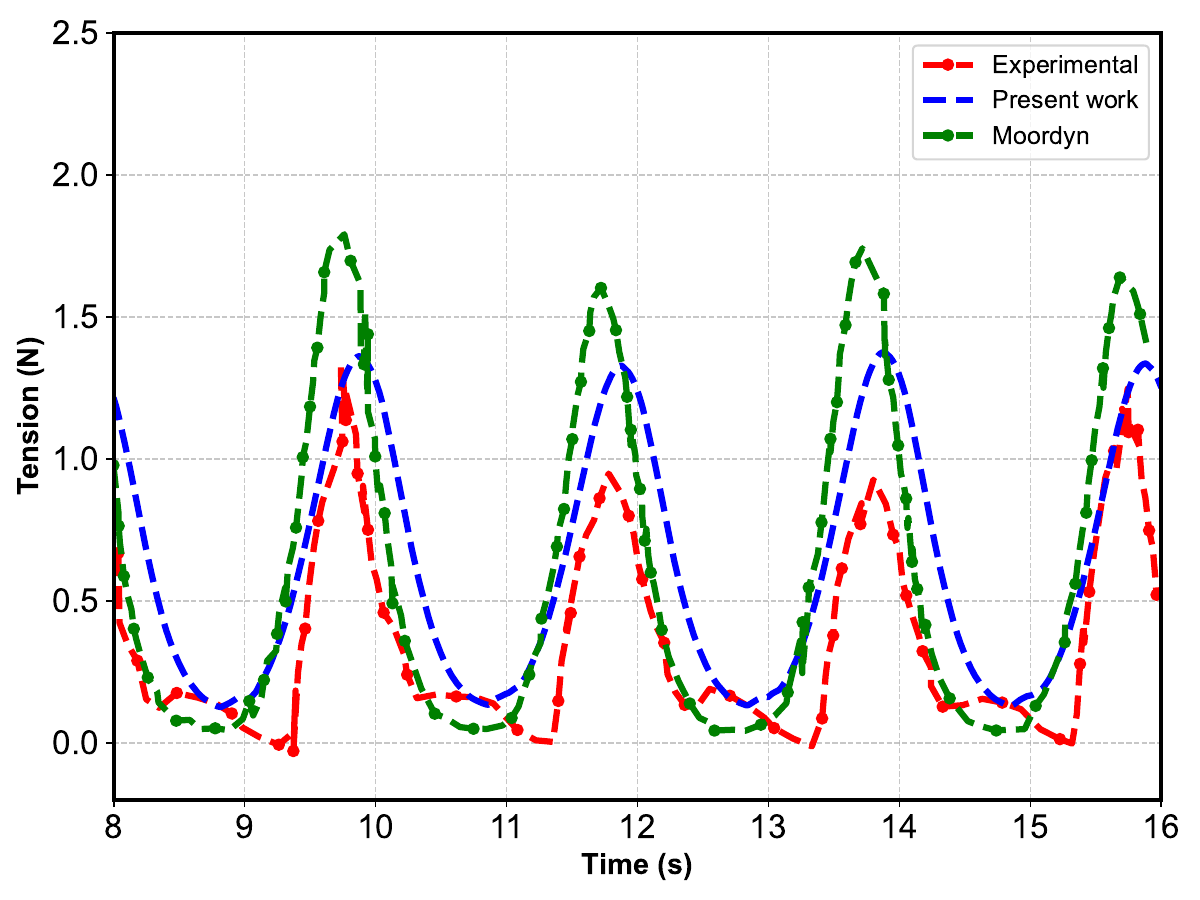}
         \caption{line1 Tension - H12T20}
         \label{fig:line1TensionH12T20}
     \end{subfigure}
     \hfill
     \begin{subfigure}[b]{0.45\textwidth}
         \centering
         \includegraphics[width=\textwidth]{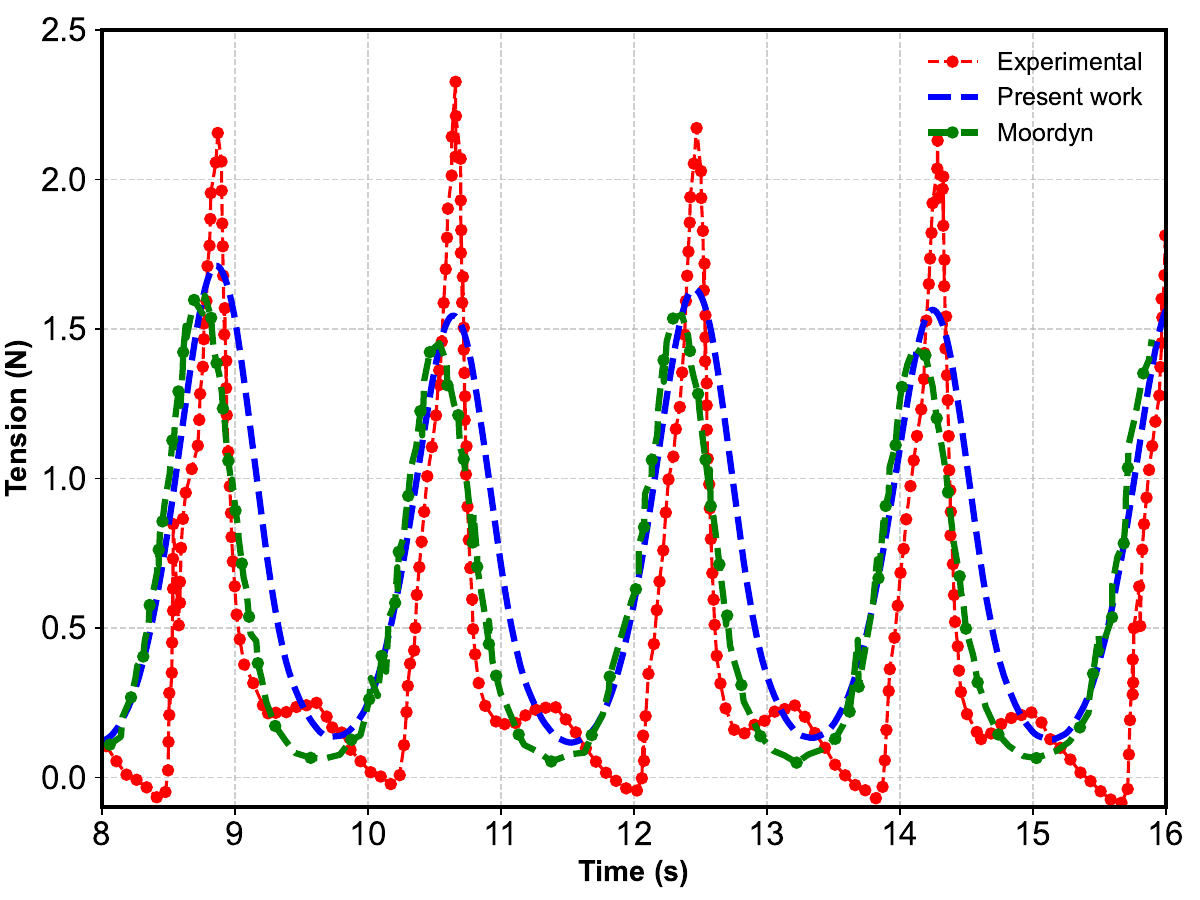}
         \caption{line1 Tension - H15T18}
         \label{fig:line1TensionH15T18}
     \end{subfigure}
     \hfill
     \begin{subfigure}[b]{0.45\textwidth}
         \centering
         \includegraphics[width=\textwidth]{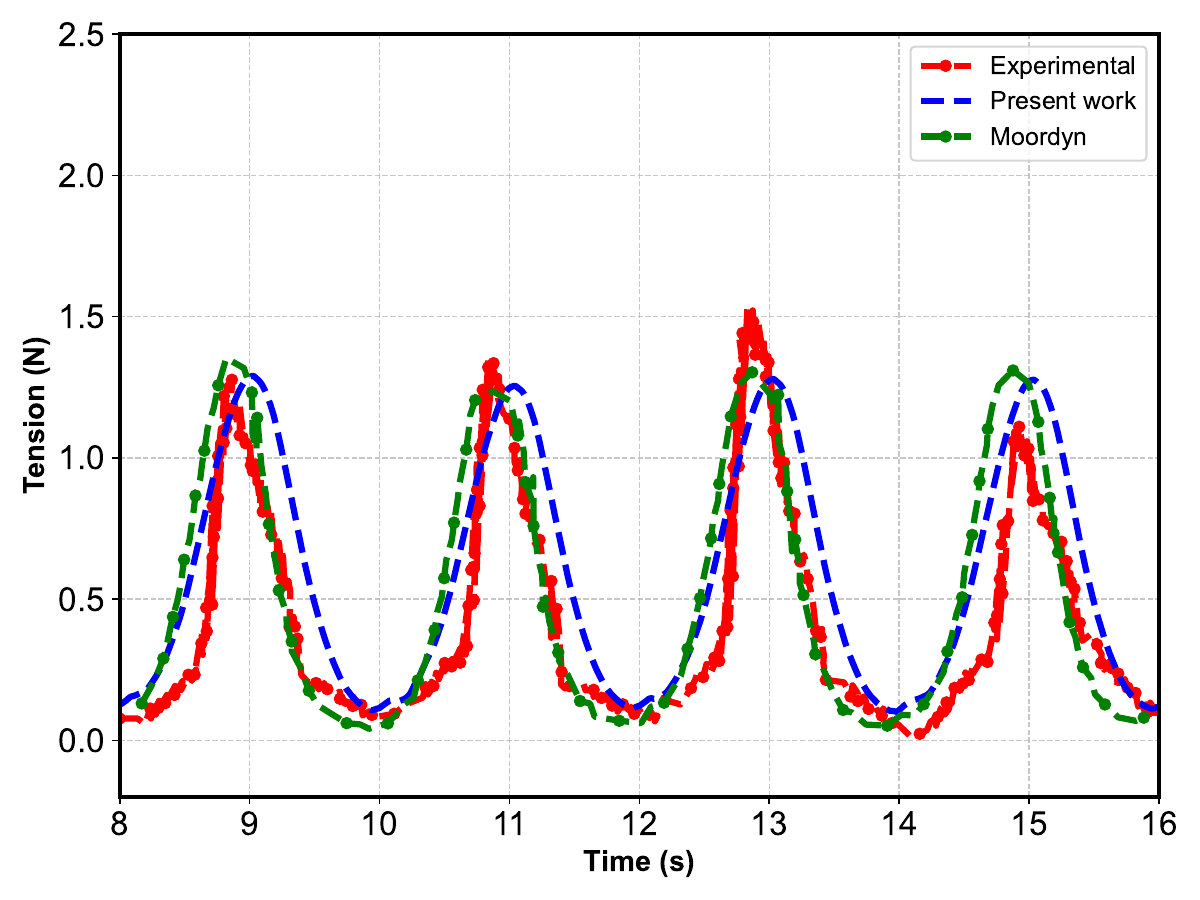}
         \caption{line3 Tension - H12T20}
         \label{fig:line3TensionH12T20}
     \end{subfigure}
     \hfill
     \begin{subfigure}[b]{0.45\textwidth}
         \centering
         \includegraphics[width=\textwidth]{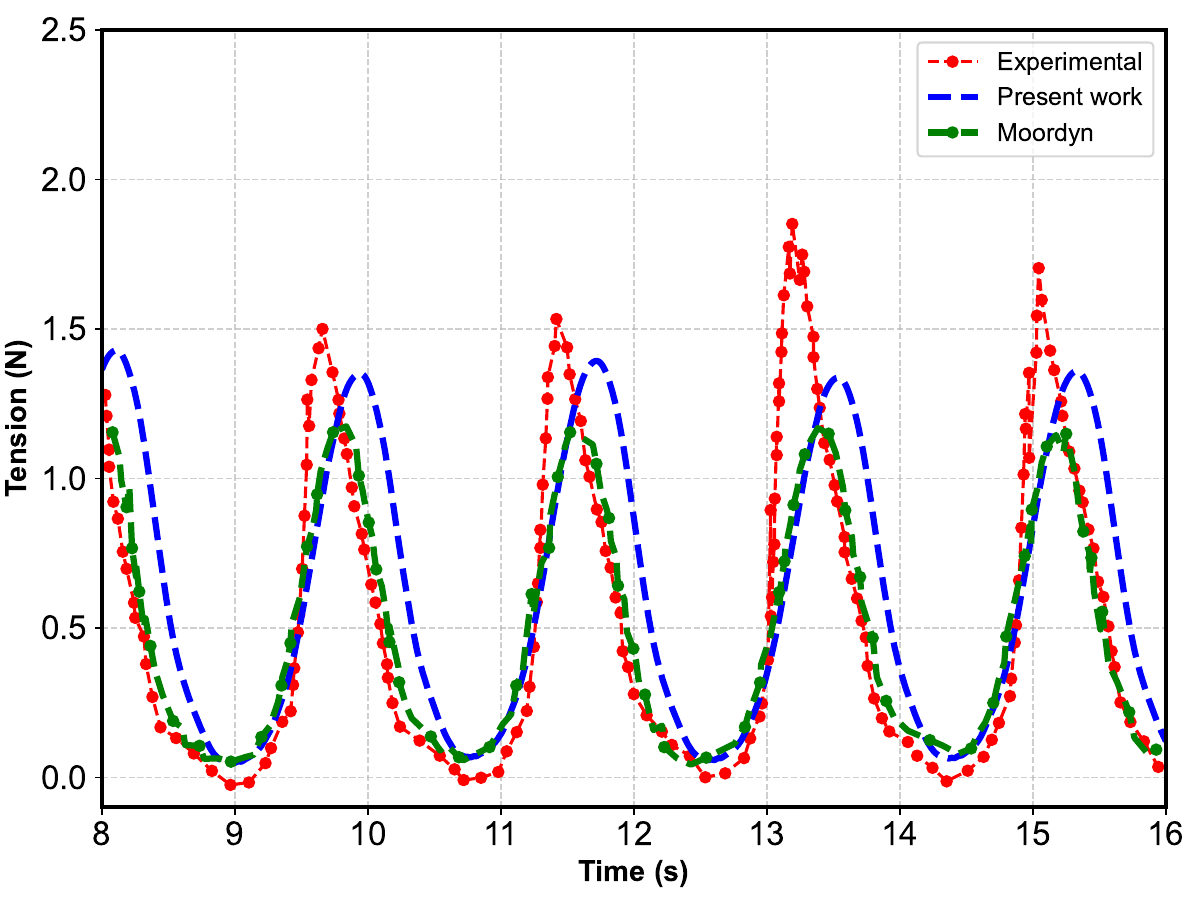}
         \caption{line3 Tension - H15T18}
         \label{fig:line3TensionH15T18}
     \end{subfigure}
     \hfill
        \caption{Mooring Line Tension}
        \label{fig:mooringLineTension}
\end{figure}

The time series of mooring line responses for the two wave cases are shown in Figure \ref{fig:mooringLineTension}. Since the system is a four-point symmetric mooring system, only the tensions for line 1 and line 3 at the anchors are shown (see Figure \ref{fig:NomDom}). The simulations and experiments show a close agreement in the shape of the tension time history. The waveward cables (line 1) experience larger tension offsets compared to the leeward cables (line 3) in both cases. This is expected, as waveward cables are more influenced by mean wave drift forces, leading to occasional snap loads.

\noindent A quantitative comparison of mooring line tension amplitudes between the present model, MoorDyn, and experimental results is provided in Table \ref{tab:mooring_tension}.
Tension amplitude is computed as the average of multiple peak-to-trough values in the time interval between 8 and 16 seconds, covering approximately 4.4 wave periods for case H15T18 and four wave periods for case H12T20. This approach provides a reliable estimate of the dominant dynamic response. The percentage error is reported relative to the experimental amplitude. For H12T20, the present model slightly underpredicts line 1 tension (Figure \ref{fig:line1TensionH12T20}) by 3.18\% compared to the experimental value. MoorDyn, on the other hand, exhibits a significantly higher deviation of 28.17\%. In line 3 (Figure \ref{fig:line3TensionH12T20}), both the present model and MoorDyn show reasonable agreement with experiments, though MoorDyn tends to underestimate the amplitude by 13.62\% as compared to the present model, which shows a deviation of 7.43\%. In the steeper wave case (H15T18), the discrepancies become more pronounced. The present model underestimates the line 1 amplitude (Figure \ref{fig:line1TensionH15T18}) by 23.62\%, while MoorDyn shows a slightly larger deviation of 26.94\%. In line 3 (Figure \ref{fig:line3TensionH15T18}), the present model exhibits a 25.25\% error, whereas MoorDyn deviates by 39.68\%, again showing a noticeably larger mismatch with the experimental data. These results suggest that the present model may offer an improvement across different wave conditions. This improvement may be mainly due to the Simo-Reissner beam formulation providing a more detailed representation of the internal forces under dynamic deformation, which may contribute to slightly closer agreement with the experimental measurements in mooring tensions. It is worth noting that no snap load events were detected in the tested cases. The selected catenary configuration and the range of wave heights (0.12–0.15 m) ensured that the mooring lines remained continuously under tension throughout the simulations.

\begin{table}[h]
    \centering
    \caption{Comparison of mooring line tension amplitudes}
    \begin{tabular}{lcccccc}
        \hline
        & \multicolumn{2}{c}{Line 1 (waveward side)} & \multicolumn{2}{c}{Line 3 (leeward side)} \\
        & Amplitude (N) & Diff. (\%) & Amplitude (N) & Diff. (\%) \\
        \hline
        \textbf{H12T20} & & & & \\
        Experimental & 1.3639 & - & 1.519 & - \\
        MoorDyn \citep{chen_cfd_2022} & 1.7481 & 28.17 & 1.3121 & 13.62 \\
        Present Model & 1.3205 & 3.18 & 1.406 & 7.43 \\
        \hline
        \textbf{H15T18} & & & & \\
        Experimental & 2.1385 & - & 1.8764 & - \\
        MoorDyn \citep{chen_cfd_2022} & 1.5623 & 26.94 & 1.1317 & 39.6876 \\
        Present Model & 1.633 & 23.624 & 1.4025 & 25.25 \\
        \hline
    \end{tabular}
    \label{tab:mooring_tension}
\end{table}

\noindent Several factors may contribute to the discrepancies observed in the predicted mooring line tensions. While mean wave drift forces likely account for part of the offset, additional modelling assumptions could also influence the results. These include the simplified representation of the mooring line connections, particularly the absence of flexibility or deformation at the chain links, which were present in the experimental setup but are not captured in the current model. Assumed values for axial stiffness and hydrodynamic coefficients, such as drag and added mass, introduce further uncertainty. Moreover, the assumption of a quiescent fluid surrounding the cables neglects the influence of local flow disturbances caused by the floating body and incoming waves. This may be particularly important in the steeper wave case (H15T18), where more energetic interactions occur near the fairlead. Similar limitations have been discussed in prior numerical studies \citep{chen_cfd_2022, dominguez_sph_2019}. While these explanations are plausible, isolating their effects remains difficult without dedicated experiments or higher-fidelity coupling. This underscores the inherent challenges of reproducing detailed experimental conditions in numerical simulations.

\begin{figure}[!ht]
     \centering
     \begin{subfigure}[b]{0.45\textwidth}
         \centering
         \includegraphics[width=\textwidth]{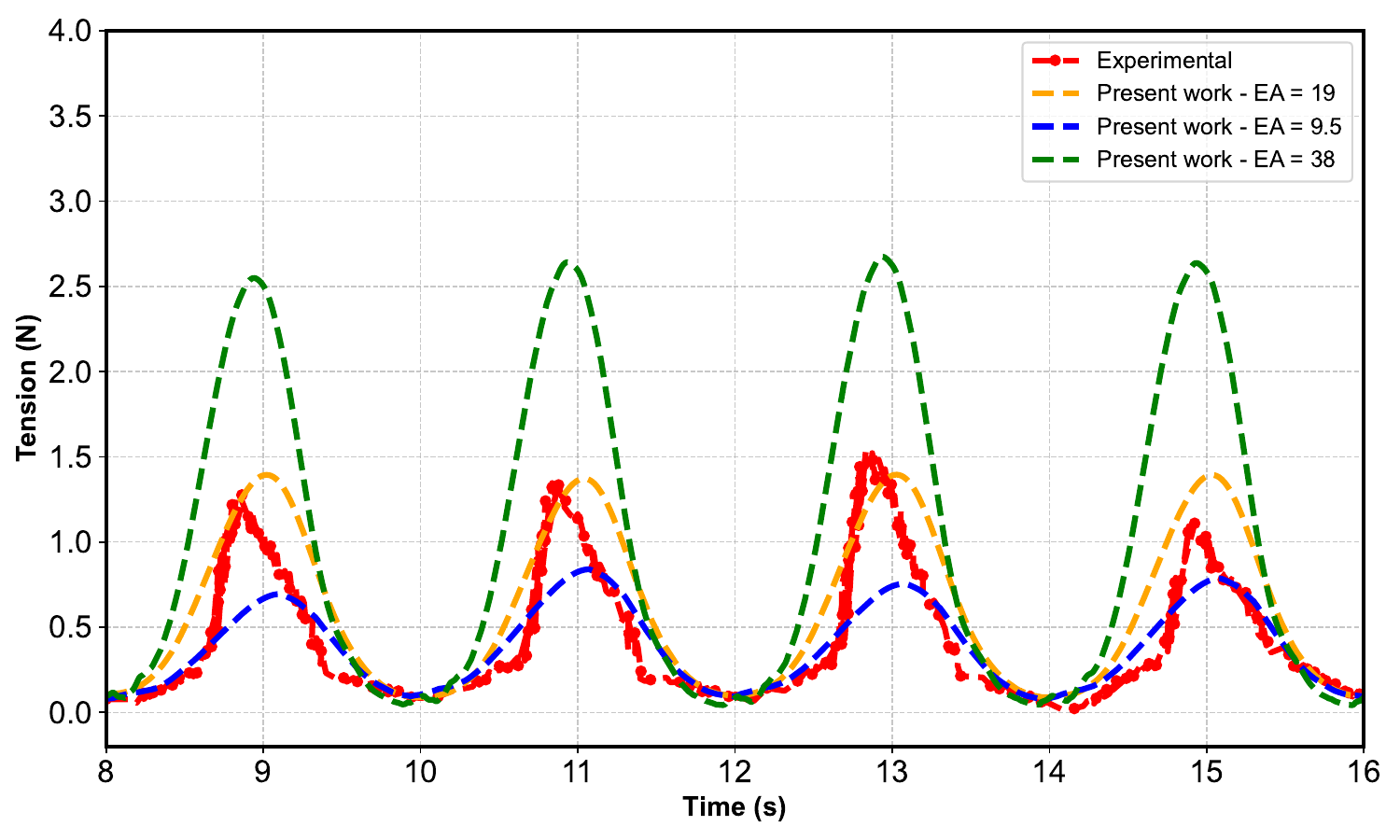}
         \caption{line1 Tension - H12T20}
         \label{fig:line1TensionH12T20}
     \end{subfigure}
     \hfill
     \begin{subfigure}[b]{0.45\textwidth}
         \centering
         \includegraphics[width=\textwidth]{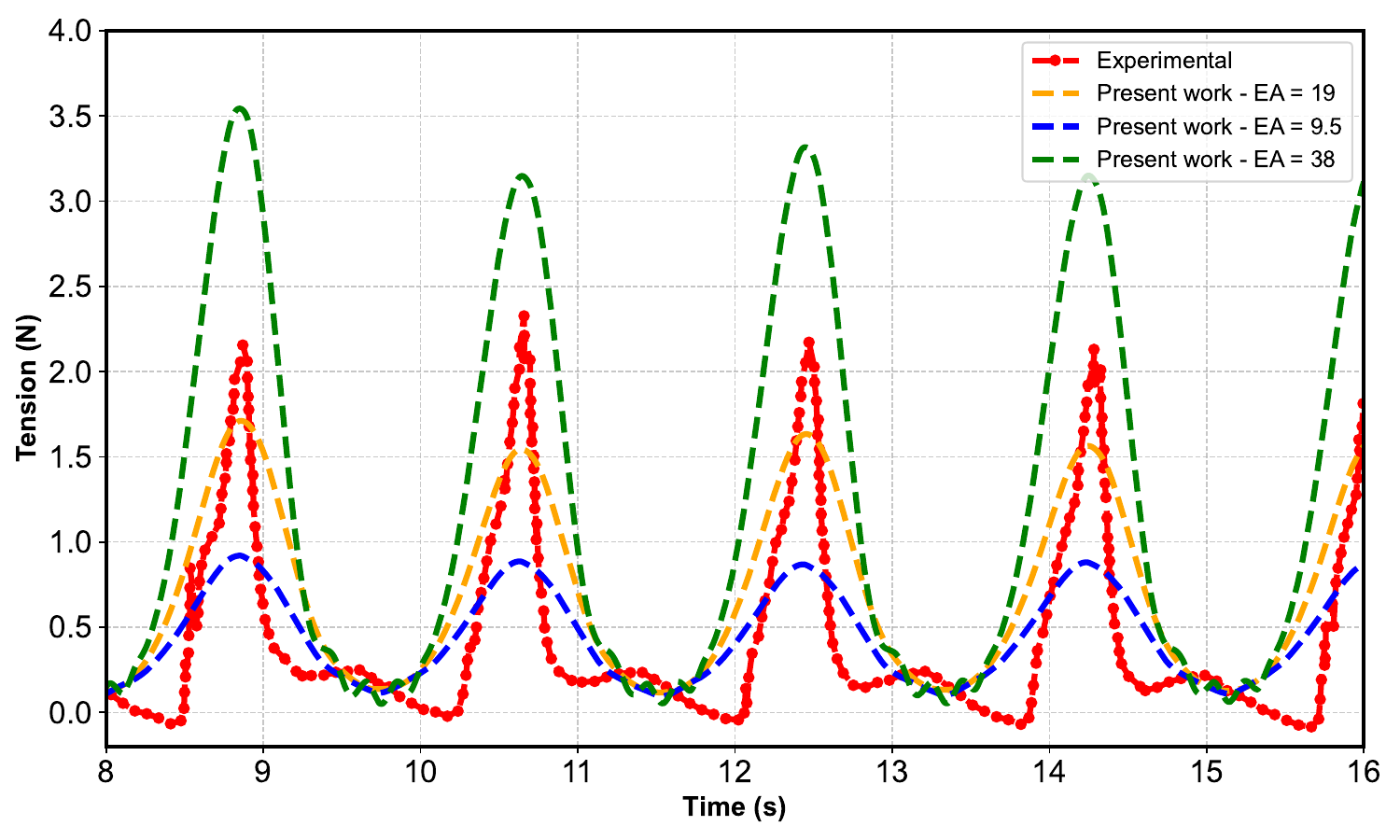}
         \caption{line1 Tension - H15T18}
         \label{fig:line1TensionH15T18}
     \end{subfigure}
     \hfill
\caption{Influence of axial stiffness($EA$) on mooring line tension}
        \label{fig:EAStudy}
\end{figure}

\noindent To investigate the sensitivity of mooring line tension to axial stiffness, the parameter $EA$ was varied from the baseline value of $19$ N to $9.5$ N and $38$ N for both H12T20 and H15T18 cases. The comparison in Figure \ref{fig:EAStudy} shows that in both wave conditions, reducing the stiffness to $9.5$ N results in a visibly damped response with lower peak tensions and broader force profiles, while increasing it to $38$ N produces sharper and more pronounced peaks. The baseline value of $EA = 19$ N consistently provides the best agreement with experimental results across both cases, especially in terms of peak magnitude and phase alignment. The H15T18 case, featuring a steeper wave profile and more pronounced wave loads, demonstrates greater sensitivity to $EA$ changes. Here, as also reported by \citet{chen_cfd_2022}, increasing $EA$ does help in partially recovering the under-predicted tension peaks. However, the current results show that overly stiff lines ($EA = 38$ N) tend to overshoot the experimental tension, introducing unrealistic spikes. These results show that axial stiffness is really a property of the whole mooring system, not just the cable material. Chain-link flexibility, contact with the seabed or fairlead, and hydrodynamic damping all change the stiffness that the line presents in practice. Calibrating $EA$ against integrated dynamic responses rather than static material tests is therefore crucial for reliable modelling of mooring line behaviour.

\section{Conclusion}\label{sec:conclusion}
This paper presents a novel unified finite-volume-based discretisation framework for modelling mooring cable dynamics, integrating multiphase CFD flow simulations and a six-degree-of-freedom motion solver. The framework combines OpenFOAM’s native six-degree-of-freedom motion solver with a finite-volume-based Simo-Reissner beam solver. This integration enables detailed analysis of the moored motion of a floating rigid box, accounting for the bending, shear, and torsional deformations of the mooring line. To validate the proposed approach, numerical results of the current model for a floating rigid box moored with four catenary lines are compared to both experimental and numerical data to assess the framework’s ability to capture complex interactions between the structure and its environment. The following observations are made from the numerical analyses:

\begin{itemize}
\item Aligning the first wave peak with the experimental data removed an initial phase lag in the surface elevation time histories and simplified direct comparisons, demonstrating that the numerical model successfully reproduces the dominant wave frequency and closely matches the amplitudes reported in \citet{wu_experimental_2019} (see Figure \ref{fig:fftofWaves}). However, the coarser mesh region near the outlet exhibits a noticeable offset in wave amplitude (Figures~\ref{fig:wg612} and \ref{fig:wg615}). This confirms the model’s ability to resolve local wave behaviour near the structure but highlights limitations in far-field accuracy.
\item Mooring line tensions on the leeward side show a strong correlation with the experimental data, while both numerical models deviate from the experimental values on the waveward side. However, the present model consistently provides a closer match to experimental values compared to MoorDyn, as seen in line 1 of H12T20 (3.18\% vs 28.17\%), line 3 of H12T20 (7.43\% vs 13.62\%), and line 3 of H15T18 (25.25\% vs 39.68\%), where MoorDyn exhibits substantially larger deviations. The improved accuracy of the present model suggests a better representation of mooring line dynamics under wave-induced forces. The remaining discrepancies are primarily attributed to the influence of mean wave drift forces, which are not fully captured in the numerical models, as well as additional assumptions regarding axial stiffness, fairlead and anchor connections, and the quiescent fluid assumption. These factors collectively contribute to the observed offset. These results indicate that the present model enhances the reliability of mooring load predictions.
\item Surge and heave motions show good agreement with experimental data, with both the present model and MoorDyn slightly overpredicting these motions. In surge, the present model shows smaller deviations than MoorDyn in both cases, with 9.09 \% compared to 27.40\% for H12T20 and 8.6\% compared to 13.04\% for H15T18. In heave, both models closely match the experimental amplitude for H12T20, with the present model deviating by 3\% and MoorDyn by 9.08\%. For H15T18, both models slightly underpredict the heave motion, with the present model showing a smaller error of 6.49\% compared to 7.14\% for MoorDyn. (see Figures \ref{fig:H12T20-Surge}, \ref{fig:H12T20-Heave}, \ref{fig:H12T20-Heave}, \ref{fig:H15T18-Heave}). This supports the validity of the present model in predicting primary translational motions of the floating body.
\item The largest discrepancies arise in pitch motion (see Figures \ref{fig:H12T20-Pitch} and \ref{fig:H15T18-Pitch}), where the present model matches the observed amplitude more consistently under the steeper wave condition (H15T18). This improved performance may be linked to the closer alignment between the wave period and the system’s natural pitch frequency, which tends to amplify any modelling inaccuracies or strengths. By contrast, MoorDyn underpredicts the amplitude in these conditions. Some of the remaining errors appear to stem from experimental uncertainties, such as a wooden plate on the box’s front face, which is omitted in the numerical model. Similar discrepancies have also been noted in prior studies \citep{chen_cfd_2022, dominguez_sph_2019}, highlighting ongoing challenges in numerically simulating moored floating structures. This indicates that the present model better captures rotational dynamics under nonlinear wave-structure interaction.
\item The proposed model demonstrates efficient scaling (see Figure \ref{fig:strongScaling}), with a notable decrease in execution time as additional CPU cores are utilised. However, beyond 256 CPU cores, the improvement in performance becomes less significant. This suggests the model is suitable for large-scale simulations, but benefits wane as the number of cells per CPU core becomes small.
\end{itemize}

\noindent Overall, the results highlight both the strengths and the current limitations of the modelling framework. In conclusion, this study reinforces that improving the representation of wave–mooring interaction is essential for advancing CFD-coupled simulations. While current frameworks perform well under simplified conditions, their ability to replicate complex experimental scenarios remains limited. Bridging this gap will require both refined coupling strategies and targeted experimental validation.

\backmatter

\bmhead{Acknowledgments}
The authors would like to acknowledge Research Ireland, NexSys project 21/SPP/3756, for funding. Vikram Pakrashi would like to acknowledge SEAI RDD/966 FlowDyn and Research Ireland MaREI RC2302-2 projects. Philip Cardiff gratefully acknowledges financial support from the Irish Research Council (IRC) through the Laureate programme, grant number IRCLA/ 2017/45 and the European Research Council (ERC) under the European Union’s Horizon 2020 research and innovation programme (Grant Agreement No. 101088740).
Additionally, the authors wish to acknowledge the Irish Centre for High-End Computing (ICHEC) for the provision of computational facilities and support (www.ichec.ie), and part of this work has been carried out using the UCD ResearchIT Sonic cluster which was funded by UCD IT Services and the UCD Research Office.
\bmhead{Availability of Data and Code}
The numerical model cases and source code used in this study are available at:  
\url{https://github.com/solids4foam/moorFV}

\bmhead{Declaration of generative AI and AI-assisted technologies in the writing process}
During the preparation of this work, the author(s) used Grammarly and ChatGPT to improve the grammar and readability of the manuscript. After using this tool/service, the author(s) reviewed and edited the content as needed and take(s) full responsibility for the content of the published article.
\appendix
\section{Grid convergence study}\label{sec:gridConvergence}
Three different mesh resolutions were selected for the grid convergence study. For this study, the H12T20 case with nominal wave height ($H=12 cm$) has been used. The medium-resolution mesh has a cell size of $0.01 \times 0.011 \times 0.01$ in the refined region of the domain. Outside this region, the mesh is gradually stretched using a grading ratio of $1.07$ (see Figure~\ref{fig:meshOverview}). The fine mesh has smaller cells than the medium mesh by a factor of $\sqrt{2}$, while the coarse mesh has larger cells by the same factor.
\begin{figure}[!ht]
     \centering
     \begin{subfigure}[b]{0.45\textwidth}
         \centering
         \includegraphics[width=\textwidth]{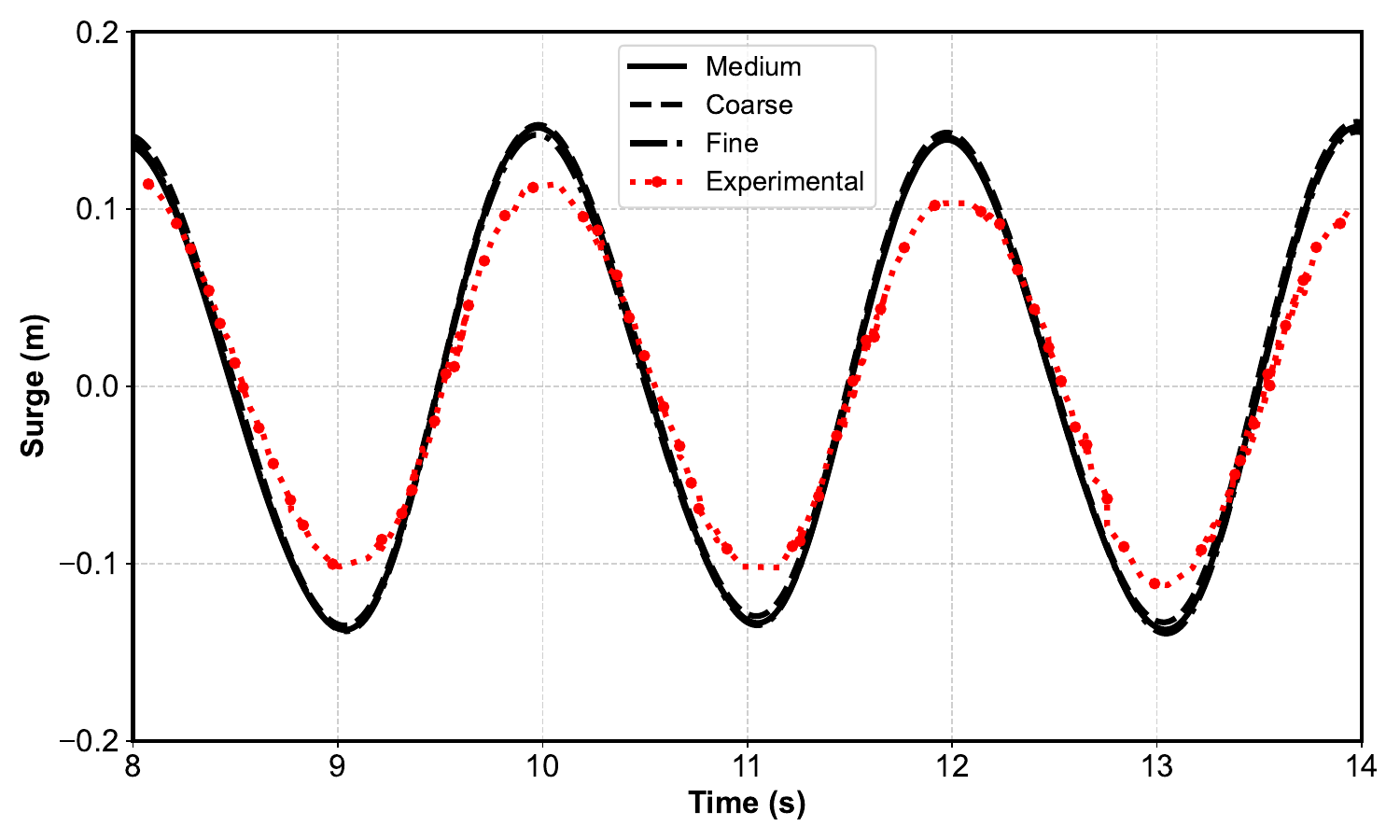}
         \caption{Surge}
         \label{fig:gcs-surge}
     \end{subfigure}
     \hfill
     \begin{subfigure}[b]{0.45\textwidth}
         \centering
         \includegraphics[width=\textwidth]{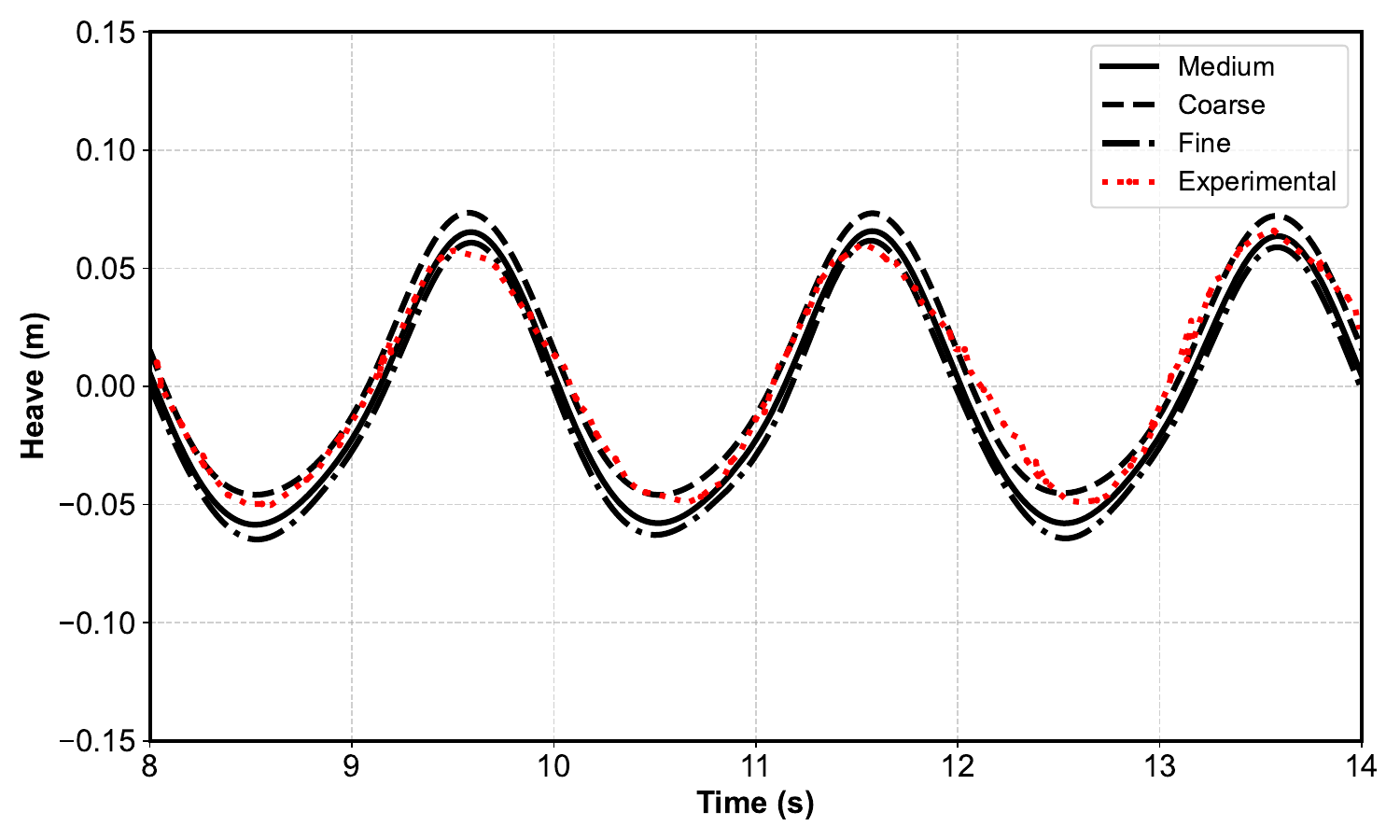}
         \caption{Heave}
         \label{fig:gcs-heave}
     \end{subfigure}
     \hfill
     \begin{subfigure}[b]{0.45\textwidth}
         \centering
         \includegraphics[width=\textwidth]{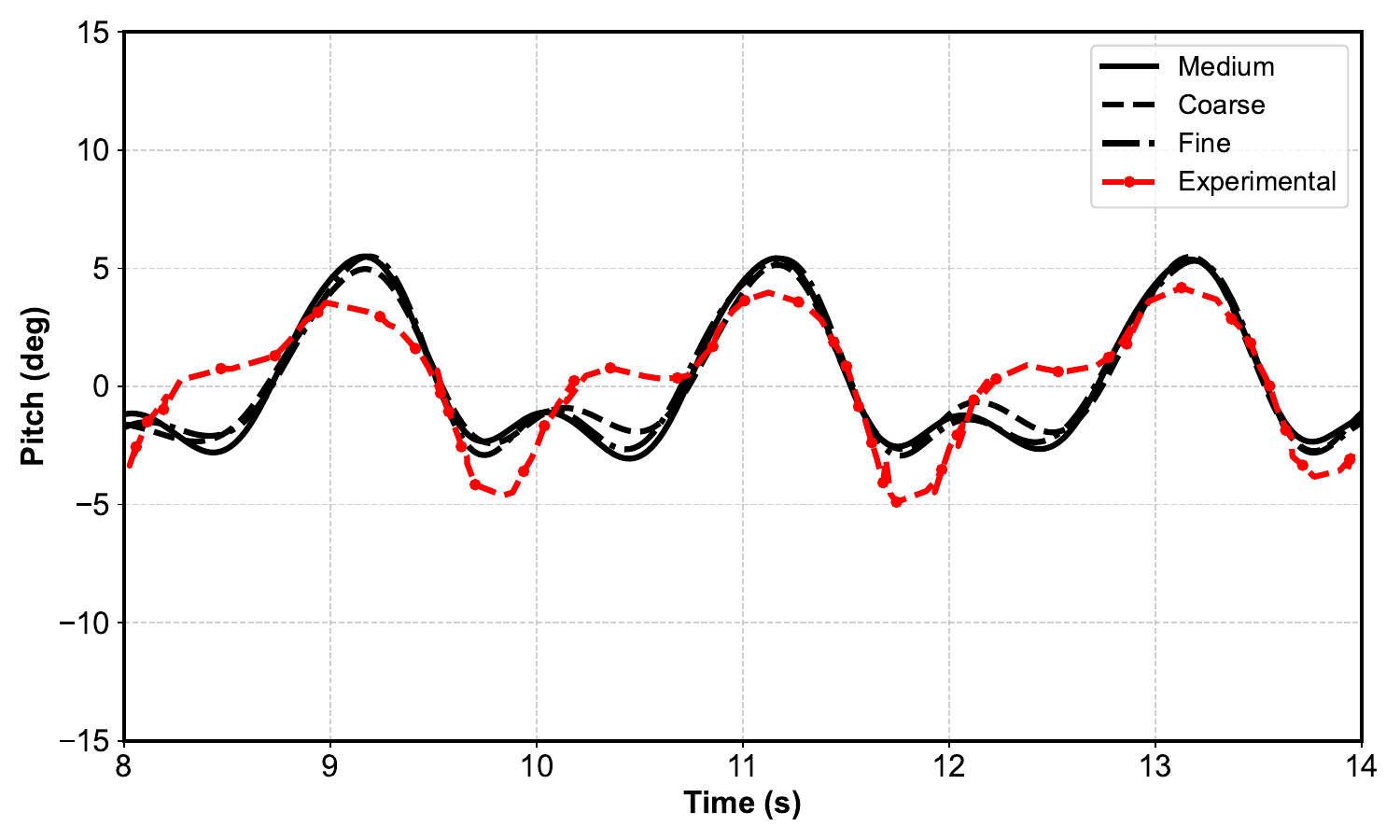}
         \caption{Pitch}
         \label{fig:gcs-pitch}
     \end{subfigure}
     \hfill
     \caption{Comparing box motion using different grid resolutions}
     \label{fig:gridStudy}
\end{figure}

\noindent The motion of the box for three different grid sizes is illustrated in Figure \ref{fig:gridStudy}. The surge motion, which is directly influenced by the restraining force of the mooring lines, remains nearly identical across all mesh resolutions. The results for heave and pitch differ across the coarse, medium, and fine meshes. In both cases, the medium-resolution results fall between those of the coarse and fine meshes, and are noticeably closer to the fine mesh results. This indicates that while there is some mesh sensitivity, the medium mesh provides a reasonable approximation of the converged solution.
\section{Time step sensitivity study}\label{sec:timeStepStudy}
In addition to the grid convergence analysis presented in Appendix~\ref{sec:gridConvergence}, a time step sensitivity analysis was carried out to assess the influence of the Courant number on the accuracy of the floating body responses. Three cases were considered, corresponding to maximum Courant numbers of 0.25, 0.50, and 0.75. Since the solver employs an adaptive time step based on the Courant condition, these values directly control the maximum allowable time step size. The body motions obtained from these simulations were compared against the experimental measurements for the H12T20 case.

\begin{figure}[!h]
    \centering
    \begin{subfigure}[b]{0.45\textwidth}
        \includegraphics[width=\textwidth]{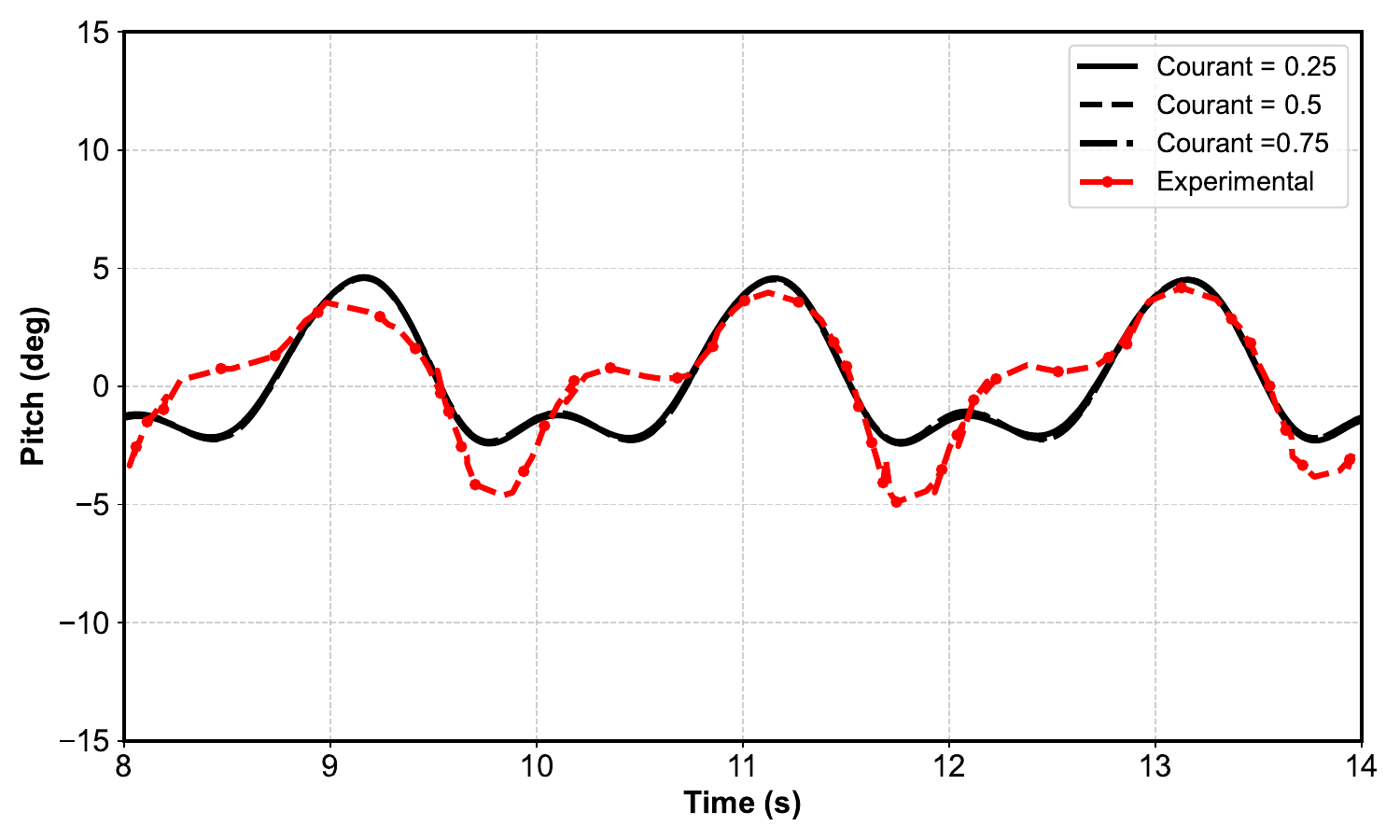}
        \caption{Pitch}
        \label{fig:tsaPitch}
    \end{subfigure}
    \hfill
    \begin{subfigure}[b]{0.45\textwidth}
        \includegraphics[width=\textwidth]{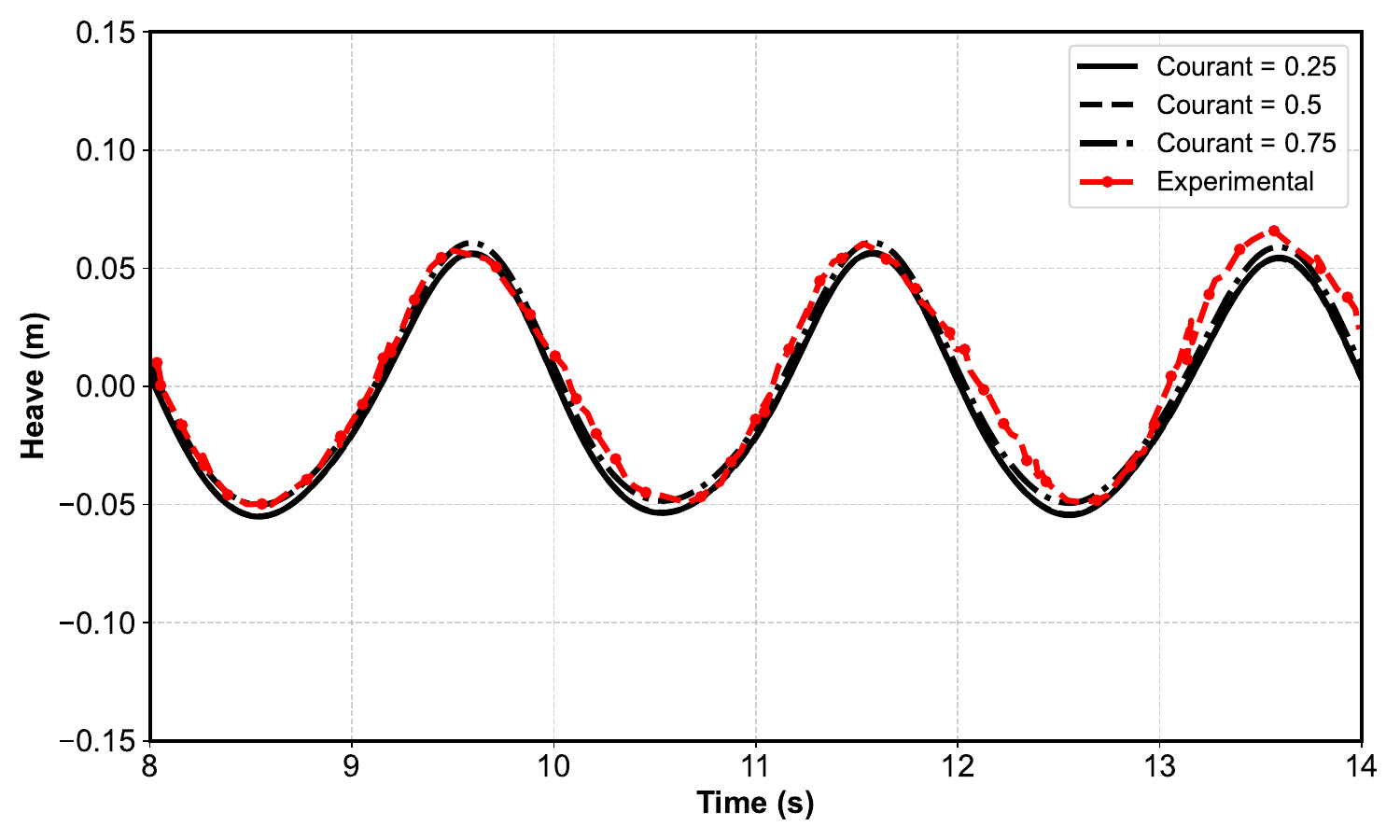}
        \caption{Heave}
        \label{fig:tsaHeave}
    \end{subfigure}
    \hfill
    \begin{subfigure}[b]{0.45\textwidth}
        \includegraphics[width=\textwidth]{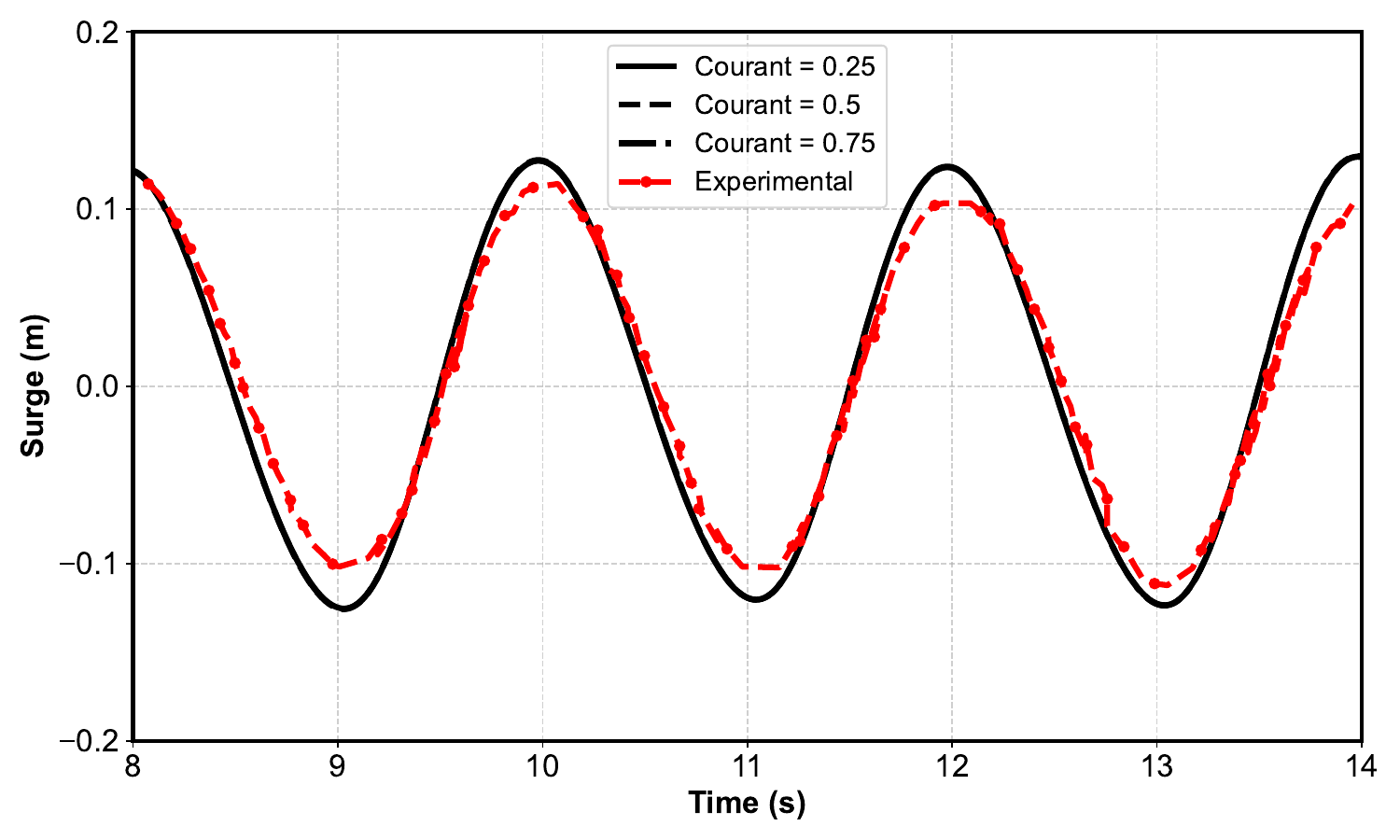}
        \caption{Surge}
        \label{fig:tsaSurge}
    \end{subfigure}
    \caption{Time step sensitivity analysis (H12T20): simulated and experimental responses in pitch, heave, and surge.}
    \label{fig:tsaCombined}
\end{figure}

\noindent Figures~\ref{fig:tsaPitch}--\ref{fig:tsaSurge} present the pitch, heave, and surge responses for the three Courant number values. Overall, the results indicate that the simulated motions are only marginally affected by the time step size. For pitch (Figure~\ref{fig:tsaPitch}), all three simulations follow the experimental trend with good phase alignment and comparable amplitudes, though the Courant = 0.25 case shows slightly sharper oscillations. For heave (Figure~\ref{fig:tsaHeave}), the three simulations are almost indistinguishable, all showing excellent agreement with the experimental data. Surge motion (Figure~\ref{fig:tsaSurge}) is similarly robust, with the Courant=0.25 and Courant=0.5 cases providing nearly identical predictions and the Courant=0.75 case showing only minor differences. 

\noindent These results suggest that the present model is not particularly sensitive to the chosen Courant number within the tested range. Based on this analysis, a maximum Courant number of 0.50 was selected for the main simulations, as it provides a good balance between numerical accuracy and computational efficiency.

\bibliography{bibiliographNew}

\end{document}